\documentclass[12pt]{amsart}
\usepackage{amscd,amssymb,amsthm,verbatim}
\newcommand{\Aut}{\operatorname{Aut}}

\newcommand{\C}{{\mathbb C}}

\newcommand{\const}{\operatorname{const.}}

\newcommand{\diam}{\operatorname{diam}}
\newcommand{\Diff}{\operatorname{Diff}}

\newcommand{\dvol}{\operatorname{dvol}}

\newcommand{\GL}{\operatorname{GL}}

\newcommand{\HH}{\operatorname{H}}

\newcommand{\Image}{\operatorname{Im}}

\newcommand{\inj}{\operatorname{inj}}

\newcommand{\Ker}{\operatorname{Ker}}

\newcommand{\nil}{\operatorname{nil}}
\newcommand{\Nil}{\operatorname{Nil}}

\newcommand{\OO}{\operatorname{O}}

\newcommand{\PSL}{\operatorname{PSL}}
\newcommand{\pt}{\operatorname{pt}}

\newcommand{\R}{{\mathbb R}}

\newcommand{\Riem}{\operatorname{Riem}}

\newcommand{\SL}{\operatorname{SL}}
\newcommand{\SO}{\operatorname{SO}}
\newcommand{\Sol}{\operatorname{Sol}}

\newcommand{\Tr}{\operatorname{Tr}}

\newcommand{\vol}{\operatorname{vol}}
\newcommand{\Z}{{\mathbb Z}}

\numberwithin{equation}{section}
\theoremstyle{plain}
\newtheorem{definition}[equation]{Definition}

\newtheorem{lemma}[equation]{Lemma}
\newtheorem{theorem}[equation]{Theorem}
\newtheorem{proposition}[equation]{Proposition}
\newtheorem{corollary}[equation]{Corollary}

\theoremstyle{remark}

\newtheorem{remark}[equation]{Remark}
\newtheorem{example}[equation]{Example}
\usepackage{graphicx}
\input{amssym.def}
\input{amssym}
\input{epsf}
\usepackage{a4wide}

\begin{document}

\title[Dimensional Reduction]
{Dimensional Reduction and the Long-Time Behavior of Ricci Flow}

\author{John Lott}
\address{Department of Mathematics\\
University of Michigan\\
Ann Arbor, MI  48109-1043\\
USA} \email{lott@umich.edu}

\thanks{This work was
supported by NSF grant DMS-0604829}
\date{October 25, 2008}

\begin{abstract}
If $g(t)$ is a three-dimensional
Ricci flow solution, with
sectional curvatures that are $O(t^{-1})$ and diameter that is
$O(t^{\frac12})$, then the pullback Ricci flow
solution on the universal cover approaches
a homogeneous expanding soliton.
\end{abstract}

\maketitle

\section{Introduction} \label{section1}

After Perelman's proof of Thurston's geometrization conjecture
\cite{Perelman1,Perelman2},
using Hamilton's Ricci flow
\cite{Hamilton (1982)}, there
are many remaining questions about three-dimensional
Ricci flow.

Since the Ricci flow is a nonlinear heat equation for the
Riemannian metric, the intuition is that it should smooth out
the metric and thereby give rise, in the long-time limit, to
the locally homogeneous pieces in the geometric decomposition.
This intuition is a bit misleading because, for example, of the
presence of singularities in the Ricci flow. 
Nevertheless, based partly on earlier work of 
Hamilton \cite{Hamilton (1999)}, Perelman showed that the
hyperbolic pieces do asymptotically appear in
the Ricci flow. Perelman's proof for the existence of
the other geometric pieces is more indirect.  Perelman showed
that the nonhyperbolic part of the evolving manifold satisfies
certain geometric conditions, from which one can show that it
is a graph manifold 
\cite{BBBMP (2007),Kleiner-Lott2,Morgan-Tian,Perelman2,Shioya-Yamaguchi (2005)}. 
By earlier work of topologists,
graph manifolds have a geometric decomposition.

It is an open question whether the Ricci flow directly
performs the geometric decomposition of a three-manifold, as
time evolves. In particular, suppose that
the geometric decomposition of the three-manifold consists of a 
single geometric piece. If this piece has Thurston type
$S^3$ or $S^1 \times S^2$ then its Ricci flow has a finite extinction time
\cite{Colding-Minicozzi (2005),Colding-Minicozzi (2007),Perelman3}. 
For the other Thurston types, one can ask whether the large-time
behavior of the  Ricci flow solution will be that of a
locally homogeneous Ricci flow, no matter what the initial metric
may be. Hamilton \cite[Section 11]{Hamilton (1993)},
Hamilton-Isenberg \cite{Hamilton-Isenberg (1993)} and Knopf
\cite{Knopf (2000)} showed
that this is true for certain manifolds of $\R^3$ or $\Sol$-type
if one assumes some extra symmetries on the initial metric.
We are interested in whether one can show asymptotic homogeneity
for a wider class of Ricci flow solutions.

To describe the results,
let $g(\cdot)$ denote a Ricci-flow-with-surgery whose initial
manifold is a closed orientable $3$-manifold. Let $M_t$ denote
the time-$t$ manifold.  (If $t$ is a surgery time then we take
$M_t$ to be the postsurgery manifold.)  From Perelman's
work \cite{Perelman3}, there is some time $T_0$ so that for all
$t \ge T_0$, each connected component $C$ of $M_t$ is $S^3$ or an
aspherical $3$-manifold. As the geometrization
conjecture holds, $C$ has a decomposition into
geometric pieces of type $S^3$, $\R^3$, $H^3$, $\Nil$, $\Sol$, 
$H^2 \times \R$ and
$\widetilde{SL_2(\R)}$; see Section \ref{section2}.

It is possible that the Ricci-flow-with-surgery involves an
infinite number of surgeries.  In the known examples,
there is a finite number of surgeries.
Furthermore, in the known examples, after all of
the surgeries are done then the sectional curvatures uniformly decay in
magnitude as $O(t^{-1})$, i.e. one has a type-III Ricci flow solution.
In order to make progress, we will consider
only Ricci-flows-with-surgery in which this is the case.  Hence,
we will consider a smooth Ricci flow $(M, g(\cdot))$, defined
for $t \in (1, \infty)$ on a closed, connected orientable $3$-manifold
$M$, with sectional curvatures that are uniformly $O(t^{-1})$.

If $M$ admits a locally homogeneous metric
modeled on a given one of the eight Thurston geometries
then we will say that
$M$ has the corresponding Thurston type.  Saying that $M$ has a certain
Thurston type is a topological statement, i.e. we allow ourselves
to consider Riemannian metrics on $M$ that
are not locally homogeneous.

In order to analyze the large-time behavior of a Ricci flow, we
use blowdown limits.

\begin{definition} \label{1.1}
For $s \ge 1$, put
$g_s(t) \: = \: \frac{1}{s} \: g(st)$. It is also a Ricci flow solution. Let
$\widetilde{g}_s(t)$ be the lift of ${g}_s(t)$ to
the universal cover $\widetilde{M}$.
\end{definition}

A time interval $[a,b]$ for $g_s$ corresponds to the time interval
$[sa,sb]$ for $g$. We are interested in the behavior as $s \rightarrow 
\infty$ of $g_s(\cdot)$ on a specified time interval $[a,b]$, 
since this gives information
about the large-time behavior of the initial Ricci flow solution
$g(\cdot)$. If there is a limiting Ricci flow solution 
$\lim_{s \rightarrow \infty} g_s(\cdot)$
then one says that it is a blowdown limit of $g(\cdot)$.

For notation,
if the Gromov-Hausdorff limit 
$\lim_{t \rightarrow \infty}
\left( M, \frac{g(t)}{t} \right)$ exists and equals a compact
metric space $X$ then we write
$\lim_{t \rightarrow \infty}
\left( M, \frac{g(t)}{t} \right) \stackrel{GH}{=} X$. 
If we write
$\lim_{s \rightarrow \infty}
\left( \widetilde{M},\widetilde{m},\widetilde{g}_s(\cdot) \right) =
\left( {M}_\infty, m_\infty, {g}_\infty(\cdot) \right)$ then
we mean that for any sequence $\{ s_j \}_{j=1}^\infty$ tending
to infinity, there is a smooth pointed limit
$\lim_{j \rightarrow \infty} 
\left( \widetilde{M}, \widetilde{m}, \widetilde{g}_{s_j}(\cdot) \right)$
of Ricci flow solutions which
equals $\left( {M}_\infty, m_\infty, {g}_\infty(\cdot) \right)$.
We recall that the notion of the limit in the statement
$\lim_{j \rightarrow \infty} 
\left( \widetilde{M}, \widetilde{m}, 
\widetilde{g}_{s_j}(\cdot) \right) \: = \:
\left( {M}_\infty, m_\infty, {g}_\infty(\cdot) \right)$
involves
$j$-dependent pointed diffeomorphisms from domains in
${M}_\infty$ to domains in $\widetilde{M}$
\cite{Hamilton (1995)}.

\begin{theorem} \label{1.2}
Let $(M, g(\cdot))$ be a smooth Ricci flow solution on a
connected closed orientable $3$-manifold, defined
for $t \in (1, \infty)$. Suppose that \\
1. The sectional curvatures
of $(M, g(t))$ are uniformly $O(t^{-1})$ and \\
2. $\diam(M, g(t)) = O(t^{\frac12})$. \\
Then $M$ is irreducible, aspherical and
its geometric decomposition contains a single geometric piece. \\
1. If $M$ has Thurston type $\R^3$ then $\lim_{t \rightarrow \infty}
\left( M, \frac{g(t)}{t} \right) \stackrel{GH}{=} \pt$. The limit
$\lim_{s \rightarrow \infty}
\left( \widetilde{M}, \widetilde{m},\widetilde{g}_s(\cdot) \right)$ 
exists and equals
the flat expanding soliton $(\R^3, g_{flat})$. \\
2. If $M$ has Thurston type $\Nil$ then $\lim_{t \rightarrow \infty}
\left( M, \frac{g(t)}{t} \right) \stackrel{GH}{=} 
\pt$. The limit
$\lim_{s \rightarrow \infty}
\left( \widetilde{M}, \widetilde{m},\widetilde{g}_s(\cdot) \right)$ 
exists and equals
the expanding soliton $\left( \R^3, \frac{1}{3 t^{\frac13}}
(dx + \frac12 y dz - \frac12 z dy)^2 + t^{\frac13} (dy^2 + dz^2) \right)$.\\
3. If $M$ has Thurston type $\Sol$ then the Gromov-Hausdorff limit
$\lim_{t \rightarrow \infty}
\left( M, \frac{g(t)}{t} \right)$ is a circle or an interval.
 The limit
$\lim_{s \rightarrow \infty}
\left( \widetilde{M}, \widetilde{m},\widetilde{g}_s(\cdot) \right)$ 
exists and equals
the expanding soliton $\left( \R^3, e^{-2z} dx^2 + e^{2z} dy^2 + 4t dz^2
\right)$.\\
4. If $M$ has Thurston type $H^2 \times \R$ then for any sequence
$\{t_j\}_{j=1}^\infty$ tending to infinity, there is a subsequence (which
we relabel as $\{t_j\}_{j=1}^\infty$) so that the
Gromov-Hausdorff limit $\lim_{j \rightarrow \infty}
\left( M, \frac{g(t_j)}{t_j} \right)$ exists and is a metric of constant curvature
$- \: \frac12$ on a closed $2$-dimensional orbifold.
The limit
$\lim_{s \rightarrow \infty}
\left( \widetilde{M}, \widetilde{m},\widetilde{g}_s(\cdot) \right)$ 
exists and equals
the expanding soliton $(H^2 \times \R, 2t g_{hyp} + g_{\R})$. \\
5. If $M$ has Thurston type $H^3$ then $\lim_{t \rightarrow \infty}
\left( M, \frac{g(t)}{t} \right) \stackrel{GH}{=} 
\left( M, 4 g_{hyp} \right)$. The limit
$\lim_{s \rightarrow \infty}
\left( \widetilde{M}, \widetilde{m},\widetilde{g}_s(\cdot) \right)$ 
exists and equals
the expanding soliton $(H^3, 4t g_{hyp})$. \\
6. If $M$ has Thurston type $\widetilde{\SL_2(\R)}$ 
then there is some sequence $\{ s_j \}_{j=1}^\infty$ tending to
infinity such that
$\lim_{j \rightarrow \infty}
\left( M, \frac{g(s_j)}{s_j} \right)$
is a metric of constant curvature
$- \: \frac12$ on a closed $2$-dimensional orbifold and
$\lim_{j \rightarrow \infty}
\left( \widetilde{M}, \widetilde{m},\widetilde{g}_{s_j}(\cdot) \right)$ is
the expanding soliton $(H^2 \times \R, 2t g_{hyp} + g_{\R})$.
\end{theorem}

\begin{corollary} \label{1.3}
In cases 1-5 of Theorem \ref{1.2}, as time becomes large the Ricci flow
solution becomes increasingly locally homogeneous.
\end{corollary}

The corollary follows from the fact that the expanding solitons
in Theorem \ref{1.2} are all homogeneous.
To state the corollary in a more precise way,
we recall that a Riemannian manifold $(M,g)$ is locally homogeneous if
and only if any function on $M$ that can be expressed as a
polynomial in the covariant derivatives of the curvature tensor $\nabla_{i_1} \nabla_{i_2}
\ldots \nabla_{i_r} R_{jklm}$ and the inverse metric
tensor $g^{ij}$, by contracting indices, is actually constant on $M$
\cite{Prufer-Tricerri-Vanhecke (1996)}. 
Corollary \ref{1.3} means that 
in cases 1-5 of Theorem \ref{1.2}, any function on $M$
which is a polynomial in the covariant derivatives of the curvature tensor and the
inverse metric tensor of the
rescaled metric $\widehat{g}(t) = \frac{g(t)}{t}$ approaches a
constant value as $t \rightarrow \infty$.

\begin{remark} \label{1.4}
The diameter condition $\diam(M, g(t)) = O(t^{\frac12})$ implies (under our
curvature assumption) that
the geometric decomposition of $M$ contains a single geometric piece;
see Proposition \ref{3.5}.
We expect that if the diameter condition is not satisfied then
the geometric decomposition of $M$ will contain more than one
geometric piece.
\end{remark}

\begin{remark} \label{1.5}
Any locally homogeneous Ricci flow solution 
$(M, g(\cdot))$ on a
closed $3$-manifold $M$, which exists for $t \in (1, \infty)$, 
does have sectional curvatures
that are uniformly $O(t^{-1})$ and
$\diam(M, g(t)) = O(t^{\frac12})$
\cite{Isenberg-Jackson (1995),Knopf-McLeod (2001)}.
Hence a Ricci flow solution on $M$
that in any reasonable sense approaches a locally homogeneous solution,
as time goes to infinity,
will satisfy the assumptions of Theorem \ref{1.2}.  In this way,
Theorem \ref{1.2} is essentially an if and only if statement.
\end{remark}

\begin{remark}
In Case 6 of Theorem \ref{1.2} we only show that we have the desired limit
for some sequence
$\{s_j \}_{j=1}^\infty$ tending to infinity, not for any such
sequence. The reason is
a technical point about local stability; see Remark \ref{6.5}.
\end{remark}

In \cite[Theorem 1.1]{Lott (2007)} 
we showed that the expanding soliton solutions
listed in Theorem \ref{1.2} are universal attractors  within the space
of homogeneous Ricci flow solutions on Thurston geometries. In
proving Theorem \ref{1.2}, we show that they are global attractors
within the space of Ricci flow solutions that satisfy the given
curvature and diameter assumptions, after passing to the universal cover.

Theorem \ref{1.2} describes the
Gromov-Hausdorff limit of the rescaled Ricci flow solution on $M$ and the
smooth pointed rescaling limit of the lifted Ricci flow solution on
$\widetilde{M}$. In the proof we show there is a rescaling limit
which is a Ricci flow solution on an object that simultaneously
encodes both the
Gromov-Hausdorff limit on $M$ and the smooth limit on $\widetilde{M}$.
This rescaling limit can be considered to give a canonical geometry
for $M$.
A similar phenomenon occurs in the work of Song and Tian 
concerning collapsing in
the K\"ahler-Ricci flow on elliptic fibrations \cite{Song-Tian (2007)}.

There are three main tools in the proof of Theorem \ref{1.2} : 
a compactness theorem, a monotonicity formula and a local stability
result.  The compactness theorem \cite[Theorem 5.12]{Lott (2007)} is
an extension of Hamilton's compactness theorem for Ricci flow
solutions \cite{Hamilton (1995)}. Hamilton's theorem allows one
to take a convergent subsequence of a sequence of pointed Ricci flow
solutions that have uniform curvature bounds on compact time
intervals and a uniform
lower bound on the injectivity radius at the basepoint. 
The rescalings of a Ricci flow solution on a manifold $M$, as considered in
Theorem \ref{1.2}, may collapse,
i.e. the Gromov-Hausdorff limit $X$ may have dimension less than
three. This means that there is no uniform lower bound on the
injectivity radius of the rescaled solution, and so there cannot
be a limiting Ricci flow solution on a $3$-manifold.
Instead, the limiting Ricci flow solution lives on a more general
object called an \'etale groupoid.  Roughly speaking, an \'etale
groupoid combines the notions of manifold and discrete group into
a single object. Its relevance for us comes from the Cheeger-Fukaya-Gromov
theory of bounded curvature collapse \cite{Cheeger-Fukaya-Gromov (1992)},
which implies that a Riemannian manifold which collapses with bounded
sectional curvature will asymptotically acquire extra symmetries.
In Section \ref{section3} we give a brief overview of how collapsing
interacts with Ricci flow.

Under the assumptions of Theorem \ref{1.2}, 
the compactness theorem of \cite{Lott (2007)}
implies that  if $\{s_j\}_{j=1}^\infty$ is a
sequence tending to infinity then after passing to a subsequence,
$\{ \left( M, g_{s_j}(\cdot) \right) \}_{j=1}^\infty$ converges
to a Ricci flow solution $\overline{g}(\cdot)$ 
on a three-dimensional 
\'etale groupoid.  It remains to understand the long-time
behavior of $\overline{g}(\cdot)$. In our case, the relevant
\'etale groupoids arise from locally free abelian group actions. In essence,
we have to understand the long-time behavior of an invariant Ricci flow
solution on the total space of a (twisted) abelian principal bundle
over a compact space $B$. Such a Ricci flow solution $\overline{g}(\cdot)$ 
becomes a coupled system of evolution equations on the lower-dimensional 
space $B$. This is the dimensional reduction part of the title of this
paper.

Our main tool to analyze the long-time behavior of
such a Ricci flow is a modification of the Feldman-Ilmanen-Ni
expanding entropy functional ${\mathcal W}_+$ 
\cite{Feldman-Ilmanen-Ni (2005)}, which in
turn is a variation on Perelman's ${\mathcal W}$-functional \cite{Perelman1}.
More generally, in Section \ref{section4} we describe
versions of the ${\mathcal F}$, ${\mathcal W}$ and ${\mathcal W}_+$ 
functionals that
are adapted for abelian actions.
Using the modified ${\mathcal W}_+$ functional, 
we show that any blowdown limit of
$\overline{g}(\cdot)$ satisfies the harmonic-Einstein equations of
\cite{Lott (2007)}. As we are in dimension three, we can solve
the harmonic-Einstein equations to find the homogeneous expanding
soliton solutions of Theorem \ref{1.2}. 

By these techniques, we show that there is some sequence 
$\{s_j\}_{j=1}^\infty$ tending to infinity so that
$\{ \left( M, g_{s_j}(\cdot) \right) \}_{j=1}^\infty$ converges in
an appropriate sense to a locally homogeneous expanding soliton
solution.  In order to get convergence for all sequences
$\{s_j\}_{j=1}^\infty$ tending to infinity, we use the local
stability of the locally homogeneous expanding solitons, along with some
further arguments.  The local stability is due to Dan Knopf
\cite{Knopf}. An important point is that we only need the local
stability of the locally homogeneous expanding soliton within
the space of Ricci flow solutions with the same abelian symmetry.
Because of this, the local stability issue reduces to an elliptic-type
analysis on the compact quotient space $B$ where one has
compact resolvents, etc. 
For the $\Nil$ and $\Sol$-expanders, the local stability in a somewhat
different sense was considered in \cite{Guenther-Isenberg-Knopf (2006)}.

The outline of this paper is as follows. In Section \ref{section2} we
make some general remarks about Ricci flow and geometrization.
In Section \ref{section3} 
we give an overview of some of the needed results from
\cite{Lott (2007)}. In Section \ref{section4}, which may be of independent
interest, we analyze Ricci flow
solutions with a locally free abelian group 
action. In Section \ref{section5} we
give the classification of the \'etale groupoids
that arise. In Section \ref{section6} we
prove Theorem \ref{1.2}. Further descriptions are given at the 
beginnings of the sections.

I thank Xiaodong Cao, Dan Knopf and Junfang Li for discussions on
the topics of this paper.  I am especially grateful to Dan 
for telling me of his local stability results \cite{Knopf}. Part of
this research was performed while attending the MSRI 2006-2007
program on Geometric Evolution Equations.  I thank MSRI and the
UC-Berkeley Mathematics Department 
for their hospitality, along with the organizers of
the MSRI program for inviting me.

\section{Geometrization Conjecture and Ricci Flow} \label{section2}

In this section we describe what one might expect for the
long-time behavior of the Ricci flow on a compact $3$-manifold $M$,
in terms of the geometric decomposition of $M$. Background information
on the geometrization conjecture is in \cite{Scott (1983)}.

Let $M$ be a connected closed orientable $3$-manifold.
The Kneser-Milnor theorem says that $M$ has a 
connected sum decomposition $M = M_1 \# M_2 \# \ldots \# M_N$ into
so-called prime factors, unique up to permutation.
Thurston's geometrization conjecture says that if $M$ is prime
then there is a (possibly empty) minimal collection of disjoint incompressible
embedded $2$-tori $\{T_i\}_{i=1}^I$ in $M$, unique up to isotopy, 
so that each connected component of $M - \bigcup_{i=1}^I T_i$ admits
a complete locally homogeneous metric of one of the following types : \\
1. A compact quotient of $S^3$, $S^2 \times \R$, $\R^3$, $\Nil$, $\Sol$,
$H^3$, $H^2 \times \R$ or $\widetilde{\SL_2(\R)}$. \\
2. A noncompact finite-volume quotient of $H^3$ or $H^2 \times \R$. \\
3. $\R \times_{\Z_2} T^2$, where the generator of
$\Z_2$ acts by $x \rightarrow -x$ on $\R$ and by the involution
on $T^2$ for which $T^2/\Z_2$ is the Klein bottle $K$.

\begin{remark} \label{2.1}
A finite-volume quotient of 
$S^3$, $S^2 \times \R$, $\R^3$, $\Nil$ or $\Sol$ is
necessarily a compact quotient. Noncompact finite-volume quotients of
$\widetilde{\SL_2(\R)}$ are not on the list, as they are 
diffeomorphic to noncompact
finite-volume quotients of $H^2 \times \R$.
\end{remark} 
\begin{remark} \label{2.2}
If we were to cut along both
2-tori and Klein bottles then we could eliminate the $\R \times_{\Z_2} T^2$
case, which is the total space of a twisted $\R$-bundle over $K$. However,
as we are dealing with orientable manifolds,
it is more natural to only cut along 2-tori.
\end{remark}

We now discuss graph manifolds.  A reference is
\cite[Chapter 2.4]{Matveev (2003)}.
We recall that a compact orientable $3$-manifold $M$ with (possibly empty)
boundary is a {\em graph manifold} if there is a collection of
disjoint embedded $2$-tori $\{T_j\}_{j=1}^J$ so that if we take the metric
completion of $M - \bigcup_{j=1}^J T_j$ (with respect to some
Riemannian metric on $M$) then each connected component is the
total space of a circle bundle over a compact surface. Clearly 
$\partial M$, if nonempty, is a disjoint union of $2$-tori.
The result of gluing two graph manifolds along boundary components
is again a graph manifold (provided that it is orientable). In addition,
the connected sum of two graph manifolds is a graph manifold. In terms
of the Thurston decomposition, a closed orientable prime $3$-manifold $M$
is a graph manifold if and only if it has no hyperbolic pieces. 

We now summarize how Perelman proved the geometrization conjecture
using Ricci flow. If $g(0)$ is an initial Riemannian metric on $M$ then
Perelman showed that there is a Ricci-flow-with-surgery
$(M_t, g(t))$ defined for all $t \in [0, \infty)$ (although 
$M_t$ may become the empty set for large $t$). A singularity in the flow
is handled by letting some connected components go extinct or
by performing surgery. If $t$ is a surgery time
then we let $M_t$ denote the postsurgery manifold $M_t^+$.
Going from a postsurgery manifold $M_t^+$ to the presurgery manifold
$M_t^-$ amounts topologically to performing connected sums on some
components of $M_t^+$, possibly along with a finite number of
$S^1 \times S^2$'s and $\R P^3$'s,
and restoring any factors that went extinct
at time $t$.
From Kneser's theorem, there is some $T_1 > 0$ so that
for a singularity time $t > T_1$, $M_t^+$ differs from
$M_t^-$ by the addition or subtraction
of some $S^3$ factors. That is, after time $T_1$, all
surgeries are topologically trivial. 

Perelman showed that any connected component which goes extinct during
the Ricci-flow-with-surgery is
diffeomorphic to $S^1 \times S^2$, $S^1 \times_{\Z_2} S^2 = \R P^3 \# \R P^3$
or $S^3/\Gamma$,
where $\Gamma$ is a finite subgroup of $\SO(4)$ that acts freely on $S^3$.
He also showed that for large $t$, any connected component $C$
of $M_t$ has a $3$-dimensional submanifold $G$ with (possibly empty)
boundary so that $G$ is a graph manifold, $\partial G$ consists of
incompressible tori in $C$ and $C - G$ admits a 
complete finite-volume hyperbolic metric. Here $G$ is allowed to be
$\emptyset$ or $C$. Using earlier results from $3$-manifold topology,
this is enough to prove the geometrization conjecture.

It is not known whether there is a finite number of surgeries, but
after some time all remaining surgeries will occur in the graph manifold part.
For example, if the original manifold $M$ admits a hyperbolic metric then
there is a finite number of surgeries, since for large time there is no
graph manifold part. We note that one can never exclude singularities
for topological reasons, as the initial metric could always contain a
pinched $2$-sphere.

In \cite{Perelman3}, 
Perelman showed that for large $t$, any connected component of
$M_t$ is aspherical or $S^3$. Thus the relevant Thurston geometries are
$S^3$,
$\R^3$, $\Nil$, $\Sol$, $H^3$, $H^2 \times \R$ and $\widetilde{\SL_2(\R)}$. 

Put $\widehat{g}(t) = \frac{g(t)}{t}$. Let us assume that there is a
finite number of surgeries, and consider the manifold $M$ to be
a connected component of the 
remaining manifold after all of the surgeries are performed.  
Based on explicit calculations for the 
Ricci flow on a locally homogeneous $3$-manifold, the most
optimistic possibility for the Gromov-Hausdorff behavior of the
long-time Ricci flow is given in the following table. Here
$X$ is the Gromov-Hausdorff limit $\lim_{t \rightarrow \infty}
(M, \widehat{g}(t))$, which we assume to exist. The ``Thurston
type'' denotes the possible geometric types in the Thurston decomposition
of $M$, but we do not assume that the metrics in the Ricci flow are
locally homogeneous.
\begin{equation}
\begin{array}{ccc}
\underline{\mbox{X}} & & \underline{\mbox{Thurston type}} \\
 & & \notag \\
\pt. & & \R^3 \text{ or } \Nil \notag \\
S^1 \mbox{ or } I & & \Sol \notag \\
\mbox{closed 2-orbifold with } K = - \: 1/2 & & H^2 \times \R \text{ or }
\widetilde{\SL_2(\R)} \notag \\
\mbox{closed 3-manifold with } K = - \: 1/4 & & H^3 \notag \\
\mbox{noncompact} & &  H^3, H^2 \times \R, \R^3
\end{array}
\end{equation}

If $X$ is noncompact then the possible geometric pieces in the geometric
decomposition of $M$ should be noncompact finite-volume quotients
of $H^3$, noncompact finite-volume quotients of $H^2 \times \R$
and copies of $\R \times_{\Z_2} T^2$. (The final $\R^3$-term in
the table refers to the latter possibility.) When discussing
Gromov-Hausdorff limits in this case, one would
have to choose a basepoint $m \in M$ and take a pointed 
Gromov-Hausdorff limit
$(X, x) \stackrel{GH}{=} 
\lim_{t \rightarrow \infty} (M, m,\widehat{g}(t))$, whose value
would depend on $m$. One would expect to get possible
Gromov-Hausdorff limits of the form \\
1. $H^3/\Gamma$, where $\Gamma$ is a torsion-free noncocompact
lattice in $\PSL(2, \C)$. \\
2. $H^2/\Gamma$, where $\Gamma$ is a noncocompact lattice in
$\PSL(2, \R)$. \\
3. $\R$.  \\
4. $[0, \infty)$. 

\begin{example} \label{2.3}
Suppose that $M = N \cup_{T_2} \overline{N}$ is the double of the 
truncation $N$ of a singly-cusped finite-volume 
hyperbolic $3$-manifold $Y$, where the metric on $N$ is perturbed to 
make it a product near $\partial N$. If $m$ is
in $N - T^2$ then one would expect that 
$\lim_{t \rightarrow \infty} (M,m, \widehat{g}(t)) 
\stackrel{GH}{=} Y$, with a
metric of constant curvature $- \: \frac14$, while if $m \in T^2$
then one would expect that 
$\lim_{t \rightarrow \infty} (M,m, \widehat{g}(t)) \stackrel{GH}{=} \R$.
\end{example}

\begin{example} \label{2.4}
Put $M^\prime =
N \cup_{T^2} (I \times_{Z_2} T^2)$, where $I \times_{\Z_2} T^2$ is
the (orientable) total space of a twisted interval
bundle over the Klein bottle $K$.
Then $M^\prime$ is double covered by $N \cup_{T^2} N$, where the
gluing is done by an orientation-reversing isometry of $T^2$.
If $m \in M^\prime -K$ then one would expect that 
$\lim_{t \rightarrow \infty} (M^\prime, m, \widehat{g}(t)) 
\stackrel{GH}{=} Y$, while
if $m \in K$ then one would expect that 
$\lim_{t \rightarrow \infty} (M^\prime, m, \widehat{g}(t)) 
\stackrel{GH}{=} \R/\Z_2 = [0, \infty)$.

This example shows why, from the point of view of Ricci flow, it
is natural to include $\R \times_{\Z_2} T^2$ as part of the
geometric decomposition; see Remark \ref{2.2}. 
(In this sense it would also be natural to
include $\R \times T^2$ as a possible piece, but such a piece would
be topologically redundant.)
\end{example}

In the collapsing case, i.e. when $\dim(X) < 3$, the 
Gromov-Hausdorff limit $X$ contains limited information about the
evolution of the $3$-dimensional geometry under the Ricci flow.
For $t$ large, any component of the time-$t$ manifold is aspherical
or $S^3$. Because of this, one natural way to get more information
about the $3$-dimensional geometry is
to look at the evolving geometry on the universal
cover. A special case is when $M$ is locally homogeneous.
In \cite[Section 3]{Lott (2007)} the Ricci flow was considered on a 
simply-connected homogeneous $3$-manifold $G/H$, where 
$G$ is a connected unimodular Lie group and $H$ is a
compact subgroup of $G$.
The Ricci flow $(G/H, g(\cdot))$ was assumed to be $G$-invariant and
exist for all positive time.
In each case, it was shown that there are pointed diffeomorphisms
$\{\phi_s\}_{s \in (0, \infty)}$ of $G/H$ so that the blowdown limit
$g_\infty(t) = \lim_{s \rightarrow \infty} \frac{1}{s} \: \phi_s^*
g(st)$ exists and is one of the expanding solitons listed in
Theorem \ref{1.2}.

\begin{remark} \label{2.5}
As an aside,
instead of looking at the rescaled Ricci flow metric
$\widehat{g}(t) \: = \: \frac{g(t)}{t}$, one could also consider
the normalized Ricci flow solution, with constant volume.
The normalized Ricci flow solution is useful in some settings but
in our case we get more uniform results, in terms of the
Thurston type, by looking at $\widehat{g}$.
For example, let $N$ be a truncated singly-cusped finite-volume
hyperbolic $3$-manifold, as in Example \ref{2.3}. Let $\Sigma_1$
and $\Sigma_2$ be 
compact connected surfaces with one boundary
component and negative Euler characteristic. Put 
$M_1 = N \cup_{T^2} (S^1 \times \Sigma_1)$ and
$M_2 = (S^1 \times \Sigma_1) \cup_{T^2} (S^1 \times \Sigma_2)$, where
the gluing of $M_2$ is such that it is not just a product
$S^1 \times (\Sigma_1 \cup_{S^1} \Sigma_2)$. Under the 
unnormalized Ricci flow,
one expects that $\vol(M_1, g(t)) \sim \const t^{3/2}$, due to the
hyperbolic piece, whereas
$\vol(M_2, g(t)) \sim \const t$.  Then the normalized Ricci flow
on $M_1$ should collapse its
$S^1 \times (\Sigma_1 - \partial \Sigma_1)$ piece, while the normalized 
Ricci flow on $M_2$ should have a three-dimensional pointed limit on its
$S^1 \times (\Sigma_1 - \partial \Sigma_1)$ piece. In contrast,
the pointed Gromov-Hausdorff limit
$\lim_{t \rightarrow \infty} (M_i, m_i, \widehat{g}(t))$,
with an appropriate choice of basepoint $m_i$ in
the $S^1 \times \Sigma_1$ piece, should be
$\Sigma_1 - \partial \Sigma_1$ with a complete finite-volume metric
of constant curvature $- \: \frac12$, independent of $i \in \{1,2\}$.  
\end{remark}

\section{Collapsing and Ricci Flow} \label{section3}

In this section we give an overview, aimed for geometers, of the
use of groupoids in collapsing theory.  More details are in
\cite[Section 5]{Lott (2007)} 
and references therein. We also show that under the
hypotheses of Theorem \ref{1.2}, the manifold has a single
geometric piece.

Suppose that $(M^n, g(\cdot))$ is a type-III Ricci flow solution that
exists for $t \in (1, \infty)$, i.e. there is some $K > 0$ so that
$\parallel \Riem(g(t)) \parallel_\infty \: \le \: \frac{K}{t}$ for all
$t > 1$. Then the rescaled metrics $\widehat{g}(t) = \frac{g(t)}{t}$
have uniformly bounded sectional curvature.  
Even if the manifolds
$(M, \widehat{g}(t))$ are collapsing in the Gromov-Hausdorff sense,
we would still like to take a limit as $t \rightarrow \infty$, in some way,
of the $n$-dimensional geometry.
To do so, it is natural to apply
the Cheeger-Fukaya-Gromov theory of bounded curvature collapse to the
Ricci flow. 

A main technique in the Cheeger-Fukaya-Gromov theory is to work 
$O(n)$-equivariantly on the orthonormal frame bundle $FM$.  This is
not very convenient when dealing with Ricci flow, as the induced
flow on $FM$ is complicated.  For this reason, we use an older
approach to collapsing with bounded sectional curvature, as described in
Gromov's book \cite{Gromov (1999)}, that deals directly with the manifold $M$.

Let $M$ be a complete $n$-dimensional Riemannian manifold with sectional
curvatures bounded in absolute value by a positive number $K$.
Given $r \in \left(0, \frac{1}{\sqrt{K}} \right)$ and $m \in M$,
we can consider the Riemannian metric 
$\exp_m^* g$ on $B(0, r) \subset T_mM$. 

Given a sequence of pointed 
complete $n$-dimensional Riemannian manifolds $\{(M_i, m_i)\}_{i=1}^\infty$
with sectional curvatures bounded in absolute value by $K$, there
is a convergent subsequence of the pointed geometries
$B(0, r) \subset T_{m_i}M_i$, whose limit is a $C^{1,\alpha}$-metric on 
an $n$-dimensional
$r$-ball $(B_\infty, m_\infty)$. If one has uniform bounds of the form
$\parallel \nabla^k \Riem (M_i) \parallel_\infty \: \le \: C(k)$
then one can assume that the limit is a $C^\infty$-metric and the
convergence is $C^\infty$.

Define an equivalence relation $\sim_i$ on 
$B \left( 0, \frac{r}{3} \right) \subset T_{m_i}M_i$
by saying that $y \sim_i z$ if $\exp_{m_i} (y)  = \exp_{m_i} (z)$.
Then $B \left( m_i, \frac{r}{3} \right) \subset M_i$ 
equals $(B \left( 0, \frac{r}{3} \right) 
\subset T_{m_i}M_i)/\sim_i$.
The equivalence relation $\sim_i$ is the equivalence relation of a
pseudogroup $\Gamma_i$ of local isometries on $B(0,r) \subset T_{m_i} M_i$,
also called the fundamental pseudogroup $\pi_1(M_i, m_i; r)$. 
One can take a convergent subsequence
of the pseudogroups, in an appropriate sense, to obtain a limit
pseudogroup $\Gamma_\infty$ of local isometries of 
$B_\infty$, which is a local Lie group. Furthermore,
a neighborhood of the identity of $\Gamma_\infty$ is isomorphic to a
neighborhood of the identity of a nilpotent Lie group. In particular,
after passing to the subsequences, the pointed Gromov-Hausdorff limit of
$\{B \left( m_i, \frac{r}{3} \right) 
\subset M_i\}_{i=1}^\infty$ is 
$\left( B \left( m_\infty, \frac{r}{3} \right) 
\subset B_\infty \right)/\Gamma_\infty$.

In this way one constructs a limiting $\frac{r}{3}$-ball. 
It has the drawback that it only describes the (lifted)
geometry near the basepoints $m_i$.
As one started with complete Riemannian manifolds, one would like to
have a limiting object which in some sense is also complete. For
example, suppose that $(M_i, m_i) = (M, m)$ for all $i$. The above
process would produce the limiting ball
$B \left( 0, \frac{r}{3} \right) \subset T_mM$, with 
$\Gamma_\infty = \pi_1(M, m; r)$.  However,
the limiting object should be all of $(M, m)$.

One way to construct a global limiting object would be to move
the basepoints to other points inside of $B \left( 0, \frac{r}{3}
\right) \subset T_{m_i}M_i$,
construct new limiting balls, repeat the process and glue all of
the ensuing balls together in a coherent way. In order to formalize
such a limiting object, the notion of a 
``Riemannian megafold'' was introduced in
\cite{Petrunin-Tuschmann (1999)}. 
This essentially consists of a pseudogroup of local
isometries of a Riemannian manifold. Another formalization was
given in \cite{Lott (2007)}, in which the limiting object is a Riemannian
groupoid.  A Riemannian groupoid is an
\'etale groupoid with a Riemannian metric on its space of units,
for which the local diffeomorphisms coming from groupoid elements are local
isometries. 
Riemannian groupoids have been extensively discussed in the literature
on foliation theory, as
they describe the transverse structure of Riemannian foliations.
For details we refer to
\cite[Section 5]{Lott (2007)} and references therein. 

(We take this opportunity to make some corrections to \cite{Lott (2007)}.
The $\infty$ on \cite[p. 629, line 42]{Lott (2007)} should read $(0, \infty)$.
The $[0,1]$ on 
\cite[p. 658, line 15]{Lott (2007)} should read
$[0,1)$.)

The upshot is that if $\{(M_i, m_i)\}_{i=1}^\infty$ is a sequence
of pointed complete $n$-dimensional Riemannian manifolds, and if
for every $k \in \Z^{\ge 0}$ and $R \in \R^+$ there is some $C(k,R) < \infty$
so that for all $i$ we have 
$| \nabla^k \Riem(M_i) | \le C(k, R)$ on $B(m_i, R) \subset M_i$,
then a subsequence converges smoothly to a pointed complete closed
effective Hausdorff $n$-dimensional Riemannian groupoid
$({\frak G}_\infty, O_{x_\infty})$
\cite[Proposition 5.9]{Lott (2007)}. 
This statement is essentially a reformulation of
results of Cheeger, Fukaya and Gromov.

Let ${\frak G}$ be a complete closed effective Hausdorff  
Riemannian groupoid. It carries a certain locally constant sheaf 
$\underline{\frak g}$ of finite dimensional Lie algebras on its
space of units ${\frak G}^{(0)}$. 
These Lie algebras act as germs of Killing vector
fields on ${\frak G}^{(0)}$. Elements of ${\frak G}$ that are 
sufficiently close to the space of units ${\frak G}^{(0)}$, 
in the $1$-jet topology, appear in the image of the 
exponentials of small local sections of $\underline{\frak g}$.
In our case, the Lie algebras are nilpotent and there is
no point $x \in {\frak G}^{(0)}$ at which all of the corresponding Killing vector
fields vanish simultaneously, unless
$\underline{\frak g} = 0$. We will say that ${\frak G}$ is 
{\em locally free} if
the isotropy groups ${\frak G}^x_x$ are finite.

\begin{remark} \label{3.1}
The locally constant sheaf $\underline{\frak g}$ is analogous to
a {\em pure} $\Nil$-structure in the sense of 
\cite{Cheeger-Fukaya-Gromov (1992)}; see 
\cite{Rong (2007)} for a recent survey. It may seem
surprising that we always get pure structures on our limiting
spaces, since a manifold that collapses with bounded
curvature generally carries a {\em mixed} $\Nil$-structure if the
diameter is not bounded during the collapse. The point is that
we are considering a completely collapsed limit.
In general, given $\epsilon, K, D > 0$, there
is a number $\delta = \delta(n,\epsilon,K,D)$ so that if 
$\parallel \Riem(M) \parallel_\infty \le K$ then
the fundamental pseudogroup $\pi_1(FM, p; \delta)$
(which is represented by loops at $p$ with length less than $\delta$) can
be continuously transported to any point $q \in B(p, D) \subset FM$,
and the result maps into $\pi_1(FM, q; \epsilon)$
\cite[Lemma 7.2]{Fukaya (1993)}. The fact that one generally cannot 
transport the short loops arbitrarily far, while keeping them
short, is responsible for the appearance of 
mixed $\Nil$-structures. As we are considering a
completely collapsed limit, we can effectively move
$\underline{\frak g}_p$ to $\underline{\frak g}_q$ for an
arbitrary value of $D$. 

The work of Cheeger-Fukaya-Gromov describes the local structure
of a Riemannian manifold with bounded sectional curvature that
is highly collapsed but not completely collapsed.  The technique
to do this, for example in \cite{Cheeger-Gromov (1990)}, is to rescale the 
highly-collapsed manifold at a point $p$ in order to make the
rescaled injectivity radius equal to $1$ and the sectional curvatures
very small. One then argues that the local geometry around
$p$ is modeled on a complete flat $n$-dimensional manifold other than
$\R^n$, giving the local $F$-structure. When dealing with Ricci flow
this rescaling is problematic, as it does not
mesh well with the flow.  For this reason, we only deal with
completely collapsed limits.
\end{remark}

The notion of smooth pointed convergence of Riemannian groupoids is
given in \cite{Lott (2007)}, which 
extends these collapsing considerations to the Ricci flow.
(Related Ricci flow limits on a single ball in a tangent space were considered in
\cite{Glickenstein (2003)}.)
In \cite{Lott (2007)} the Ricci flow on an \'etale groupoid was considered.
This consists of a Ricci flow $g(t)$ on the space of units 
${\frak G}^{(0)}$,
in the usual sense, so that for each $t$ 
the local diffeomorphisms (arising from elements of ${\frak G}$)
act by isometries. One has the following
compactness theorem.

\begin{theorem} \label{3.2} \cite[Theorem 5.12]{Lott (2007)}
Let $\{(M_i, p_i, g_i(\cdot))\}_{i=1}^\infty$ be a sequence of
Ricci flow solutions on pointed $n$-dimensional
manifolds $(M_i, p_i)$. We assume that
there are numbers $-\infty \: \le A \: < \: 0$ and $0 \: < \:\Omega
\: \le \: \infty$ so that \\
1. Each Ricci flow solution $(M_i, p_i, g_i(\cdot))$ is defined on the
time interval $(A, \Omega)$. \\
2. For each $t \in (A,\Omega)$, $g_i(t)$ is a complete Riemannian metric
on $M_i$. \\
3. For each compact interval $I  \subset (A, \Omega)$
there is some $K_{I} \: < \: \infty$ so that $|\Riem(g_i)(x, t)| \: \le \:
K_{I}$ for all $x \in M_i$ and $t \in I$. 

Then after passing to a subsequence,
the Ricci flow solutions $g_i (\cdot)$ converge
smoothly to
a Ricci flow solution $g_\infty(\cdot)$ on a pointed $n$-dimensional
\'etale groupoid
$\left( {\frak G}_\infty, O_{x_\infty} \right)$, defined again for
$t \in (A, \Omega)$.
\end{theorem}

This theorem is an analog of Hamilton's compactness theorem
\cite{Hamilton (1995)}, except without the assumption of
a uniform positive lower bound on the
injectivity radius at $p_i \in (M_i, g_i(0))$. In Hamilton's theorem
one obtains a limiting Ricci flow on a manifold, which is a special
type of \'etale groupoid.  The proof of \cite[Theorem 5.12]{Lott (2007)} is
essentially the same as the proof of Hamilton's compactness theorem, when
transplanted to the groupoid setting.

\begin{remark} \label{3.3}
If $\{\diam(M_i, g_i(0))\}_{i=1}^\infty$ is uniformly bounded above then
$({\frak G}_\infty, g_\infty(0))$ 
has finite diameter and we do not have to talk
about basepoints.
\end{remark}

An immediate consequence of Theorem \ref{3.2} is the following.

\begin{corollary} \label{3.4} \cite[Corollary 5.15]{Lott (2007)} 
Given $K > 0$, the space of pointed Ricci flow solutions on 
$n$-dimensional manifolds, with
$\sup_{t \in (1, \infty)} \: t \: \parallel \Riem(g_t) \parallel_\infty
\: \le \: K$, is relatively compact among Ricci flows on
pointed $n$-dimensional \'etale
groupoids, defined for $t \in (1, \infty)$.
\end{corollary}

The next proposition will be used in later sections.

\begin{proposition} \label{3.5}
Let $(M, g(\cdot))$ be a Ricci flow solution on a closed 
orientable $3$-manifold that
is defined for all $t \in [0, \infty)$. Suppose that \\
1. The sectional curvatures
of $(M, g(t))$ are uniformly $O(t^{-1})$ and \\
2. $\diam(M, g(t)) = O(t^{\frac12})$. 

Then $M$ is irreducible, aspherical and
its geometric decomposition contains a single geometric piece.
\end{proposition}
\begin{proof}
As mentioned in Section \ref{section2}, since the Ricci flow exists for all
$t \in [0, \infty)$ it follows that $M$ is aspherical.
The validity of the Poincar\'e Conjecture then implies that
$M$ is irreducible \cite[Theorem 2]{Milnor (1962)}.

Put $\widehat{g}(t) = \frac{g(t)}{t}$. From the evolution equation
for the scalar curvature $R$
and the maximum principle applied to $R \: + \: \frac{3}{2t}$, it
follows that
$\vol(M, \widehat{g}(t))$ is nonincreasing in $t$; see, for example,
\cite[(1.7)]{Feldman-Ilmanen-Ni (2005)}. Suppose that
$\lim_{t \rightarrow \infty} \vol(M, \widehat{g}(t)) > 0$. Then
$(M, \widehat{g}(t))$ is noncollapsing. 
Recall Definition \ref{1.1}.
If $\{s_j \}_{j=1}^\infty$
is a sequence tending to infinity 
then Hamilton's compactness theorem \cite{Hamilton (1995)} implies that
after passing to a subsequence, there is a limiting three-dimensional Ricci flow
solution 
$(M_\infty, g_\infty(\cdot)) = \lim_{j \rightarrow \infty} 
\left( M, g_{s_j}(\cdot) \right)$. From the diameter assumption, 
$M_\infty$ is diffeomorphic to $M$. Using monotonic quantities,
one can show that $(M_\infty, g_\infty(t))$ has constant
sectional curvature $- \: \frac{1}{4t}$; see 
\cite[Section 1]{Feldman-Ilmanen-Ni (2005)} and references therein.
Thus $M$ has an $H^3$-structure.

Now suppose that
$\lim_{t \rightarrow \infty} \vol(M, \widehat{g}(t)) = 0$.
Then $(M, \widehat{g}(t))$ collapses with bounded sectional
curvature and bounded diameter.  There will be a sequence
$t_i \rightarrow \infty$ such that the Gromov-Hausdorff limit
$\lim_{i \rightarrow \infty} (M, \widehat{g}(t_i))$ exists and
equals some compact metric space $X$ of dimension less than
three. In what follows we use some results about
bounded curvature collapsing from \cite{Rong (2007)} and references therein. 

If $\dim(X) = 0$ then $M$ is an almost flat manifold
and so has an $\R^3$ or $\Nil$-structure \cite{Gromov (1978)}. 

If $\dim(X) = 2$ then
$X$ is a closed orbifold and $M$ is the total space of an orbifold
circle bundle over $X$
\cite[Proposition 11.5]{Fukaya (1990)},
from which it follows that $M$ has a
geometric structure. 

Finally, suppose that $\dim(X) = 1$.
First, $X$ is $S^1$ or an interval. If $X = S^1$ then
$M$ is the total space of a torus bundle over $S^1$ and
hence carries a geometric structure. If $X$ is an interval $[0,L]$
then there is a Gromov-Hausdorff approximation
$\pi \: : \: M \rightarrow X$ with 
$\pi^{-1}(0,L) = (0,L) \times T^2$. Now $X$ is locally
the quotient of $M$ by a fixed-point free $T^2$-action
\cite{Cheeger-Gromov (1986)}.
If the action is locally free then $[0,L]$ is an
orbifold. As the orbifold $[0,L]$ is double covered by $S^1$, 
the manifold $M$ is double covered by a $T^2$-bundle
over $S^1$. Hence in this case, $M$ has a geometric
structure
\cite{Meeks-Scott (1986)}.

Suppose that the $T^2$-action is not locally free, say
on $\pi^{-1}[0, \delta)$, with $\delta$ small. From the slice theorem, a
neighborhood of $\pi^{-1}(0)$ is equivariantly diffeomorphic
to $T^2 \times_H \R^N$, where $H$ is the isotropy group. As
the $T^2$-action has no fixed points, $H$ must be a virtual
circle group. However, since $M$ is aspherical, the map
$\pi_1(T^2) \rightarrow \pi_1(M)$ must be injective
\cite[Remark 0.9]{Cheeger-Rong (1995)}. This is a contradiction.
Similarly, the $T^2$-action must be locally free on
$\pi^{-1}(L - \delta, L]$.
\end{proof}

\begin{remark} \label{3.6}
In our case, one can
see directly that there is a contradiction if $H$ is a virtual
circle group. Suppose so.
Then $\pi^{-1}([0, \delta])$ 
(or $\pi^{-1}([L- \delta,L])$)
is diffeomorphic
to $S^1 \times D^2$.
If the $T^2$-action fails to be locally free on both
$\pi^{-1}([0, \delta])$ and $\pi^{-1}([L- \delta,L])$ then
$M$ is the union of two solid tori and so is diffeomorphic to $S^3$,
$S^1 \times S^2$ or a lens space. If it fails to be locally free
on exactly one of 
$\pi^{-1}([0, \delta])$ and $\pi^{-1}([L- \delta,L])$ then a
double cover of $M$ is diffeomorphic to $S^3$, $S^1 \times S^2$ or
a lens space. In either case, $M$ fails to be aspherical.
\end{remark}

\begin{remark} \label{3.7}
By the argument of the proof of Proposition  \ref{3.5}, we can say the
following about aspherical $3$-manifolds that collapse
 with bounded curvature and bounded diameter. 
If $M$ carries an $H^3$-structure then it cannot collapse. If
$M$ carries an $H^2 \times \R$ or $\widetilde{\SL_2(\R)}$-structure 
then it can only collapse to a two-dimensional
orbifold of negative Euler characteristic. If $M$ carries
a $\Sol$-structure then it can only collapse to $S^1$ or an
interval. However, if $M$ carries an $\R^3$ or $\Nil$-structure
then {\it a priori} it could collapse to a two-dimensional
orbifold with vanishing Euler characteristic, a circle,
an interval or a point.  We will show that under the Ricci flow,
with our curvature and diameter assumptions it can only 
collapse to a point. 
\end{remark}

\begin{remark}
Some results about three-dimensional type-IIb Ricci flow solutions,
i.e. Ricci flow solutions defined on $[0, \infty)$ with 
$\limsup_{t \rightarrow \infty} t \: \parallel \Riem(g(t)) \parallel_\infty
\: = \: \infty$, were obtained in \cite{Chow-Glickenstein-Lu (2006)}.
Although phrased differently, the collapsing results in 
\cite{Chow-Glickenstein-Lu (2006)} can be
considered to be results about Ricci flow solutions with nonnegative
sectional curvature on three-dimensional
\'etale groupoids.
\end{remark}

\section{Dimensional Reduction} \label{section4}

In this section we consider a Ricci flow $(M, \overline{g}(\cdot))$ which is
invariant under local actions of a connected abelian Lie group
on $M$. We first define the notion of a twisted principal bundle
and write out the Ricci flow equation for an invariant metric
$\overline{g}$ on the total space $M$. The Ricci flow equation becomes
as a coupled system of equations on the base $B$ of the twisted
principal bundle. We construct modified ${\mathcal F}$, ${\mathcal W}$
and ${\mathcal W}_+$ functionals for $\overline{g}(\cdot)$ and show
that they are monotonic. We use ${\mathcal W}_+$ to show that
any blowdown limit of $\overline{g}(\cdot)$ satisfies the harmonic-Einstein
equations of \cite{Lott (2007)}.

Related functionals were considered independently by
Bernhard List \cite{List (2006)} and Jeff Streets \cite{Streets (2007)}. 
(I thank Gerhard Huisken and Gang Tian for these references.) 
In \cite{List (2006)} the modified ${\mathcal F}$-functional is
considered in the special case when $N = 1$ and $A^i_\alpha = 0$.
The motivation comes from the static Einstein equation.
In \cite{Streets (2007)}, modified ${\mathcal F}$ and 
${\mathcal W}$-functionals are considered for a certain invariant flow on
the total space of a principal bundle, with the fiber geometry being
fixed under the flow.

\subsection{Twisted principal bundles} \label{subsection4.1}

Let ${\mathcal G}$ be a Lie group, with Lie algebra
${\frak g}$. Let $B$ be a connected $n$-dimensional smooth manifold.
Let ${\mathcal E}$ be a local system on $B$ of Lie groups isomorphic to ${\mathcal G}$.
Fixing a basepoint $b_0 \in B$ and an isomorphism ${\mathcal E}_{b_0} \cong {\mathcal G}$
of the stalk over $b_0$, the
local system is specified by a homomorphism $\rho \: : \: \pi_1(B,b) \rightarrow
\Aut({\mathcal G})$. Equivalently, we have a ${\mathcal G}$-bundle 
$E = {\mathcal G} \times_\rho \widetilde{B}$ over $B$, with a flat connection,
which gives the \'etale space of the locally constant sheaf ${\mathcal E}$. 
Put $e = {\frak g} \times_{\rho} \widetilde{B}$, a flat
${\frak g}$-vector bundle on $B$. We will write
$\Lambda^{max} e = \Lambda^{max} {\frak g} \times_{\rho} \widetilde{B}$
for the corresponding flat real line bundle
on $B$ of fiberwise volume forms, and
$|\Lambda^{max} e| = |\Lambda^{max} {\frak g}| \times_{\rho} \widetilde{B}$
for the flat $\R^{\ge 0}$-bundle of fiberwise densities.

Hereafter we assume that
the density bundle $|\Lambda^{max}e|$ is a flat product bundle
$\R^{\ge 0} \times B$. (Some of the subsequent results do not need
this assumption, but for simplicity we will assume uniformly that
it holds.)

\begin{example} \label{4.1}
If ${\mathcal G} = \R^N$ then $\Aut({\mathcal G}) = \GL(N, \R)$ and
$E=e$ is a flat $\R^N$-bundle over $B$. 
The assumption on $|\Lambda^{max}e|$ means that the holonomy of
$e$ lies in $\det^{-1}(\pm 1)$.
If ${\mathcal G} = T^N$ then
$\Aut({\mathcal G}) = \GL(N, \Z)$, $E$ is a flat
$T^N$-bundle over $B$ and $e$ is a flat $\R^N$-bundle over $B$.
In this case the assumption on $|\Lambda^{max}e|$ holds automatically.
\end{example}

Let $\pi \: : \: M \rightarrow B$ be a fiber bundle with fiber ${\mathcal G}$.
We write $E_b$ for the fiber of $E$ over $b \in B$ and
$M_b$ for the fiber of $M$ over $b \in B$. 
Consider the fiber product $E \times_B M = \bigcup_{b \in B} E_b \times M_b$.
We assume that there is a smooth map $E \times_B M \rightarrow M$
so that
over a point $b \in B$, the map $E_b \times M_b \rightarrow M_b$
gives a free transitive action of ${\mathcal G} \cong E_b$ on $M_b$.
The action must be consistent with the 
flat connection on $E$ in the sense that 
if $U \subset B$ is such that
$E \Big|_U \cong U \times {\mathcal G}$ is a local trivialization of the flat
${\mathcal G}$-bundle $E$ then $\pi^{-1}(U)$ has a free ${\mathcal G}$-action, and so is the total 
space of a principal ${\mathcal G}$-bundle over $U$. In this way, $M$ can be
considered to be a twisted principal ${\mathcal G}$-bundle over $B$, with
the twisting coming from the flat ${\mathcal G}$-bundle $E$. 
There is a natural isomorphism between the vertical tangent
bundle $T^{vert}M = \Ker(d\pi)$ and $\pi^* e$.

An isomorphism of two
twisted principal ${\mathcal G}$-bundles $\pi \: : \: M \rightarrow B$ and 
$\pi^\prime \: : \: M^\prime \rightarrow B^\prime$ is given by a
diffeomorphism $\eta \: : \: B \rightarrow B^\prime$, an isomorphism
$\hat{\phi} \: : \: E \rightarrow E^\prime$ of flat 
${\mathcal G}$-bundles that covers $\eta$, and a diffeomorphism 
$\phi \: : \: M \rightarrow M^\prime$ that covers $\eta$ with the
property that for all
$m \in M$ and $x \in E_{\pi(m)}$, we have
$\phi(x \cdot m) \: = \: \widehat{\phi}(x) \cdot \phi(m)$.

It makes
sense to talk about a connection
$A \in \Omega^1(M; \pi^* e)$ on
a twisted principal ${\mathcal G}$-bundle $M$. 
The restriction of $A$ to $\pi^{-1}(U)$ is a ${\frak g}$-valued
connection in the usual sense. 

We assume that $M$
has a Riemannian metric $\overline{g}$ with a local free isometric
${\mathcal G}$-action.  This means that if $E \Big|_U \cong U \times 
{\mathcal G}$ is
a local trivialization of $E$ as above then the action of
${\mathcal G}$ on $\pi^{-1}(U)$ is isometric.

Hereafter we assume that ${\mathcal G}$ is a connected $N$-dimensional
abelian Lie group.

Suppose that $U$ is also small enough so that 
$U$ is a
coordinate chart for $B$ with local parametrization 
$\{x^\alpha\}_{\alpha = 1}^n \rightarrow 
\rho(x^\alpha) \in U$. Take a section
$s \: : \: U \rightarrow \pi^{-1}(U)$.
Choosing a basis
$\{e_i\}_{i=1}^N$ of ${\frak g}$,
we obtain coordinates on 
$\pi^{-1}(U)$ by
$(x^\alpha, x^i) \rightarrow \exp \left( \sum_{i=1}^N x^i e_i \right) \cdot 
s(\rho(x^\alpha))$. In terms of these coordinates we can write
\begin{equation} \label{4.2}
\overline{g} \: = \: \sum_{i,j=1}^N G_{ij} \: (dx^i + A^i) (dx^j + A^j) \: + \:
\sum_{\alpha, \beta = 1}^n g_{\alpha \beta} \: dx^\alpha dx^\beta.
\end{equation}
Here
$G_{ij}$ is the local expression of a Euclidean inner product on
$e$,
$\sum_{\alpha, \beta = 1}^n g_{\alpha \beta} \: dx^\alpha dx^\beta$ is
the local expression of a Riemannian metric $g_B$ on $B$ and
$A^i = \sum_{\alpha} A^i_\alpha dx^\alpha$ are the components of 
$s^* A$.
A change of section $s$ changes $A^i$ by
an exact form. The curvatures $F^i = dA^i$ form an element of
$\Omega^2(B; e)$.

If $M$ and $M^\prime$ are two 
twisted principal ${\mathcal G}$-bundles
then an isomorphism
$\phi \: : \: M \rightarrow M^\prime$
can be written in local coordinates as
\begin{equation} \label{4.3}
\phi(y^\gamma, y^k) \: = \:
(x^\alpha (y^\gamma), \sum_k T^i_{\: \: k} y^k \: + \: f^i(y^\gamma)),
\end{equation}
where the $T^i_{\: \: k}$'s are constants. It covers a diffeomorphism
$\eta \: : \: B \rightarrow B^\prime$. 
The isomorphism 
$\widehat{\phi} \: : \: E \rightarrow E^\prime$
of flat ${\mathcal G}$-bundles
is represented locally by the functions $T^i_{\: \: k}$.

A locally ${\mathcal G}$-invariant Ricci flow is a $1$-parameter family of such
Riemannian metrics $(M, \overline{g}(\cdot))$ that satisfies the Ricci flow
equation. We will consider a basepoint for such a solution to be a point
$p \in B$.

Let $\{(M_i, p_i, \overline{g}_i(\cdot))\}_{i=1}^\infty$ be a sequence of
locally ${\mathcal G}$-invariant Ricci flow solutions defined for 
$t \in (1, \infty)$.
We say that $\lim_{i \rightarrow \infty}
(M_i, p_i, \overline{g}_i(\cdot)) \: = \: (M_\infty, p_\infty, \overline{g}_\infty(\cdot))$ 
if there are \\
1. A sequence of open subsets $\{U_j\}_{j=1}^\infty$ of $B_\infty$, containing
$p_\infty$, so that any compact subset of $B_\infty$ eventually lies in all $U_j$, and \\
2. Open subsets $V_{i,j} \subset B_i$ containing $p_i$ and 
isomorphisms $\phi_{i,j} \: : \: \pi_\infty^{-1}(U_j) \rightarrow
\pi_i^{-1}(V_{i,j})$
sending $\pi_\infty^{-1}(p_\infty)$ to
$\pi_i^{-1}(p_i)$ so that \\
3. For all $j$, $\lim_{i \rightarrow \infty} \phi_{i,j}^* \: \overline{g}_i(\cdot) \: = \:
\overline{g}_\infty(\cdot)$ smoothly on $\pi_\infty^{-1}(U_j) \times 
[1 + j^{-1}, 1+j]$.

If $B_\infty$ is compact then we can remove the reference to basepoints.

\subsection{Ricci flow on twisted principal bundles} \label{subsection4.2}

In what follows, we use the Einstein summation convention freely.
Let $(x^\alpha, x^i)$
be local coordinates on $\pi^{-1}(U)$ as in Subsection 
\ref{subsection4.1}.
Writing $A^i \: = \: 
\sum_{\alpha=1}^n A^i_\alpha \: dx^\alpha$,
put $F^i_{\alpha \beta} = \partial_\alpha A^i_\beta - 
\partial_\beta A^i_\alpha$.
We also write
\begin{equation} \label{4.4}
G_{ij;\alpha \beta} \: = \: G_{ij,\alpha \beta} \: - \:
\Gamma^{\sigma}_{\: \: \alpha \beta} \: G_{ij, \sigma},
\end{equation}
where $\{\Gamma^{\sigma}_{\: \: \alpha \beta}\}$ are the
Christoffel symbols for the metric $g_{\alpha \beta}$ on $B$.

Given $b \in U$, it is convenient to choose the section $s$ so that
$A^i(b) = 0$.
Then the curvature tensor 
$\overline{R}_{IJKL}$ of $M$ is given in terms of the curvature tensor
$R_{\alpha \beta \gamma \delta}$ of $B$, the $2$-forms $F^i_{\alpha \beta}$
and the metrics $G_{ij}$ by
\begin{align} \label{4.5}
\overline{R}_{ijkl} \: & = \: - \: \frac14 \: 
g^{\alpha \beta} \:
G_{ik,\alpha} \: G_{jl,\beta} \:  + \: \frac14 \:
g^{\alpha \beta} \:
G_{il,\alpha} \: G_{jk, \beta} \\
\overline{R}_{ijk\alpha} \: & = \: \frac14 \: g^{\beta \gamma} \:
G_{jm} \:  G_{ik,\beta} \: F^m_{\alpha \gamma} \: - \: \frac14 \:
g^{\beta \gamma} \: G_{im} \: G_{jk,\beta} \: F^m_{\alpha \gamma} \notag \\
\overline{R}_{ij \alpha \beta} \: & = \: - \: \frac14 \:
G^{mk} \: G_{im,\alpha} \: G_{kj,\beta} \: + \: \frac14 \:
G^{mk} \: G_{im,\beta} \: G_{kj,\alpha} \: - \: \frac14 \:
g^{\gamma \delta} \: G_{im} \: G_{jk} \: F^m_{\alpha \gamma} \: 
F^k_{\beta \delta} \: + \:
\frac14 \:
g^{\gamma \delta} \: G_{im} \: G_{jk} \: F^m_{\beta \gamma} \: 
F^k_{\alpha \delta} \notag \\
\overline{R}_{i\alpha j \beta} \: & = \: - \: \frac12 \:
G_{ij;\alpha \beta} \: + \: \frac14 \: G^{kl} \:
G_{ik, \beta} \: G_{jl, \alpha} \: + \: \frac14 \:
g^{\gamma \delta} \: G_{ik} \: G_{jl} \: F^k_{\alpha \gamma} \: 
F^l_{\beta \delta} \notag \\
\overline{R}_{i \alpha \beta \gamma} \: & = \:
\frac12 \: G_{ij} \: F^j_{\beta \gamma; \alpha} \: + \: \frac12 \: 
G_{ij,\alpha} \: F^j_{\beta \gamma} \: + \: \frac14 \:
G_{ij, \beta} \: F^j_{\alpha \gamma} \: - \: \frac14 \:
G_{ij, \gamma} \: F^j_{\alpha \beta} \notag \\
\overline{R}_{\alpha \beta \gamma \delta} \: & = \:
{R}_{\alpha \beta \gamma \delta} \: - \: \frac12 \: G_{ij} \:
F^i_{\alpha \beta} \: F^j_{\gamma \delta} - \: \frac14 \: G_{ij} \:
F^i_{\alpha \gamma} \: F^j_{\beta \delta} + \: \frac14 \: G_{ij} \:
F^i_{\alpha \delta} \: F^j_{\beta \gamma}. \notag
\end{align}
The Ricci tensor is given by
\begin{align} \label{4.6}
\overline{R}_{ij} \: & = \: - \: \frac12 \: g^{\alpha \beta} \: 
G_{ij; \alpha \beta} \: - \: \frac14 \: g^{\alpha \beta} \:
G^{kl} \: G_{kl, \alpha} \: G_{ij, \beta} \: + \:
\frac12 \: g^{\alpha \beta} \: G^{kl} \: G_{ik, \alpha} \:
G_{lj, \beta} \: + \: \frac14 \: g^{\alpha \gamma} \: g^{\beta \delta} \:
G_{ik} \: G_{jl} \: F^k_{\alpha \beta} \: F^l_{\gamma \delta} \\
\overline{R}_{i \alpha} \: & = \: \frac12 \: g^{\gamma \delta} \:
G_{ik} \: F^k_{\alpha \gamma; \delta} \: + \: \frac12 \: g^{\gamma \delta} \:
G_{ik, \gamma} \: F^k_{\alpha \delta} \: + \: \frac14 \:
g^{\gamma \delta} \: G_{im} \: G^{kl} \: G_{kl, \gamma} \: F^m_{\alpha \delta} 
\notag \\
\overline{R}_{\alpha \beta} \: & =  \: R_{\alpha \beta} \: - \: 
\frac12 \: G^{ij} \: G_{ij; \alpha \beta} \: + \: \frac14 \:
G^{ij} \: G_{jk,\alpha} \: G^{kl} \: G_{li,\beta} \: - \:
\frac12 \: g^{\gamma \delta} \: \: G_{ij} \: F^i_{\alpha \gamma} \:
F^j_{\beta \delta}. \notag
\end{align}
The scalar curvature is
\begin{equation} \label{4.7}
\overline{R} \: = \: 
R \: - \: g^{\alpha \beta} G^{ij}  \: G_{ij; \alpha \beta} \: + \: 
\frac34 \: g^{\alpha \beta} \: G^{ij} \: G_{jk, \alpha} \: G^{kl} \:
G_{li, \beta} \: - \: \frac14 \: g^{\alpha \beta} \: G^{ij} \: 
G_{ij, \alpha} \: G^{kl} \: G_{kl, \beta} \: - \: \frac14 \:
g^{\alpha \gamma} \: g^{\beta \delta} \: G_{ij} \:
F^i_{\alpha \beta} \: F^j_{\gamma \delta}.
\end{equation}

Consider a $1$-parameter family of such Riemannian metrics 
$\overline{g}(\cdot)$ on $M$. Writing 
$G_{ij}(t)$, $A^i_\alpha(t)$ and
$g_{\alpha \beta}(t)$ as
functions of $t$, the Ricci flow equation becomes
\begin{align} \label{4.8}
& \frac{d}{dt} \left( G_{ij} (dx^i + A^i)(dx^j + A^j) \: + \: 
g_{\alpha \beta} dx^\alpha dx^\beta \right)
\: = \\
&  -2 \overline{R}_{ij} \: (dx^i + A^i)(dx^j + A^j) \: - \: 4 \: 
\overline{R}_{i \alpha} \:
(dx^i + A^i) dx^\alpha  \: - \: 2 \: 
\overline{R}_{\alpha \beta} dx^\alpha dx^\beta. \notag
\end{align}
Equivalently,
\begin{align} \label{4.9}
\frac{\partial G_{ij}}{\partial t} \: & = \: g^{\alpha \beta} \: 
G_{ij; \alpha \beta} \: + \: \frac12 \: g^{\alpha \beta} \:
G^{kl} \: G_{kl, \alpha} \: G_{ij, \beta} \: - \:
g^{\alpha \beta} \: G^{kl} \: G_{ik, \alpha} \:
G_{lj, \beta} \: - \: \frac12 \: g^{\alpha \gamma} \: g^{\beta \delta} \:
G_{ik} \: G_{jl} \:  F^k_{\alpha \beta} \: F^l_{\gamma \delta} \\ \notag
\frac{\partial A^i_{\alpha}}{\partial t} \: & = \: 
- \: g^{\gamma \delta} \:
F^i_{\alpha \gamma; \delta} \: - \: g^{\gamma \delta} \: G^{ij} \:
G_{jk, \gamma} \: F^k_{\alpha \delta} \: - \: \frac12 \:
g^{\gamma \delta} \: G^{kl} \: G_{kl, \gamma} \: F^i_{\alpha \delta} 
\notag \\
\frac{\partial g_{\alpha \beta}}{\partial t} \: & =  \: 
-2 R_{\alpha \beta} \: + \: 
G^{ij} \: G_{ij; \alpha \beta} \: - \: \frac12 \:
G^{ij} \: G_{jk,\alpha} \: G^{kl} \: G_{li,\beta} \: + \:
g^{\gamma \delta} \: \: G_{ij} \: F^i_{\alpha \gamma} \:
F^j_{\beta \delta}. \notag
\end{align}
Adding a Lie derivative with respect to $- \: \nabla \ln \sqrt{\det(G_{ij})}$
to the
right-hand side, and adding an exact form to the right-hand side
of the equation for $\frac{\partial A^i_\alpha}{\partial t}$, gives
a new equivalent set of equations :
\begin{align} \label{4.10}
\frac{\partial G_{ij}}{\partial t} \: & = \: g^{\alpha \beta} \: 
G_{ij; \alpha \beta} \: - \:
g^{\alpha \beta} \: G^{kl} \: G_{ik, \alpha} \:
G_{lj, \beta} \: - \: \frac12 \: g^{\alpha \gamma} \: g^{\beta \delta} \:
G_{ik} \: G_{jl} \:  F^k_{\alpha \beta} \: F^l_{\gamma \delta} \\
\frac{\partial A^i_{\alpha}}{\partial t} \: & = \: 
- \: g^{\gamma \delta} \:
F^i_{\alpha \gamma; \delta} \: - \: g^{\gamma \delta} \: G^{ij} \:
G_{jk, \gamma} \: F^k_{\alpha \delta}
\notag \\
\frac{\partial g_{\alpha \beta}}{\partial t} \: & =  \: 
-2 R_{\alpha \beta} \: + \:  \frac12 \:
G^{ij} \: G_{jk,\alpha} \: G^{kl} \: G_{li,\beta} \: + \:
g^{\gamma \delta} \: \: G_{ij} \: F^i_{\alpha \gamma} \:
F^j_{\beta \delta}. \notag
\end{align}

The equations in (\ref{4.10}) consist of a heat type
equation for $G_{ij}$, a Yang-Mills gradient flow type equation
for $A^i_\alpha$ and a Ricci flow type equation for $g_{\alpha \beta}$.
If $B$ is closed then an extension of the DeTurck trick
\cite{DeTurck (1983)}
to our setting shows short-time existence and uniqueness
for the system (\ref{4.10}). 

\subsubsection{Modified ${\mathcal F}$-functional} \label{subsubsection4.2.1}

We now assume that $B$ is closed.  

\begin{definition} \label{4.11} Given $f \in C^\infty(B)$, put
\begin{align} \label{4.12}
& {\mathcal F}(G_{ij}, A^i_\alpha, g_{\alpha \beta}, f) \: = \\
& \int_B \left( |\nabla f|^2 \: + \: R \: - \: \frac14 \:
g^{\alpha \beta} \: G^{ij} \: G_{jk, \alpha} \: G^{kl} \:
G_{li, \beta} \: - \: \frac14 \: g^{\alpha \gamma} \: g^{\beta \delta} \:
G_{ij} \: F^i_{\alpha \beta} \: F^j_{\gamma \delta} \right)
\: e^{-f} \: \dvol_B. \notag
\end{align}
\end{definition}

If $N = 0$, i.e. if $M = B$, then this is the same as
Perelman's ${\mathcal F}$-functional \cite{Perelman1}.
Otherwise, the expression in (\ref{4.12}) differs from Perelman's
${\mathcal F}$-functional by the subtraction of terms
corresponding to a Dirichlet energy of the field $G$ and
a Yang-Mills action for the connection $A$.

We now compute the variation of ${\mathcal F}$.

\begin{lemma} \label{4.13}
Given a smooth $1$-parameter family $\{(G_{ij}(s), A^i_\alpha(s), 
g_{\alpha \beta}(s), 
f(s))\}_{s \in
(-\epsilon, \epsilon)}$, write 
$\dot{G}_{ij} \: = \: \frac{dG_{ij}}{ds} \Big|_{s=0}$,
$\dot{A}^i_\alpha \: = \: \frac{dA^i_\alpha}{ds} \Big|_{s=0}$,
$\dot{g}_{\alpha \beta} \: = \: \frac{dg_{\alpha \beta}}{ds} \Big|_{s=0}$ and
$\dot{f} \: = \: \frac{df}{ds} \Big|_{s=0}$.
Then
\begin{align} \label{4.14}
& \frac{d}{ds} \Big|_{s=0} {\mathcal F}(G_{ij}, A^i_\alpha, 
g_{\alpha \beta}, f) 
\: = \\
& - \: \int_B \dot{G}_{kl} \: G^{ik} \: G^{jl} \notag \\
& \left(
- \: \frac12 \: g^{\alpha \beta} \: 
G_{ij; \alpha \beta} \: + \:
\frac12 \: g^{\alpha \beta} \: G^{kl} \: G_{ik, \alpha} \:
G_{lj, \beta} \: + \: \frac14 \: g^{\alpha \gamma} \: g^{\beta \delta} \:
G_{ik} \: G_{jl} \: F^k_{\alpha \beta} \: F^l_{\gamma \delta}
\: + \: \frac12 \: g^{\alpha \beta} \: G_{ij, \alpha} \: f_{,\beta}
\right) e^{-f} \: \dvol_B \: - \notag \\
& 2 \int_B \dot{A}^j_\beta \: g^{\alpha \beta} \: G_{ij} \left(
\frac12 \: g^{\gamma \delta} \:
F^i_{\alpha \gamma; \delta} \: + \: \frac12 \: g^{\gamma \delta} \:
G^{ij} \:
G_{jk, \gamma} \: F^k_{\alpha \delta} \: - \: \frac12 \: g^{\gamma \delta} \:
f_{, \gamma} \: F^k_{\alpha \delta} \right) \:
e^{-f} \: \dvol_B \: - \notag \\
& \int_B \dot{g}^{\alpha \beta} \left(
R_{\alpha \beta} \: - \: \frac14 \:
G^{ij} \: G_{jk,\alpha} \: G^{kl} \: G_{li,\beta} \: - \:
\frac12 g^{\gamma \delta} \: \: G_{ij} \: F^i_{\alpha \gamma} \:
F^j_{\beta \delta} \: + \: f_{;\alpha \beta} \right) \: e^{-f} \: \dvol_B \: + \notag \\
& 
\int_B \left( \frac12
g^{\alpha \beta} \dot{g}_{\alpha \beta} \: - \: \dot{f} \right) \notag \\
& \left( 2 \nabla^2 f - |\nabla f|^2 + 
R \: - \: \frac14 \:
g^{\alpha \beta} \: G^{ij} \: G_{jk, \alpha} \: G^{kl} \:
G_{li, \beta} \: - \: \frac14 \: g^{\alpha \gamma} \: g^{\beta \delta} \:
G_{ij} \: F^i_{\alpha \beta} \: F^j_{\gamma \delta}
\right) \: e^{-f} \: \dvol_B. \notag
\end{align}
\end{lemma}
\begin{proof}
This follows from a calculation along the lines of the corresponding
calculation for Perelman's ${\mathcal F}$-functional; see
\cite[Section 5]{Kleiner-Lott}.
\end{proof}

As a consequence of Lemma \ref{4.13}, we can show that
${\mathcal F}$ is nondecreasing under a certain flow.

\begin{corollary} \label{4.15}
Under the flow equations
\begin{align} \label{4.16}
\frac{\partial G_{ij}}{\partial t} \: & = \: 
 g^{\alpha \beta} \: 
G_{ij; \alpha \beta} \: - \:
g^{\alpha \beta} \: G^{kl} \: G_{ik, \alpha} \:
G_{lj, \beta} \: - \: \frac12 \: g^{\alpha \gamma} \: g^{\beta \delta} \:
G_{ik} \: G_{jl} \: F^k_{\alpha \beta} \: F^l_{\gamma \delta}
\: - \: g^{\alpha \beta} \: G_{ij, \alpha} \: f_{,\beta} \\
\frac{\partial A^i_\alpha}{\partial t} \: & = \: 
- \: g^{\gamma \delta} \:
F^i_{\alpha \gamma; \delta} \: - \: g^{\gamma \delta} \:
G^{ij} \:
G_{jk, \gamma} \: F^k_{\alpha \delta} \: + \: g^{\gamma \delta} \:
f_{, \gamma} \: F^k_{\alpha \delta} \notag \\
\frac{\partial g_{\alpha \beta}}{\partial t} \: & = \: - \: 2 \:
R_{\alpha \beta} \: + \: \frac12 \:
G^{ij} \: G_{jk,\alpha} \: G^{kl} \: G_{li,\beta} \: + \:
g^{\gamma \delta} \: \: G_{ij} \: F^i_{\alpha \gamma} \:
F^j_{\beta \delta} \: - \: 2 \: f_{;\alpha \beta} \notag \\
\frac{\partial f}{\partial t} \: & = \: - \:
R \: + \: \frac14 \:
g^{\alpha \beta} \: G^{ij} \: G_{jk, \alpha} \: G^{kl} \:
G_{li, \beta} \: + \: \frac12 \: g^{\alpha \gamma} \: g^{\beta \delta} \:
G_{ij} \: F^i_{\alpha \beta} \: F^j_{\gamma \delta} \: - \:
\nabla^2 f \notag
\end{align}
one has
\begin{align} \label{4.17}
& \frac{d}{dt} {\mathcal F}(G_{ij}, A^i_\alpha, g_{\alpha \beta}, f) 
\: = \\
& \frac12 \: \int_B \left|
g^{\alpha \beta} \: 
G_{ij; \alpha \beta} \: - \:
g^{\alpha \beta} \: G^{kl} \: G_{ik, \alpha} \:
G_{lj, \beta} \: - \: \frac12 \: g^{\alpha \gamma} \: g^{\beta \delta} \:
G_{ik} \: G_{jl} \: F^k_{\alpha \beta} \: F^l_{\gamma \delta}
\: - \: g^{\alpha \beta} \: G_{ij, \alpha} \: f_{,\beta}
\right|^2 e^{-f} \: \dvol_B \: + \notag \\
& \int_B \left|
g^{\gamma \delta} \:
F^i_{\alpha \gamma; \delta} \: + \: g^{\gamma \delta} \:
G^{ij} \:
G_{jk, \gamma} \: F^k_{\alpha \delta} \: - \: g^{\gamma \delta} \:
f_{, \gamma} \: F^k_{\alpha \delta} \right|^2 \:
e^{-f} \: \dvol_B \: + \notag \\
& 2 \int_B \left|
R_{\alpha \beta} \: - \: \frac14 \:
G^{ij} \: G_{jk,\alpha} \: G^{kl} \: G_{li,\beta} \: - \:
\frac12 g^{\gamma \delta} \: \: G_{ij} \: F^i_{\alpha \gamma} \:
F^j_{\beta \delta} \: + \: f_{;\alpha \beta} \right|^2 
\: e^{-f} \: \dvol_B. \notag
\end{align}
\end{corollary}
\begin{proof}
This is an immediate consequence of Lemma \ref{4.13}.
\end{proof}

As with Perelman's ${\mathcal F}$-functional,
we now perform an infinitesimal diffeomorphism to decouple
the equation for $f$ and obtain the Ricci flow on $M$.

\begin{corollary} \label{4.18}
Under the flow equations
\begin{align} \label{4.19}
\frac{\partial G_{ij}}{\partial t} \: & = \: 
 g^{\alpha \beta} \: 
G_{ij; \alpha \beta} \: - \:
g^{\alpha \beta} \: G^{kl} \: G_{ik, \alpha} \:
G_{lj, \beta} \: - \: \frac12 \: g^{\alpha \gamma} \: g^{\beta \delta} \:
G_{ik} \: G_{jl} \: F^k_{\alpha \beta} \: F^l_{\gamma \delta} \\
\frac{\partial A^i_\alpha}{\partial t} \: & = \: 
- \: g^{\gamma \delta} \:
F^i_{\alpha \gamma; \delta} \: - \: g^{\gamma \delta} \:
G^{ij} \:
G_{jk, \gamma} \: F^k_{\alpha \delta} \notag \\
\frac{\partial g_{\alpha \beta}}{\partial t} \: & = \: - \: 2 \:
R_{\alpha \beta} \: + \: \frac12 \:
G^{ij} \: G_{jk,\alpha} \: G^{kl} \: G_{li,\beta} \: + \: 
g^{\gamma \delta} \: \: G_{ij} \: F^i_{\alpha \gamma} \:
F^j_{\beta \delta} \notag \\
\frac{\partial(e^{-f})}{\partial t} \: & = \:
- \: \nabla^2 \: e^{-f} \: + \:
\left( R \: - \: \frac14 \:
g^{\alpha \beta} \: G^{ij} \: G_{jk, \alpha} \: G^{kl} \:
G_{li, \beta} \: - \: \frac12 \: g^{\alpha \gamma} \: g^{\beta \delta} \:
G_{ij} \: F^i_{\alpha \beta} \: F^j_{\gamma \delta} \right) e^{-f} \notag
\end{align}
one has
\begin{align} \label{4.20}
& \frac{d}{dt} {\mathcal F}(G_{ij}, A^i_\alpha, g_{\alpha \beta}, f) 
\: = \\
& \frac12 \: \int_B \left|
g^{\alpha \beta} \: 
G_{ij; \alpha \beta} \: - \:
g^{\alpha \beta} \: G^{kl} \: G_{ik, \alpha} \:
G_{lj, \beta} \: - \: \frac12 \: g^{\alpha \gamma} \: g^{\beta \delta} \:
G_{ik} \: G_{jl} \: F^k_{\alpha \beta} \: F^l_{\gamma \delta}
\: - \: g^{\alpha \beta} \: G_{ij, \alpha} \: f_{,\beta}
\right|^2 e^{-f} \: \dvol_B \: + \notag \\
& \int_B \left|
g^{\gamma \delta} \:
F^i_{\alpha \gamma; \delta} \: + \: g^{\gamma \delta} \:
G^{ij} \:
G_{jk, \gamma} \: F^k_{\alpha \delta} \: - \: g^{\gamma \delta} \:
f_{, \gamma} \: F^k_{\alpha \delta} \right|^2 \:
e^{-f} \: \dvol_B \: + \notag \\
& 2 \int_B \left|
R_{\alpha \beta} \: - \: \frac14 \:
G^{ij} \: G_{jk,\alpha} \: G^{kl} \: G_{li,\beta} \: - \:
\frac12 g^{\gamma \delta} \: \: G_{ij} \: F^i_{\alpha \gamma} \:
F^j_{\beta \delta} \: + \: f_{;\alpha \beta} \right|^2 
\: e^{-f} \: \dvol_B. \notag
\end{align}
\end{corollary}
\begin{proof}
This follows because the right-hand sides of (\ref{4.16}) and 
(\ref{4.19}) differ
by a Lie derivative with respect to $\nabla f$.
\end{proof}

Note that the first three equations in (\ref{4.19}) are the same as
(\ref{4.10}).

We now analyze what it means for ${\mathcal F}$ to be constant
along the flow (\ref{4.19}).

\begin{proposition} \label{4.21}
If ${\mathcal F}(G_{ij}, A^i_\alpha, 
g_{\alpha \beta}, f)$ is constant in $t$ then
$F^i_{\alpha \beta} = 0$, $\det(G_{ij})$ is constant and
\begin{align} \label{4.22}
g^{\alpha \beta} \: 
G_{ij; \alpha \beta} \: - \:
g^{\alpha \beta} \: G^{kl} \: G_{ik, \alpha} \:
G_{lj, \beta} \: & = \: 0,
\\
R_{\alpha \beta} \: - \: \frac14 \:
G^{ij} \: G_{jk,\alpha} \: G^{kl} \: G_{li,\beta}
& = \: 0. \notag
\end{align}
\end{proposition}
\begin{proof}
From (\ref{4.20}), we have
\begin{equation} \label{4.23}
g^{\alpha \beta} \: 
G_{ij; \alpha \beta} \: - \:
g^{\alpha \beta} \: G^{kl} \: G_{ik, \alpha} \:
G_{lj, \beta} \: - \: \frac12 \: g^{\alpha \gamma} \: g^{\beta \delta} \:
G_{ik} \: G_{jl} \: F^k_{\alpha \beta} \: F^l_{\gamma \delta}
\: - \: g^{\alpha \beta} \: G_{ij, \alpha} \: f_{,\beta} \: = \: 0
\end{equation}
and
\begin{equation} \label{4.24}
R_{\alpha \beta} \: - \: \frac14 \:
G^{ij} \: G_{jk,\alpha} \: G^{kl} \: G_{li,\beta} \: - \:
\frac12 g^{\gamma \delta} \: \: G_{ij} \: F^i_{\alpha \gamma} \:
F^j_{\beta \delta} \: + \: f_{;\alpha \beta} \: = \: 0.
\end{equation}
Multiplying (\ref{4.23}) by $G^{ij}$ and summing over indices gives
\begin{equation} \label{4.25}
\nabla^2 \ln \det(G_{ij}) \: - \: \langle \nabla f, \nabla
\ln \det(G_{ij}) \rangle \: - \: \frac12 \: 
g^{\alpha \gamma} \: g^{\beta \delta} \:
G_{ij} \: F^i_{\alpha \beta} \: F^j_{\gamma \delta} \: = \: 0.
\end{equation}
(Here we are using the trivialization of 
$|\Lambda^{max}e|$ to think of $\det(G_{ij})$ as a function
on $B$, defined up to multiplication by a positive constant.)
Equivalently,
\begin{equation} \label{4.26}
\nabla^\alpha \left( e^{-f} \nabla_\alpha \ln \det(G_{ij}) \right) 
\: - \: \frac12 \: e^{-f} 
g^{\alpha \gamma} \: g^{\beta \delta} \:
G_{ij} \: F^i_{\alpha \beta} \: F^j_{\gamma \delta} \: = \: 0.
\end{equation}
Integrating (\ref{4.26}) over $B$ gives $F^i_{\alpha \beta} = 0$.
Then multiplying (\ref{4.26}) by $\ln \det(G_{ij})$ and integrating over
$B$ gives $\nabla \ln \det(G_{ij}) \: = \: 0$, so $\ln \det(G_{ij})$
is spatially constant.

Given that $F^i_{\alpha \beta} = 0$, the equation for
$G^{ij} \frac{\partial G_{ij}}{\partial t}$ implies
\begin{equation} \label{4.27}
\frac{\partial}{\partial t} \ln \det(G_{ij}) \: = \:
\nabla^2 \ln \det(G_{ij}). 
\end{equation}
Thus $\ln \det(G_{ij})$ is also temporally constant.

As $\det(G_{ij})$ is spatially constant, we have
\begin{equation} \label{4.28}
G^{ij} \: G_{ij;\alpha \beta} \: - \: G^{ij} \: G_{jk, \alpha} \:
G^{kl} \: G_{li, \beta} \: = \: 0.
\end{equation}
Along with the fact that $F^i_{\alpha \beta} = 0$, it follows that
\begin{align} \label{4.29}
\overline{R}_{ij} \: & = \: - \: \frac12 \: g^{\alpha \beta} \: 
G_{ij; \alpha \beta} \: + \:
\frac12 \: g^{\alpha \beta} \: G^{kl} \: G_{ik, \alpha} \:
G_{lj, \beta} \\
\overline{R}_{i \alpha} \: & = \: 0 \notag \\
\overline{R}_{\alpha \beta} \: & =  \: R_{\alpha \beta} \: - \: 
\frac14 \:
G^{ij} \: G_{jk,\alpha} \: G^{kl} \: G_{li,\beta} \notag
\end{align}
and
\begin{equation} \label{4.30}
\overline{R} \: = \: R \: - \: \frac14 \:
g^{\alpha \beta} \: G^{ij} \: G_{jk, \alpha} \: G^{kl} \:
G_{li, \beta}.
\end{equation}
From equation (\ref{4.24}),
\begin{equation} \label{4.31}
\int_B \overline{R} \: \dvol_B \: = \: 0.
\end{equation}

On $M$, the evolution of the scalar curvature is given by
\begin{equation} \label{4.32}
\frac{\partial \overline{R}}{\partial t} \: = \:
\overline{\nabla}^2 \overline{R} \: + \: 2 \: |\overline{R}_{IJ}|^2.
\end{equation}
In our case, and using the fact that $\det(G_{ij})$ is spatially
constant, this becomes
\begin{equation} \label{4.33}
\frac{\partial \overline{R}}{\partial t} \: = \:
{\nabla}^2 \overline{R} \: + \: 2 \: |\overline{R}_{ij}|^2
 \: + \: 2 \: |\overline{R}_{\alpha \beta}|^2.
\end{equation}
From (\ref{4.19}), (\ref{4.23}) and (\ref{4.24}), the flow equations are
\begin{align} \label{4.34}
\frac{\partial G_{ij}}{\partial t} \: & = \: g^{\alpha \beta} \: G_{ij,\alpha} \:
f_{,\beta} \\
\frac{\partial g_{\alpha \beta}}{\partial t} \: & = \: 2 \:
f_{;\alpha \beta}. \notag
\end{align}
As the right-hand side of (\ref{4.34}) 
is given by Lie derivatives with respect
to $\nabla f$, it follows that
\begin{equation} \label{4.35}
\frac{\partial \overline{R}}{\partial t} \: = \: \langle \nabla f, \nabla \overline{R}
\rangle.
\end{equation}

Thus
\begin{equation} \label{4.36}
{\nabla}^2 \overline{R} \: + \: 2 |\overline{R}_{ij}|^2 \: + \:
2 |\overline{R}_{\alpha \beta}|^2 
\: = \: \langle {\nabla} {f}, {\nabla} 
\overline{R} \rangle,
\end{equation}
or
\begin{equation} \label{4.37}
{\nabla}^2 \overline{R} \: + \: 2 |\overline{R}_{ij}|^2 \: + \:
2 |\overline{R}_{\alpha \beta} \: - \: \frac{1}{n} \: 
\overline{R} \: g_{\alpha \beta}|^2 \: + \: \frac{2}{n} \:
\overline{R}^2 
\: = \: \langle {\nabla} {f}, {\nabla} 
\overline{R} \rangle.
\end{equation}
From (\ref{4.31}), either $\overline{R} = 0$ or $\overline{R}_{min} < 0$.
If $\overline{R}_{min} < 0$ then we obtain a contradiction to the
minimum principle, applied to (\ref{4.37}).  Thus
$\overline{R} = 0$. Equation (\ref{4.37}) now implies that
$\overline{R}_{ij} = \overline{R}_{\alpha \beta} = 0$, which proves
the proposition.
\end{proof}

From (\ref{4.19}), under the conclusion of Proposition \ref{4.21} it
follows that $G_{ij}$, 
$A^i_\alpha$ and $g_{\alpha \beta}$ are time-independent.
The Ricci flow solution $\overline{g}_\infty(\cdot)$ on
$M$ is Ricci-flat. In the case $N=0$ the proof of Proposition
\ref{4.21} essentially reduces to the standard proof that a steady gradient
soliton on a compact manifold is Ricci-flat; see, for example,
\cite[Chapter 1]{Chowetal}.

With $\det^{-1}(\pm 1) \subset \GL(N, \R)$, we can write
$\det^{-1}(\pm 1)/\OO(N) = \SL(N, \R)/\SO(N)$. 
From \cite[Proposition 4.17]{Lott (2007)}, the first equation in (\ref{4.22}) 
says that the map $b \rightarrow G_{ij}(b)$ 
describes
a (twisted) harmonic map $G \: : \: B \rightarrow 
\det^{-1}(\pm 1)/\OO(N)$. The
twisting refers to the fact that if the flat $\R^N$-bundle $e$ has
holonomy representation $\rho \: : \: \pi_1(B, b_0) \rightarrow
\det^{-1}(\pm 1)$ then we really have a harmonic map $\widetilde{G} \: : \:
\widetilde{B} \rightarrow  \det^{-1}(\pm 1)/\OO(N)$ which satisfies
$\widetilde{G}(\gamma \widetilde{b}) \: = \: 
\rho(\gamma) \: \widetilde{G}(\widetilde{b})$ for $\gamma \in \pi_1(B, b)$ and
$\widetilde{b}$ in the universal cover $\widetilde{B}$. 
After passing to a double cover of $B$ if necessary, we can assume that
$\rho$ takes value in $\SL(N, \R)$.  
For simplicity,
we will make this assumption
hereafter and consider $\widetilde{G}$ to be a twisted harmonic map
from $B$ to $\SL(N, \R)/\SO(N)$.
Information on such twisted harmonic maps appears in
\cite{Corlette (1988)}, 
\cite[Section 1.2]{Jost-Zuo (1997)} and \cite{Labourie (1991)}.
Given $\rho$, 
such a twisted harmonic map $G$ exists if and only if the
Zariski closure of $\Image(\rho)$ is reductive in $\SL(N, \R)$.
Given $\rho$, if there are two such equivariant harmonic maps 
$\widetilde{G}_1$ and $\widetilde{G}_2$ then there is a 
$1$-parameter family $\{ \widetilde{G}_t \}_{t \in [1,2]}$ of such
equivariant harmonic maps, all with the same quotient energy,
so that for each $\widetilde{b} \in \widetilde{B}$ the map
$t \rightarrow \widetilde{G}_t(\widetilde{b})$ is a 
constant-speed geodesic arc, whose length is independent of 
$\widetilde{b}$.

If the second equation in (\ref{4.22}) is satisfied then $B$ clearly has
nonnegative Ricci curvature.

We now look at the solutions of (\ref{4.22}).

\begin{proposition} \label{4.38}
Any solution $\overline{g}$ of (\ref{4.22}) is a locally product 
metric on a Ricci-flat base $B$.
\end{proposition}
\begin{proof}
From the second equation in (\ref{4.22}), $B$ has nonnegative Ricci
curvature. For some $r$, the universal cover $\widetilde{B}$ is an isometric
product of $\R^r$ and $W$, where $W$ is a simply-connected closed
$(n-r)$-dimensional manifold of nonnegative Ricci curvature
\cite{Cheeger-Gromoll (1971)}.
As before, let $\widetilde{G} \: : \: \widetilde{B} \rightarrow
\SL(N, \R)/\SO(N)$ denote the lift of
$G$ to $\widetilde{B}$.

Let $x^1, \ldots, x^r$ be Cartesian coordinates on $\R^r$ and let
$x^{r+1}, \ldots, x^n$ be local coordinates on $W$. 
From the second equation of (\ref{4.22}), 
$\widetilde{G}_{ij, \alpha} = 0$ for
$1 \le \alpha \le r$. That is, $\widetilde{G}$ is constant in the
$\R^k$-directions. Then the first equation of (\ref{4.22}) implies that
for each $y \in \R^r$, the restriction of $\widetilde{G}$
to $\{y\} \times W$ is a harmonic
map from $W$ to $\SL(N, \R)/\SO(N)$. It follows that 
for each $y \in \R^r$, the restriction of $\widetilde{G}$
to $\{y\} \times W$ is a point map.  Thus $\widetilde{G}$ is constant.
From the second equation of (\ref{4.22}), $\widetilde{B}$ is
Ricci flat. The conclusion is that
$B$ is Ricci flat and $G$ is locally constant.
\end{proof}

In the next proposition we use ${\mathcal F}$ to analyze a
long-time limit of a locally ${\mathcal G}$-invariant Ricci flow
solution.  The method of proof is along the lines of the proof of
\cite[Theorem 1.3]{Feldman-Ilmanen-Ni (2005)}.

\begin{proposition} \label{4.39}
Suppose that $(M, \overline{g}(\cdot))$ is a locally ${\mathcal G}$-invariant Ricci flow
defined for all $t \in [0, \infty)$.
Let $\{ s_i \}_{i=1}^\infty$ be a sequence of positive numbers
tending to
infinity.  Put $\overline{g}_i(t) \: = \: \overline{g}(t+s_i)$. Suppose
that $\lim_{i \rightarrow \infty} \overline{g}_i(\cdot)$ exists and equals
$\overline{g}_\infty(\cdot)$ in the sense of Subsection \ref{subsection4.1}, 
for
a locally ${\mathcal G}$-invariant Ricci flow $\overline{g}_\infty(\cdot)$ with a compact base
$B_\infty$. Writing $\overline{g}_\infty(\cdot) \: \equiv \: 
(G_{ij,\infty}(\cdot), A^i_{\alpha,\infty}(\cdot), 
g_{\alpha \beta,\infty}(\cdot))$, we conclude that \\
1. The curvatures $F^i_{\alpha \beta, \infty}$ vanish. \\
2. $\det(G_{ij,\infty})$ is
constant. \\
3. Equations (\ref{4.22}) are satisfied for 
$G_{ij,\infty}(\cdot)$ and $g_{\alpha \beta,\infty}(\cdot)$.
\end{proposition}
\begin{proof}
We first construct a positive solution of the conjugate heat equation
\begin{equation} \label{4.40}
\frac{\partial u}{\partial t} \: = \:
- \: \nabla^2 u \: + \:
\left( R \: - \: \frac14 \:
g^{\alpha \beta} \: G^{ij} \: G_{jk, \alpha} \: G^{kl} \:
G_{li, \beta} \: - \: \frac12 \: g^{\alpha \gamma} \: g^{\beta \delta} \:
G_{ij} \: F^i_{\alpha \beta} \: F^j_{\gamma \delta} \right) u.
\end{equation}
that exists for all $t \in [0, \infty)$. Note that if $u$ is a 
solution to (\ref{4.40}) then $\int_B u \: \dvol_B$ is constant in $t$.
Let $\{t_j\}_{j=1}^\infty$ be a sequence of times going to infinity.
Let $\widetilde{u}_j(\cdot)$ 
be a solution to (\ref{4.40}) on the interval $[0,t_j]$ with initial
condition 
$\widetilde{u}_j(t_j) \: = \: \frac{1}{\vol(B, g_{\alpha \beta}(t_j))}$.
For any $T > 0$, we claim that a subsequence of the $\widetilde{u}_j$'s 
converges smoothly on the time interval $[0, T]$. 
To see this, at time $T+1$ we know
that if $t_j \ge T+1$ then
$\widetilde{u}_j(T+1) \ge 0$ and 
$\int_B \widetilde{u}_j(T+1) \: \dvol_B \: = \: 1$. Solving
the conjugate heat equation with initial data at time $T+1$,
and restricting the solution to the time interval $[0,T]$, gives a 
smoothing operator from
the space of initial data 
$\{\widetilde{u} \in L^1(B) \: : \: \widetilde{u} \ge 0, 
\int_B \widetilde{u} \: \dvol_B(T+1) 
= 1\}$ to $C^\infty([0, T] \times B)$. 
Thus we have the derivative bounds needed to extract a 
subsequence of the $\widetilde{u}_j$'s that converges smoothly on $[0, T]$. 
By a diagonal
argument, we can extract a subsequence
of the $\widetilde{u}_j$'s that converges smoothly on
compact subsets of $[0, \infty)$ 
to a nonzero solution $\widetilde{u}_\infty(\cdot)$ of (\ref{4.40}), defined
for $t \in [0, \infty)$. 

One can show, as in \cite[Pf. of Proposition 7.5]{Kleiner-Lott}, 
that $\widetilde{u}_\infty(\cdot)
> 0$.
If $\widetilde{f}_\infty(t)$ is given by
$\widetilde{u}_\infty(t) \: = \: e^{- \: 
\widetilde{f}_\infty(t)}$ then
${\mathcal F}(G_{ij}(t), A^i_\alpha(t), g_{\alpha \beta}(t), 
\widetilde{f}_\infty(t))$ is
nondecreasing in $t$. We write ${\mathcal F}_\infty \: = \:
\lim_{t \rightarrow \infty}
{\mathcal F}(G_{ij}(t), A^i_\alpha(t), g_{\alpha \beta}(t), 
\widetilde{f}_\infty(t))$, which is
possibly infinite for the moment.

Next, put $u_i(t) \: = \: \widetilde{u}_\infty(t + s_i)$. By assumption,
$\lim_{i \rightarrow \infty} \overline{g}_i(\cdot) \: = \:
\overline{g}_\infty(\cdot)$ in the sense of Subsection \ref{subsection4.1}.
Then by the same smoothing argument as above, there is a subsequence
of $\{u_i(\cdot)\}_{i=1}^\infty$ that converges smoothly on compact
subsets of $[0, \infty)$ to a solution $u_\infty(\cdot)$ of (\ref{4.40}) on
$B_\infty$, where (\ref{4.40}) is now written in terms of
$G_{ij,\infty}(\cdot)$, $A^i_{\alpha,\infty}(\cdot)$ and
$g_{\alpha \beta,\infty}(\cdot)$. 
(When taking a convergent subsequence, 
we perform the same diffeomorphisms on the $u_i$'s as are used
in forming the limit $\lim_{i \rightarrow \infty} \overline{g}_i(\cdot)$.)
Define $f_\infty(t)$ by
${u}_\infty(t) \: = \: e^{- \: {f}_\infty(t)}$.
Then after passing to a subsequence,
\begin{align} \label{4.41}
{\mathcal F}(G_{ij,\infty}(t), A^i_{\alpha,\infty}(t), 
g_{\alpha \beta,\infty}(t), 
f_\infty(t)) \: & = \:
\lim_{i \rightarrow \infty}
{\mathcal F}(G_{ij}(t+s_i), A^i_\alpha(t+s_i), 
g_{\alpha \beta}(t+s_i), \widetilde{f}_\infty(t+s_i)) \\
& = \: {\mathcal F}_\infty. \notag
\end{align}
This shows that ${\mathcal F}_\infty < \infty$ and that
${\mathcal F}(G_{ij,\infty}(t), A^i_{\alpha,\infty}(t), 
g_{\alpha \beta,\infty}(t), 
f_\infty(t))$ is constant in $t$. The proposition now follows from
Proposition \ref{4.21}.
\end{proof}

Junfang Li pointed out that the modified ${\mathcal F}$-functional
has an $(n+N)$-dimensional interpretation.  Namely, for
$\overline{f} \in C^\infty(B)$,
put
\begin{equation} \label{4.42}
\overline{\mathcal F}(G_{ij}, A^i_\alpha, 
g_{\alpha \beta}, \overline{f}) \: = \:
\int_B \left( |\nabla \overline{f}|^2 \: + \: \overline{R} \right)
\: e^{- \overline{f}} \: \sqrt{\det(G_{ij})} \: \dvol_B.
\end{equation}
This is a renormalized version of Perelman's ${\mathcal F}$-functional
on $M$. 
\begin{proposition} \label{4.43}
Put $f = \overline{f} - \ln \sqrt{\det(G_{ij})}$. Then
\begin{equation} \label{4.44}
\overline{\mathcal F}(G_{ij}, A^i_\alpha, g_{\alpha \beta}, \overline{f}) \: = 
{\mathcal F}(G_{ij}, A^i_\alpha, g_{\alpha \beta},f).
\end{equation}
\end{proposition}
\begin{proof}
We have
\begin{equation} \label{4.45}
\int_B |\nabla \overline{f} |^2
\: e^{- \overline{f}} \: \sqrt{\det(G_{ij})} \: \dvol_B \: = \\
\int_B \left| \nabla f \: + \: \nabla \ln \sqrt{\det(G_{ij})} 
\right|^2
\: e^{- f} \: \dvol_B
\end{equation}
and
\begin{align} \label{4.46}
& \int_B \left| \nabla f \: + \: \nabla \ln \sqrt{\det(G_{ij})} 
\right|^2
\: e^{- f} \: \dvol_B \: = \\
&\int_B \left( \left| \nabla f \right|^2 \: + \: 
2 \: \langle \nabla f, \nabla \ln \sqrt{\det(G_{ij})} \rangle
\: + \: \left| \nabla \ln \sqrt{\det(G_{ij})} \right|^2 \right)
\: e^{- f} \: \dvol_B \: = \: \notag \\ 
&\int_B \left( \left| \nabla f \right|^2 \: + \: 2 \: 
\nabla^2 \ln \sqrt{\det(G_{ij})}
\: + \: \left| \nabla \ln \sqrt{\det(G_{ij})} \right|^2 \right)
\: e^{- f} \: \dvol_B \: = \: \notag \\ 
&\int_B \left( \left| \nabla f \right|^2 \: + \: 
g^{\alpha \beta} G^{ij}  \: G_{ij; \alpha \beta} \: - \: 
g^{\alpha \beta} \: G^{ij} \: G_{jk, \alpha} \: G^{kl} \:
G_{li, \beta} \: + \: \frac14 \: g^{\alpha \beta} \: G^{ij} \: 
G_{ij, \alpha} \: G^{kl} \: G_{kl, \beta} \right) \: 
e^{- f} \: \dvol_B. \notag
\end{align}
Combining this with (\ref{4.7}) gives
\begin{align} \label{4.46.5}
& \int_B \left( \left| \nabla f \: + \: \nabla \ln \sqrt{\det(G_{ij})} 
\right|^2 \: + \: \overline{R} \right)
\: e^{- f} \: \dvol_B \: = \\
& \int_B \left( |\nabla f|^2 \: + \: R \: - \: \frac14 \:
g^{\alpha \beta} \: G^{ij} \: G_{jk, \alpha} \: G^{kl} \:
G_{li, \beta} \: - \: \frac14 \: g^{\alpha \gamma} \: g^{\beta \delta} \:
G_{ij} \: F^i_{\alpha \beta} \: F^j_{\gamma \delta} \right)
\: e^{-f} \: \dvol_B, \notag
\end{align}
which proves the proposition.
\end{proof}

\subsubsection{Modified ${\mathcal W}$-functional} \label{subsubsection4.2.2}

\begin{definition} \label{4.47}
Given $f \in C^\infty(B)$ and $\tau \in \R^+$, put
\begin{align} \label{4.48}
& {\mathcal W}(G_{ij},A^i_\alpha,g_{\alpha \beta},f,\tau) \: = \\
&\int_B \left[ \tau
\left( |\nabla f|^2 \: + \: 
R \: - \: \frac14 \:
g^{\alpha \beta} \: G^{ij} \: G_{jk, \alpha} \: G^{kl} \:
G_{li, \beta} \: - \: \frac14 \: g^{\alpha \gamma} \: g^{\beta \delta} \:
G_{ij} \: F^i_{\alpha \beta} \: F^j_{\gamma \delta} \right)
\: + \: f \: - \: n \right] \notag \\
& (4\pi \tau)^{- \: \frac{n}{2}}
\: e^{-f} \: \dvol_B. \notag
\end{align}
\end{definition}

If $N = 0$, i.e. if $M = B$, then this is the same as
Perelman's ${\mathcal W}$-functional \cite{Perelman1}.
The next proposition says how ${\mathcal W}$ varies along the
Ricci flow.

\begin{proposition} \label{4.49}
Under the flow equations
\begin{align} \label{4.50}
\frac{\partial G_{ij}}{\partial t} \: & = \: 
 g^{\alpha \beta} \: 
G_{ij; \alpha \beta} \: - \:
g^{\alpha \beta} \: G^{kl} \: G_{ik, \alpha} \:
G_{lj, \beta} \: - \: \frac12 \: g^{\alpha \gamma} \: g^{\beta \delta} \:
G_{ik} \: G_{jl} \: F^k_{\alpha \beta} \: F^l_{\gamma \delta} \\
\frac{\partial A^i_\alpha}{\partial t} \: & = \: 
- \: g^{\gamma \delta} \:
F^i_{\alpha \gamma; \delta} \: - \: g^{\gamma \delta} \:
G^{ij} \:
G_{jk, \gamma} \: F^k_{\alpha \delta} \notag \\
\frac{\partial g_{\alpha \beta}}{\partial t} \: & = \: - \: 2 \:
R_{\alpha \beta} \: + \: \frac12 \:
G^{ij} \: G_{jk,\alpha} \: G^{kl} \: G_{li,\beta} \: + \: 
g^{\gamma \delta} \: \: G_{ij} \: F^i_{\alpha \gamma} \:
F^j_{\beta \delta} \notag \\
\frac{\partial(e^{-f})}{\partial t} \: & = \:
- \: \nabla^2 \: e^{-f} \: + \:
\left( R \: - \: \frac14 \:
g^{\alpha \beta} \: G^{ij} \: G_{jk, \alpha} \: G^{kl} \:
G_{li, \beta} \: - \: \frac12 \: g^{\alpha \gamma} \: g^{\beta \delta} \:
G_{ij} \: F^i_{\alpha \beta} \: F^j_{\gamma \delta} \: - \:
\frac{n}{2 \tau} \right) e^{-f} \notag \\
\frac{\partial \tau}{\partial t} \: & = \: -1 \notag
\end{align}
one has
\begin{align} \label{4.51}
& \frac{d}{dt} {\mathcal W}(G_{ij}, A^i_\alpha, g_{\alpha \beta}, f, \tau) 
\: = \\
& \frac{\tau}{2} \: \int_B \left|
g^{\alpha \beta} \: 
G_{ij; \alpha \beta} \: - \:
g^{\alpha \beta} \: G^{kl} \: G_{ik, \alpha} \:
G_{lj, \beta} \: - \: \frac12 \: g^{\alpha \gamma} \: g^{\beta \delta} \:
G_{ik} \: G_{jl} \: F^k_{\alpha \beta} \: F^l_{\gamma \delta}
\: - \: g^{\alpha \beta} \: G_{ij, \alpha} \: f_{,\beta}
\right|^2 \notag \\
& (4 \pi \tau)^{- \: \frac{n}{2}} \: 
e^{-f} \: \dvol_B \: + \notag \\
& \tau \int_B \left|
g^{\gamma \delta} \:
F^i_{\alpha \gamma; \delta} \: + \: g^{\gamma \delta} \:
G^{ij} \:
G_{jk, \gamma} \: F^k_{\alpha \delta} \: - \: g^{\gamma \delta} \:
f_{, \gamma} \: F^k_{\alpha \delta} \right|^2 \:
 (4 \pi \tau)^{- \: \frac{n}{2}} \:
e^{-f} \: \dvol_B \: + \notag \\
& 2 \tau \int_B \left|
R_{\alpha \beta} \: - \: \frac14 \:
G^{ij} \: G_{jk,\alpha} \: G^{kl} \: G_{li,\beta} \: - \:
\frac12 \: g^{\gamma \delta} \: \: G_{ij} \: F^i_{\alpha \gamma} \:
F^j_{\beta \delta} \: + \: f_{;\alpha \beta} \: - \: \frac{1}{2\tau}
\: g_{\alpha \beta} \right|^2 
\: (4 \pi \tau)^{- \: \frac{n}{2}} \: e^{-f} \: \dvol_B \: - \: \notag \\
& \frac14 \: \int_B g^{\alpha \gamma} \: g^{\beta \delta} \: G_{ij} \:
F^i_{\alpha \beta} \: F^j_{\gamma \delta} \:
 (4 \pi \tau)^{- \: \frac{n}{2}} \:
e^{-f} \: \dvol_B. \notag
\end{align}
\end{proposition}
\begin{proof}
The proof stands in relation to the proof of Corollary \ref{4.18} as
the corresponding statements about Perelman's ${\mathcal W}$-functional
vs. Perelman's ${\mathcal F}$-functional; see
\cite[Section 12]{Kleiner-Lott}.
\end{proof}

Note that the 
$\frac14 \: \int_B g^{\alpha \gamma} \: g^{\beta \delta} \: G_{ij} \:
F^i_{\alpha \beta} \: F^j_{\gamma \delta} \:
 (4 \pi \tau)^{- \: \frac{n}{2}} \:
e^{-f} \: \dvol_B$ term occurs on the right-hand
side of (\ref{4.51}) with a negative sign. We now look at what it
means for ${\mathcal W}$ to be constant in $t$, under the
assumption that $F^i_{\alpha \beta}$ vanishes.

\begin{proposition} \label{4.52}
Suppose that $F^i_{\alpha \beta} = 0$. If 
${\mathcal W}(G_{ij}, A^i_\alpha, g_{\alpha \beta}, f, \tau)$ is constant
in $t$ then $\det(G_{ij})$ is constant and
\begin{align} \label{4.53}
g^{\alpha \beta} \: 
G_{ij; \alpha \beta} \: - \:
g^{\alpha \beta} \: G^{kl} \: G_{ik, \alpha} \:
G_{lj, \beta} \: - \: g^{\alpha \beta} \: G_{ij, \alpha} \: f_{,\beta}
& = \: 0, \\
R_{\alpha \beta} \: - \: \frac14 \:
G^{ij} \: G_{jk,\alpha} \: G^{kl} \: G_{li,\beta} \: + \:
f_{;\alpha \beta} \: - \: \frac{1}{2\tau} \: g_{\alpha \beta}
& = \: 0. \notag
\end{align}
\end{proposition}
\begin{proof}
The same argument as in the proof of Proposition \ref{4.21} shows that
$\det(G_{ij})$ is constant. Then (\ref{4.53}) follows from 
(\ref{4.51}). 

Unlike in Proposition \ref{4.21}, we cannot conclude that $f$ is constant,
because of the existence of nontrivial compact gradient shrinking solitons.
\end{proof}

\begin{remark} \label{4.54}
The term $\frac14 \: \int_B g^{\alpha \gamma} \: g^{\beta \delta} \: G_{ij} \:
F^i_{\alpha \beta} \: F^j_{\gamma \delta} \:
 (4 \pi \tau)^{- \: \frac{n}{2}} \:
e^{-f} \: \dvol_B$ occurs on the right-hand
side of (\ref{4.51}) with a useless sign.  This is not surprising, as can
be seen by looking at the Ricci flow on a round $3$-sphere $M$,
which we consider to be the total space of a circle bundle over
$S^2$. We shift the time parameter 
so that the $3$-sphere disappears at time zero.
As the $3$-sphere gives a gradient shrinking soliton,
the functional ${\mathcal W}$ is constant in $t$. However,
the circle bundle has nonvanishing curvature.  Hence having
${\mathcal W}$ constant in $t$ cannot imply that $F^i_{\alpha \beta}$ 
vanishes.
\end{remark}

We now look at some special cases of (\ref{4.53}).

\begin{proposition}
Under the hypotheses of Proposition \ref{4.52}, if $1 \le \dim(B) \le 2$
then the only solutions of (\ref{4.53}) occur when $B$ is
$S^2$ or $\R P^2$.
\end{proposition}
\begin{proof}
The second equation in (\ref{4.53}) implies that
\begin{equation}
\int_B R \: \dvol_B \: - \: \frac14 \: \int_B g^{\alpha \beta} \: 
G^{ij} \: G_{jk,\alpha} \: G^{kl} \: G_{li,\beta} \: \dvol_B 
\: - \: \frac{n}{2\tau} \: \vol(B) \: = \: 0,
\end{equation}
from which the proposition follows.
\end{proof}

We now use ${\mathcal W}$ to analyze a blowup limit.

\begin{proposition} \label{4.55}
Suppose that $(M, \overline{g}(\cdot))$ is a locally ${\mathcal G}$-invariant Ricci flow
defined for all $t \in (- T, 0)$, with $T \le \infty$.
Suppose that $F^i_{\alpha \beta} = 0$. 
Put $\tau \: = \: - \: t$.
Let $\{ s_i \}_{i=1}^\infty$ be a sequence of positive numbers
tending to infinity.
Put $\overline{g}_i(\tau) \: = \: s_i \: \overline{g}(s_i^{-1} \tau)$. Suppose
that $\lim_{i \rightarrow \infty} \overline{g}_i(\cdot)$ exists and equals
$\overline{g}_\infty(\cdot)$ in the sense of Subsection \ref{subsection4.1}, 
for
a locally ${\mathcal G}$-invariant Ricci flow $\overline{g}_\infty(\cdot)$ with a compact base
$B_\infty$, defined for $\tau \in (0, \infty)$.
Writing $\overline{g}_\infty(\cdot) \: \equiv \: 
(G_{ij,\infty}(\cdot), g_{\alpha \beta,\infty}(\cdot))$, we
conclude that \\
1. $\det(G_{ij,\infty})$ is
constant. \\
2. Equations (\ref{4.53}) are satisfied for 
$G_{ij,\infty}(\cdot)$ and $g_{\alpha \beta,\infty}(\cdot)$.
\end{proposition}
\begin{proof}
The proof is along the lines of the proof of Proposition \ref{4.39}.
\end{proof}

\subsubsection{Modified ${\mathcal W}_+$-functional} 
\label{subsubsection4.2.3}

\begin{definition} \label{4.56}
Given $f \in C^\infty(B)$ and $t \in \R^+$, put
\begin{align} \label{4.57}
& {\mathcal W}_+(G_{ij},A^i_\alpha,g_{\alpha \beta},f,t) \: = \\
&\int_B \left[ t
\left( |\nabla f|^2 \: + \: 
R \: - \: \frac14 \:
g^{\alpha \beta} \: G^{ij} \: G_{jk, \alpha} \: G^{kl} \:
G_{li, \beta} \: - \: \frac14 \: g^{\alpha \gamma} \: g^{\beta \delta} \:
G_{ij} \: F^i_{\alpha \beta} \: F^j_{\gamma \delta} \right)
\: - \: f \: + \: n \right] \notag \\
& (4\pi t)^{- \: \frac{n}{2}}
\: e^{-f} \: \dvol_B. \notag
\end{align}
\end{definition}

If $N = 0$, i.e. if $M = B$, then this is the same as the
Feldman-Ilmanen-Ni ${\mathcal W}_+$-functional 
\cite{Feldman-Ilmanen-Ni (2005)}.

In what follows, 
we will need a lower bound for ${\mathcal W}_+$ in terms of the
scalar curvature of $M$ and the volume of $B$.

\begin{lemma} \label{4.58}
If $(4 \pi t)^{- \: \frac{n}{2}} \int_B e^{-f} \: \dvol_B \: = \: 1$
then 
\begin{equation} \label{4.59}
{\mathcal W}_+(G_{ij},A^i_\alpha,g_{\alpha \beta},f,t) \: \ge \:
t \: \overline{R}_{min} \: + \: n \: + \: \frac{n}{2} \: \ln(4\pi)
\: - \: \ln \left( t^{- \: \frac{n}{2}} \: \vol(B, g_{\alpha \beta}(t)) 
\right). 
\end{equation}
\end{lemma}
\begin{proof}
From (\ref{4.46.5}),
\begin{align} \label{4.60}
& {\mathcal W}_+(G_{ij},A^i_\alpha,g_{\alpha \beta},f,t) \: = \\ 
& \int_B \left[t \left( \left| \nabla f \: + \: \nabla \ln \sqrt{\det(G_{ij})} 
\right|^2 \: + \: \overline{R} \right) \: - \: f \: + \: n \right] \:
(4\pi t)^{- \: \frac{n}{2}} \: e^{- f} \: \dvol_B \: \ge \notag \\
& t \: \overline{R}_{min} \: + \: n \: - \: 
(4 \pi t)^{- \: \frac{n}{2}} \int_B f \: e^{-f} \: \dvol_B \: \ge \: \notag \\
& t \: \overline{R}_{min} \: + \: n \: + \: \frac{n}{2} \: \ln(4\pi)
\: - \: \ln \left( t^{- \: \frac{n}{2}} \: \vol(B, g_{\alpha \beta}(t))
\right), \notag
\end{align}
where we used Jensen's inequality. This proves the lemma.
\end{proof}

The next proposition says that if $f$ satisfies a conjugate
heat equation then ${\mathcal W}_+$ is monotonic under the Ricci flow.

\begin{proposition} \label{4.61}
Under the flow equations
\begin{align} \label{4.62}
\frac{\partial G_{ij}}{\partial t} \: & = \: 
 g^{\alpha \beta} \: 
G_{ij; \alpha \beta} \: - \:
g^{\alpha \beta} \: G^{kl} \: G_{ik, \alpha} \:
G_{lj, \beta} \: - \: \frac12 \: g^{\alpha \gamma} \: g^{\beta \delta} \:
G_{ik} \: G_{jl} \: F^k_{\alpha \beta} \: F^l_{\gamma \delta} \\
\frac{\partial A^i_\alpha}{\partial t} \: & = \: 
- \: g^{\gamma \delta} \:
F^i_{\alpha \gamma; \delta} \: - \: g^{\gamma \delta} \:
G^{ij} \:
G_{jk, \gamma} \: F^k_{\alpha \delta} \notag \\
\frac{\partial g_{\alpha \beta}}{\partial t} \: & = \: - \: 2 \:
R_{\alpha \beta} \: + \: \frac12 \:
G^{ij} \: G_{jk,\alpha} \: G^{kl} \: G_{li,\beta} \: + \: 
g^{\gamma \delta} \: \: G_{ij} \: F^i_{\alpha \gamma} \:
F^j_{\beta \delta} \notag \\
\frac{\partial(e^{-f})}{\partial t} \: & = \:
- \: \nabla^2 \: e^{-f} \: + \:
\left( R \: - \: \frac14 \:
g^{\alpha \beta} \: G^{ij} \: G_{jk, \alpha} \: G^{kl} \:
G_{li, \beta} \: - \: \frac12 \: g^{\alpha \gamma} \: g^{\beta \delta} \:
G_{ij} \: F^i_{\alpha \beta} \: F^j_{\gamma \delta}
\: + \: \frac{n}{2t} \right) e^{-f} \notag
\end{align}
one has
\begin{align} \label{4.63}
& \frac{d}{dt} {\mathcal W}_+(G_{ij}, A^i_\alpha, g_{\alpha \beta}, f, t) 
\: = \\
& \frac{t}{2} \: \int_B \left|
g^{\alpha \beta} \: 
G_{ij; \alpha \beta} \: - \:
g^{\alpha \beta} \: G^{kl} \: G_{ik, \alpha} \:
G_{lj, \beta} \: - \: \frac12 \: g^{\alpha \gamma} \: g^{\beta \delta} \:
G_{ik} \: G_{jl} \: F^k_{\alpha \beta} \: F^l_{\gamma \delta}
\: - \: g^{\alpha \beta} \: G_{ij, \alpha} \: f_{,\beta}
\right|^2 \notag \\
& (4 \pi t)^{- \: \frac{n}{2}} \: 
e^{-f} \: \dvol_B \: + \notag \\
& t \int_B \left|
g^{\gamma \delta} \:
F^i_{\alpha \gamma; \delta} \: + \: g^{\gamma \delta} \:
G^{ij} \:
G_{jk, \gamma} \: F^k_{\alpha \delta} \: - \: g^{\gamma \delta} \:
f_{, \gamma} \: F^k_{\alpha \delta} \right|^2 \:
 (4 \pi t)^{- \: \frac{n}{2}} \:
e^{-f} \: \dvol_B \: + \notag \\
& 2 t \int_B \left|
R_{\alpha \beta} \: - \: \frac14 \:
G^{ij} \: G_{jk,\alpha} \: G^{kl} \: G_{li,\beta} \: - \:
\frac12 \: g^{\gamma \delta} \: \: G_{ij} \: F^i_{\alpha \gamma} \:
F^j_{\beta \delta} \: + \: f_{;\alpha \beta} \: + \: \frac{1}{2t}
\: g_{\alpha \beta} \right|^2 
\: (4 \pi t)^{- \: \frac{n}{2}} \: e^{-f} \: \dvol_B \: + \: \notag \\
& \frac14 \: \int_B g^{\alpha \gamma} \: g^{\beta \delta} \: G_{ij} \:
F^i_{\alpha \beta} \: F^j_{\gamma \delta}
\: (4 \pi t)^{- \: \frac{n}{2}} \: e^{-f} \: \dvol_B. \notag
\end{align}
\end{proposition}
\begin{proof}
The proof is along the lines of the proof of Corollary \ref{4.18}.
\end{proof}

We now look at what it means for ${\mathcal W}_+$ to be constant
along the flow (\ref{4.62}).

\begin{proposition} \label{4.64}
If ${\mathcal W_+}(G_{ij}, A^i_\alpha, g_{\alpha \beta}, f, t)$ 
is constant in $t$ then
$F^i_{\alpha \beta} = 0$, $\det(G_{ij})$ is constant and
\begin{align} \label{4.65}
g^{\alpha \beta} \: 
G_{ij; \alpha \beta} \: - \:
g^{\alpha \beta} \: G^{kl} \: G_{ik, \alpha} \:
G_{lj, \beta} \: & = \: 0, \\
R_{\alpha \beta} \: - \: \frac14 \:
G^{ij} \: G_{jk,\alpha} \: G^{kl} \: G_{li,\beta} \: + \: 
\frac{1}{2t} \: g_{\alpha \beta}
& = \: 0. \notag
\end{align}
\end{proposition}
\begin{proof}
From (\ref{4.63}), we see first that $F^i_{\alpha \beta} = 0$. Then
we also see that
\begin{equation} \label{4.66}
g^{\alpha \beta} \: 
G_{ij; \alpha \beta} \: - \:
g^{\alpha \beta} \: G^{kl} \: G_{ik, \alpha} \:
G_{lj, \beta}
\: - \: g^{\alpha \beta} \: G_{ij, \alpha} \: f_{,\beta} \: = \: 0
\end{equation}
and
\begin{equation} \label{4.67}
R_{\alpha \beta} \: - \: \frac14 \:
G^{ij} \: G_{jk,\alpha} \: G^{kl} \: G_{li,\beta} \:
+ \: f_{;\alpha \beta} \: + \: \frac{1}{2t} \: g_{\alpha \beta} \: = \: 0.
\end{equation}
As in the proof of Proposition \ref{4.21}, we can show from (\ref{4.66}) that
$\det(G_{ij})$ is constant. Then equations (\ref{4.29}) and (\ref{4.30}) hold.
From (\ref{4.67}), we have 
\begin{equation} \label{4.68}
\int_B \left( \overline{R} + \frac{n}{2t} \right)
\: \dvol_B \: = \: 0.
\end{equation}
As in the proof of Proposition \ref{4.21}, we have
\begin{equation} \label{4.69}
\frac{\partial \overline{R}}{\partial t} \: = \:
{\nabla}^2 \overline{R} \: + \: 2 \: |\overline{R}_{ij}|^2
 \: + \: 2 \: |\overline{R}_{\alpha \beta}|^2.
\end{equation}
From (\ref{4.62}), (\ref{4.66}) and (\ref{4.67}), the flow equations are
\begin{align} \label{4.70}
\frac{\partial G_{ij}}{\partial t} \: & = \: g^{\alpha \beta} \: G_{ij,\alpha} \:
f_{,\beta} \\
\frac{\partial g_{\alpha \beta}}{\partial t} \: & = \: 2 \:
f_{;\alpha \beta} \: + \: \frac{1}{t} \: g_{\alpha \beta}. \notag
\end{align}
It follows that
\begin{equation} \label{4.71}
\frac{\partial \overline{R}}{\partial t} \: = \: \langle \nabla f, \nabla \overline{R}
\rangle \: - \: \frac{\overline{R}}{t}.
\end{equation}
Thus
\begin{equation} \label{4.72}
{\nabla}^2 \overline{R} \: + \: 2 |\overline{R}_{ij}|^2 \: + \:
2 |\overline{R}_{\alpha \beta}|^2 \: + \: \frac{\overline{R}}{t}
\: = \: \langle {\nabla} {f}, {\nabla} 
\overline{R} \rangle.
\end{equation}
Then
\begin{equation} \label{4.73}
{\nabla}^2 
\left( \overline{R} + \frac{n}{2t} \right) 
\: + \: 2 |\overline{R}_{ij}|^2 \: + \:
2 |\overline{R}_{\alpha \beta} + \frac{1}{2t} g_{\alpha \beta}|^2 
\: - \: \frac{1}{t} \: \left( \overline{R} + \frac{n}{2t} \right) 
\: = \: \langle {\nabla} {f}, {\nabla} 
\left( \overline{R} + \frac{n}{2t} \right)  \rangle.
\end{equation}
From (\ref{4.68}), either $\overline{R} + \frac{n}{2t} \: = \: 0$ or
$\overline{R}_{min} + \frac{n}{2t} < 0$. If
$\overline{R}_{min} + \frac{n}{2t} < 0$ then we obtain a contradiction
to the minimum principle, applied to (\ref{4.73}). Thus
$\overline{R} + \frac{n}{2t} = 0$. From (\ref{4.73}), it follows that
$\overline{R}_{ij} \: = \: 
\overline{R}_{\alpha \beta} + \frac{1}{2t} g_{\alpha \beta}
\: = \: 0$. This proves the proposition.
\end{proof}

\begin{lemma} \label{4.74}
Under the conclusion of Proposition \ref{4.64}, 
$G_{ij}$ and $A^i_\alpha$ are time-independent,
and $g_{\alpha \beta}$ is proportionate to $t$.
\end{lemma}
\begin{proof}
This follows from (\ref{4.62}) and (\ref{4.65}).
\end{proof}

\begin{remark} \label{4.75}
Equations (\ref{4.65}) were called the harmonic-Einstein equations in
\cite{Lott (2007)}, where they were used as an ansatz to construct
expanding soliton solutions on the total spaces of flat vector
bundles.
\end{remark}

We now use ${\mathcal W}_+$ to analyze blowdown limits.

\begin{proposition} \label{4.76}
Suppose that $(M, \overline{g}(\cdot))$ is a locally ${\mathcal G}$-invariant Ricci flow
defined for all $t \in (0, \infty)$.
Let $\{ s_i \}_{i=1}^\infty$ be a sequence of positive numbers
tending to infinity.
Put $\overline{g}_i(t) \: = \: s_i^{-1} \: 
\overline{g}(s_i t)$. Suppose
that $\lim_{i \rightarrow \infty} \overline{g}_i(\cdot)$ exists and equals
$\overline{g}_\infty(\cdot)$ in the sense of Subsection \ref{subsection4.1}, 
for
a locally ${\mathcal G}$-invariant Ricci flow $\overline{g}_\infty(\cdot)$ with a compact base
$B_\infty$, defined for $t \in (0, \infty)$.
Writing $\overline{g}_\infty(\cdot) \: \equiv \: 
(G_{ij,\infty}(\cdot), A^i_{\alpha,\infty}(\cdot), 
g_{\alpha \beta,\infty}(\cdot))$,
we conclude that \\
1. $F^i_{\alpha \beta,\infty} \: = \: 0$. \\ 
2. $\det(G_{ij,\infty})$ is
constant. \\
3. Equations (\ref{4.65}) are satisfied for 
$G_{ij,\infty}(\cdot)$ and $g_{\alpha \beta,\infty}(\cdot)$.
\end{proposition}
\begin{proof}
The proof is along the lines of the proof of Proposition \ref{4.39}.
\end{proof}

We now look at some special solutions of (\ref{4.65}).
Recall that $\rho \: : \: \pi_1(B, b) \rightarrow \SL(N)$ is
the holonomy representation.

\begin{proposition} \label{4.77}
Under the assumptions of Proposition \ref{4.76}, if $N = 0$ then
$\overline{g}_\infty(t) = t g_{Ein}$, where $g_{Ein}$ is an
Einstein metric on $M = B$ with Einstein constant $- \: \frac12$.
If $N = 1$ then
$\overline{g}_\infty(t)$ is locally an isometric product of
$\R$ or $S^1$ with $(B, t g_{Ein})$, where $g_{Ein}$ is an Einstein metric
on $B$ with Einstein constant $- \: \frac12$. For any $N$, if $\dim(B) = 1$
then with an appropriate choice of section $s$, we can locally write
$G_{ij}(b) \: = \: (e^{bX})_{ij}$ and $g_B \: = \: \frac{t}{2} \: \Tr(X^2) \:
db^2$, where $X$ is a real diagonal $(N \times N)$-matrix with vanishing 
trace.

If $\dim(B) = 2$ and $N = 2$ then $B$ has negative Euler
characteristic. Also, \\
1. $\overline{g}$ is a locally product metric and $B$ has 
sectional curvature $- \: \frac{1}{2t}$, or \\
2. $\rho$ fixes no point of the boundary of
$\SL(2, \R)/\SO(2) = H^2$ and with the right choice of orientation of
$\widetilde{B}$, the map $\widetilde{G} \: : \:
\widetilde{B} \rightarrow H^2$ is holomorphic.
\end{proposition}
\begin{proof}
The $N=0$ case is clear.
As $\det(G_{ij})$ is constant, if $N = 1$ then we are in a local
product situation.  For any $N$, if $\dim(B) = 1$ then
the map $b \rightarrow G_{ij}(b)$ describes a geodesic in
$\SL(N, \R)/\SO(N, \R)$, from which the proposition follows.
(See \cite[Example 4.27]{Lott (2007)}).

If $\dim(B) = 2$ and $N =2$ then we can consider $\widetilde{G}$ to be
a $\rho$-equivariant harmonic map $u \: : \: \widetilde{B} \rightarrow
H^2$. Choosing an orientation of $\widetilde{B}$,
we use a local complex coordinate $z$ on $\widetilde{B}$.
There is a solution to
the first equation in (\ref{4.65}) if and only if the representation 
$\rho \: : \: \pi_1(B) \rightarrow \SL(2, \R)$ is not conjugate
to a (nondiagonal) representation by upper triangular matrices
\cite{Jost-Zuo (1997),Labourie (1991)}.
If there is a solution to the first equation in (\ref{4.65}) then 
looking at the $dz^2$-component of the second equation in (\ref{4.65})
gives
\begin{equation} \label{4.77.5}
u_{z} 
\overline{u_{\overline{z}}} = 0.
\end{equation} 

We consider the subset of $\partial H^2$, the boundary at infinity
of $H^2$, which
is pointwise
fixed by $\Image(\rho)$. It is either all of $\partial H^2$, two points in 
$\partial H^2$, one point in $\partial H^2$ or the empty set.
If all of $\partial H^2$ 
is fixed by $\Image(\rho)$ then $\rho$ is the identity representation,
$u$ descends to a harmonic function on $B$ (which must be constant)
and $B$ has constant sectional curvature $- \: \frac{1}{2t}$.
If $\Image(\rho)$ fixes exactly two points of
$\partial H^2$ then $\rho$ is conjugate to a diagonal representation and 
$u$ maps to a nontrivial geodesic in $H^2$.
We can assume that $u$ is real-valued. Then equation (\ref{4.77.5})
implies that $u$ is constant, which is a contradiction.
As has been said, there is no solution to
the first equation in (\ref{4.65}) if $\Image(\rho)$ fixes a
single point of $\partial H^2$.
Finally, suppose that $\Image(\rho)$ fixes no point of
$\partial H^2$. Then $u$ is constant or $du$ has generic rank two. 
If $u$ is constant then $\overline{g}$ is a locally product metric.
Suppose that $u$ is nonconstant.
As $du$ has generic rank $2$,
equation (\ref{4.77.5}) implies
that $u$ is holomorphic or antiholomorphic. If $u$ is antiholomorphic then 
we change the orientation of $\widetilde{B}$ to make $u$ holomorphic. 
As $u$ is nonconstant,
Liouville's theorem implies that $B$ has negative Euler characteristic.
\end{proof}

\begin{remark} \label{4.78}
The solutions with $\dim(B) = 1$,
$G_{ij}(b) \: = \: (e^{bX})_{ij}$ and $g_B \: = \: \frac{t}{2} \: \Tr(X^2) \:
db^2$ are generalized Sol-solutions.
\end{remark}

\begin{remark} When $\dim(B) = 2$ and $N = 2$, the equations (\ref{4.65})
arose independently in the paper \cite{Song-Tian (2007)} on
K\"ahler-Ricci flow. In that paper, which is in the holomorphic setting,
the map $G$ arises as the classifying map for the torus bundle of an
elliptic fibration. The term 
$\frac14 \: G^{ij} \: G_{jk,\alpha} \: G^{kl} \: 
G_{li,\beta}$ of (\ref{4.65}) is called the Weil-Petersson term.
The second equation of (\ref{4.65}), in the K\"ahler case, is considered to be a
generalized K\"ahler-Einstein equation for the geometry of a collapsing
limit.
\end{remark}

\begin{remark} \label{4.79}
All of the results of this section extend
to the case when $B$ is an orbifold, $E$ is a flat orbifold 
${\mathcal G}$-bundle
over $B$, a manifold $M$ is the total space of an orbifold fiber bundle
$\pi \: : \: M \rightarrow B$ and ${\mathcal G}$ acts
locally freely on $M$ (via a map $E \times_B M \rightarrow M$)
with orbifold quotient $B$.
\end{remark}

\section{Equivalence classes of \'etale groupoids} \label{section5}

Let ${\frak G}$ be a complete effective path-connected 
Hausdorff \'etale groupoid that admits an invariant
Riemannian metric on the space of units $G^{(0)}$.
We assume that \\
1. ${\frak G}$ equals its closure $\overline{\frak G}$. \\
2. The local symmetry sheaf $\underline{\frak g}$
of ${\frak G}$ is a locally constant
sheaf of abelian
Lie algebras isomorphic to $\R^N$.

\begin{example} \label{5.1}
Let $M$ be the total space of a twisted abelian principal ${\mathcal G}$-bundle
as in Subsection \ref{subsection4.1}. We can take ${\frak G} \: = \:
E \times_B M$, where the flat bundle $E$ has the \'etale
topology, with ${\frak G}^{(0)} = M$. The local symmetry sheaf
comes from the flat vector bundle $\pi^* e$ on $M$. 

We can perform a similar construction in the setting
of Remark \ref{4.79}, where $M$ is a manifold and $B$ is an orbifold.
\end{example}

The results of Section \ref{section4}
extend to the setting of a Ricci flow on 
${\frak G}$, under the analogous curvature and diameter assumptions,
provided that ${\frak G}$ is locally free.
The reason is that the local structure of
such an \'etale groupoid is the same as the local structure
considered in Section \ref{section4} 
\cite[Corollary 3.2.2]{Haefliger (1985)}. 
We can then perform the integrals of
Section \ref{section4} over the orbit space of ${\frak G}$ and derive
the same consequences as in Section \ref{section4}.

It will be useful to determine the global structure
of such \'etale groupoids, at least in low dimensions.

\begin{proposition} \label{5.2}
Suppose that ${\frak G}$ is locally free.
Then the orbit space ${\mathcal O}$ is an orbifold. 
There is a flat (orbifold) $\R^N$-bundle $e$ on ${\mathcal O}$ associated
to ${\frak G}$.

If $\dim({\mathcal O}) = 1$ then ${\frak G}$ is classified by
the isomorphism class of $e$.

In general, if $e$ is trivial then
${\frak G}$ is equivalent to the groupoid
of a principal bundle over ${\mathcal O}$.
It is classified up to groupoid equivalence by the orbits of
$\GL(N, \R)$ on $\HH^2({\mathcal O}; \R^N$).
\end{proposition}
\begin{proof}
The proof is similar to the classification in 
\cite{Haefliger-Salem (1988)} of the transverse structure
of Riemannian foliations with low-codimension leaves. 
(As the paper \cite{Haefliger-Salem (1988)} considers 
Riemannian groupoids that may not equal their closure,
there is an additional step in \cite{Haefliger-Salem (1988)}
which consists of analyzing the restriction of
the groupoid to an orbit closure. Since we only deal
with \'etale groupoids that equal their closures,
we do not have to deal with this complication.)

Given $x \in {\frak G}^{(0)}$, let ${\mathcal O}_x$ be its orbit.
There is an invariant neighborhood of the orbit whose groupoid
structure is described by \cite[Corollary 3.2.2]{Haefliger (1985)}. 
In particular,
the point in the orbit space ${\mathcal O}$, corresponding to
${\mathcal O}_x$, has a neighborhood $U$ that is homeomorphic to 
$V/{\frak G}_x^x$, where $V$ is a representation space for
the isotropy group
${\frak G}_x^x$. This gives the orbifold structure on the orbit space.

The classification of such \'etale groupoids comes from the
bundle theory developed in \cite[Section 2.3]{Haefliger (1985)},
which we now follow. For notation,
if $G$ is a topological group then let $G_\delta$ denote $G$ with the
discrete topology.
Suppose first that the isotropy groups ${\frak G}_x^x$ are trivial,
so the orbifold ${\mathcal O}$ is a
manifold. Let $U \subset {\mathcal O}$ be a neighborhood of 
${\mathcal O}_x$ as above. Let $\pi \: : \: {\frak G}^{(0)} \rightarrow
{\mathcal O}$ be the quotient map. 
By \cite[Corollary 3.2.2]{Haefliger (1985)},
the restriction of ${\frak G}$ to $\pi^{-1}(U)$ is equivalent to
the cross-product groupoid $(\R^N \times U) \rtimes \R^N_\delta$,
where $\R^N_\delta$ acts
on $\R^N$ by translation and acts trivially on $U$.
This gives the local structure of ${\frak G}$. It remains to determine
the possible ways to glue these local structures together.

To follow the notation of \cite[Section 2.1]{Haefliger (1985)},
put $\Gamma = \R^N_\delta \subset \Diff(\R^N)_\delta$. The normalizer $N^\Gamma$ of $\Gamma$ in 
$\Diff(\R^N)$ is $\R^N \widetilde{\times} \GL(N, \R)$ and the centralizer
is $C^\Gamma = \R^N$. We give $N^\Gamma$ the topology
$\R^N \widetilde{\times} \GL(N, \R)_\delta$. 

Following the discussion in \cite[Section 2.1]{Haefliger (1985)},
suppose that $U \subset {\mathcal O}$ is an open set.
Consider the cross-product groupoid $(\R^N \times U) \rtimes \R^N_\delta$.
Let
${\mathcal E}(U)$ be the self-equivalences of 
$(\R^N \times U) \rtimes \R^N_\delta$ that
project onto the identity of $U$. This forms a sheaf
$\underline{\mathcal E}$ on ${\mathcal O}$. 
We can cover ${\mathcal O}$ by open sets $U$ such that
$\pi^{-1}(U)$ is equivalent to $(\R^N \times U) \rtimes \R^N_\delta$.
It follows that 
the \'etale groupoids in question are classified by the set
$\HH^1({\mathcal O}; \underline{\mathcal E})$ 
\cite[Proposition 2.3.2]{Haefliger (1985)}.

To compute $\HH^1({\mathcal O}; \underline{\mathcal E})$,
let $\underline{\R^N}$
be the sheaf on ${\mathcal O}$ for which $\underline{\R^N}(U)$ consists
of smooth maps $U \rightarrow \R^N$, let
$\R^N_\delta$ (also) denote the constant sheaf on
${\mathcal O}$ with stalk $\R^N_\delta$ and
let
$\GL(N, \R)_\delta$ (also) denote the constant sheaf on
${\mathcal O}$ with stalk $\GL(N, \R)_\delta$.
As in \cite[(2.4.2)]{Haefliger (1985)} there is a short exact sequence
of sheaves
\begin{equation} \label{5.3}
0 \longrightarrow \underline{\R^N}/\R^N_\delta \longrightarrow
\underline{\mathcal E} \longrightarrow
\GL(N, \R)_\delta \longrightarrow 0.
\end{equation} 

From \cite[Th\'eor\`eme 1.2]{Frenkel (1957)}, this short exact sequence
of sheaves
gives rise to an exact sequence of pointed sets
\begin{equation} \label{5.4}
\ldots \longrightarrow
\HH^0({\mathcal O}; \GL(N, \R)_\delta) \longrightarrow
\HH^1({\mathcal O}; \underline{\R^N}/\R^N_\delta) \longrightarrow
\HH^1({\mathcal O}; \underline{\mathcal E}) \longrightarrow
\HH^1({\mathcal O}; \GL(N, \R)_\delta).
\end{equation}
The set $\HH^1({\mathcal O}; \GL(N, \R)_\delta)$ is the same as
the set of homomorphisms $\pi_1({\mathcal O}) \rightarrow
\GL(N, \R)$ modulo conjugation by elements of $\GL(N, \R)$ or,
equivalently, the set of equivalence classes of flat $\R^N$-vector
bundles on ${\mathcal O}$. 
The image of the classifying element of ${\frak G}$,
under the map $\HH^1({\mathcal O}; \underline{\mathcal E}) \longrightarrow
\HH^1({\mathcal O}; \GL(N, \R)_\delta)$, classifies the flat $\R^N$-vector
bundle $e$ mentioned in Proposition \ref{5.2}. 
More explicitly, the transition functions of
$e$ come from the image under 
$\underline{\mathcal E} \longrightarrow
\GL(N, \R)_\delta$ of the transition functions of ${\frak G}$.

The short exact sequence
\begin{equation} \label{5.5}
0 \longrightarrow 
\R^N_\delta \longrightarrow
\underline{\R^N} \longrightarrow
\underline{\R^N}/\R^N_\delta \rightarrow 0
\end{equation}
of sheaves of abelian groups gives a long exact sequence
\begin{equation} \label{5.6}
\ldots \longrightarrow
\HH^1({\mathcal O}; \underline{\R^N}) \longrightarrow
\HH^1({\mathcal O}; \underline{\R^N}/\R^N_\delta) \longrightarrow
\HH^2({\mathcal O}; \R^N_\delta) \longrightarrow
\HH^2({\mathcal O}; \underline{\R^N}) \longrightarrow \ldots
\end{equation}
of abelian groups.
As $\underline{\R^N}$ is a fine sheaf, it follows from (\ref{5.6}) that
$\HH^1({\mathcal O}; \underline{\R^N}/\R^N_\delta) \cong
\HH^2({\mathcal O}; \R^N_\delta) \: = \: 
\HH^2({\mathcal O}; \R^N)$.

As $\HH^0$ consists of global sections, 
(\ref{5.4}) gives an exact sequence of pointed sets
\begin{equation} \label{5.7}
\GL(N, \R) \longrightarrow
\HH^2({\mathcal O}; \R^N) \longrightarrow
\HH^1({\mathcal O}; \underline{\mathcal E}) \longrightarrow
\HH^1({\mathcal O}; \GL(N, \R)_\delta).
\end{equation}

If $\dim({\mathcal O}) = 1$ then from (\ref{5.7}), the map
$\HH^1({\mathcal O}; \underline{\mathcal E}) \longrightarrow
\HH^1({\mathcal O}; \GL(N, \R)_\delta)$ is injective.
Thus ${\frak G}$ is determined up to groupoid equivalence
by the isomorphism class of the flat vector bundle $e$. 

If ${\mathcal O}$ has arbitrary dimension, suppose that
$e$ is trivial.
Consider the preimage under
$\HH^1({\mathcal O}; \underline{\mathcal E}) \longrightarrow
\HH^1({\mathcal O}; \GL(N, \R)_\delta)$ of the element
in $\HH^1({\mathcal O}; \GL(N, \R)_\delta)$ corresponding to the
identity representation. By (\ref{5.7}), this preimage can be
identified with the orbit space for the action of
$\GL(N, \R)$ on $\HH^2({\mathcal O}; \R^N)$. Any such
orbit contains an element of 
$\Image( \HH^2({\mathcal O}; \Z^N) \rightarrow \HH^2({\mathcal O}; \R^N))$,
which implies that ${\frak G}$ is equivalent to the \'etale groupoid
arising from some principal $T^N$-bundle on ${\mathcal O}$.

The preceding considerations extend to the case when the
(finite) isotropy groups 
${\frak G}^x_x$ are not all trivial.
In that case, ${\mathcal O}$ is an orbifold and the argument
extends to the orbifold setting.  
For example, $\HH^*({\mathcal O}; \R^N)$ has to
be interpreted as an orbifold cohomology group.
\end{proof}

\begin{remark} \label{5.8}
If one starts with an (untwisted) principal
${\mathcal G}$-bundle, with ${\mathcal G}$ abelian,
then the triviality of the corresponding
\'etale groupoid is determined by whether or not
$\{\int_B F^i\}_{i=1}^N$ vanishes in $\HH^2(B; \R^N)$.

Suppose that the \'etale groupoid is nontrivial
and $\{\overline{g}_j\}_{j=1}^\infty$ is a sequence of invariant metrics 
on the
principal ${\mathcal G}$-bundle, so that there is a limiting invariant
metric $\overline{g}_\infty$. 
It is possible that
the curvatures $\{F^i\}_{i=1}^N$
approach zero in norm as $j \rightarrow \infty$.  If this is the case 
then $\overline{g}_\infty$ will live on a distinct \'etale groupoid,
as its curvature $\{F^i\}_{i=1}^N$ vanishes.  This phenomenon
occurs in the rescaled Ricci flow on the unit circle bundle of a
surface of constant negative curvature.  

On the other hand, if
we start with a trivial \'etale groupoid and 
$\{\overline{g}_j\}_{j=1}^\infty$ is a noncollapsing
sequence of invariant metrics on the
principal ${\mathcal G}$-bundle then any limiting invariant metric
$\overline{g}_\infty$ will necessarily 
be on the same \'etale groupoid.

The relevance of Proposition \ref{5.2} is that for
\'etale groupoids which satisfy its hypotheses, we can discuss convergence
of Ricci flow solutions on such \'etale groupoids
in terms of convergence of invariant
Ricci flow solutions on twisted principal bundles.
\end{remark}

\begin{example} \label{5.9}
Suppose that $M$ is the total space of a principal $S^1$-bundle over 
a compact oriented surface $B$.
Given a subgroup $\Z_k \subset S^1$, let $M/\Z_k$ be the quotient
space. It is also the total space of a principal $S^1$-bundle over $B$.

The (discrete) $S^1$-action on a principal $S^1$-bundle gives 
an \'etale groupoid. The map $M \rightarrow M/\Z_k$ gives an equivalence of
\'etale groupoids, in the sense of
\cite[Chapter III.${\mathcal G}$.2.4]{Bridson-Haefliger}.
However, the Euler class of the circle bundle
$M/\Z_k \rightarrow B$ is $k$ times that of the circle bundle
$M \rightarrow B$. This shows that the 
Euler class of the circle bundle is not an invariant of 
the groupoid equivalence class.
Instead, all that is relevant is
whether or not the rational Euler class vanishes.
\end{example}

\begin{example} \label{5.10}
If $\dim({\mathcal O}) = 1$ then
any homomorphism $\alpha \: : \: \pi_1({\mathcal O}) \rightarrow
\GL(N, \R)$ gives rise to an \'etale groupoid with unit space
${\frak G}^{(0)} = \R^N \times_\alpha \widetilde{\mathcal O}$.

If ${\mathcal O}$ is a closed orientable $2$-dimensional orbifold
then
$\HH^2({\mathcal O}; \R^N) \cong \R^N$ and 
the action of
$\GL(N, \R)$ on $\HH^2({\mathcal O}; \R^N)$ has
two orbits, namely the zero element and the nonzero elements.
Thus if $e$ is trivial then there are two equivalence classes
of such groupoids with orbit space ${\mathcal O}$, 
one corresponding to a vanishing
``Euler class'' and one corresponding to a nonvanishing
``Euler class''.
\end{example}

Suppose that $M$ is the total space of a twisted principal $\R^N$-bundle.
Let $\overline{g}$ be an invariant metric on $M$. We recall that 
there are two distinct connections in this situation, the flat connection on
the twisting bundle ${\mathcal E}$ and the connection $A$ on the twisted
principal bundle.  We will use the following lemma later.

\begin{lemma} \label{5.11}
Let
$\pi \: : \: M \rightarrow B$ be a twisted principal $\R^N$-bundle.
Given $G_{ij}$ and $g_{\alpha \beta}$, let 
$A_1$ and $A_2$ be two flat connections on $M$. Let
$\overline{g}_1$ and $\overline{g}_2$ be the corresponding invariant
metrics on $M$. Then their underlying Riemannian groupoids are equivalent.
\end{lemma}
\begin{proof}
Let $\{U_i\}$ be a covering of $B$ by open contractible sets.
Let ${\mathcal U} = \{ \pi^{-1}(U_i) \}$ be the corresponding
covering of $M$
and let ${\frak G}_{\mathcal U}$ be the
localization of ${\frak G}$ \cite[Section 5.2]{Lott (2007)}.
In our case, elements of ${\frak G}_{\mathcal U}$ are quadruples
$(i,p_i,p_j,j)$ with $p_i \in \pi^{-1}(U_i)$,
$p_j \in \pi^{-1}(U_j)$ and $\pi(p_i) = \pi(p_j)$.
The multiplication is $(i,p_i,p_j,j) \cdot (j,p_j,p_k,k) \: = \:
(i,p_i,p_k,k)$.
The units ${\frak G}_{\mathcal U}^{(0)}$ are quadruples $(i,p_i,p_i, i)$
and the source and range maps are
$s(i,p_i,p_j,j) = (j,p_j,p_j,j)$ and
$r(i,p_i,p_j,j) = (i,p_i,p_i,i)$.
Let $s^1_i \: : \: U_i \rightarrow \pi^{-1}(U_i)$ be a section
for which $(s^1_i)^* A_1 = 0$. Similarly, let
$s^2_i \: : \: U_i \rightarrow \pi^{-1}(U_i)$ be a section
for which $(s^2_i)^* A_2 = 0$. Define a map 
$F \: : \: {\frak G}_{\mathcal U} \rightarrow {\frak G}_{\mathcal U}$ by 
$F(i,p_i,p_j,j) = (i, p_i + s^2_i(u_i) - s^1_i(u_i), 
p_j + s^2_j(u_j) - s^1_j(u_j), j)$, where 
$u_i = \pi(p_i)$, $u_j = \pi(p_j)$ and we write the action of
$\R^N$ additively. Then $F$ is a groupoid isomorphism. 
On the space of units, 
$F(i,p_i,p_i,i) = (i, p_i + s^2_i(u_i) - s^1_i(u_i), 
p_i + s^2_i(u_i) - s^1_i(u_i), i)$ and so $F$ sends the section
$s^1_i$ to $s^2_i$. It follows that $F$ is an isomorphism of
Riemannian groupoids.
\end{proof}

\section{Convergence arguments and universal covers} \label{section6}

In this section we prove Theorem \ref{1.2}. In Subsection \ref{subsection6.1} 
we
prove convergence to a
locally homogeneous Ricci flow on an \'etale groupoid. In
Subsection \ref{subsection6.2} we promote this to convergence on the
universal cover of $M$. 

\subsection{Convergence arguments} \label{subsection6.1}

In this subsection we show that under the hypotheses of Theorem \ref{1.2},
there is a rescaling limit which is a locally homogeneous expanding
soliton solution on an \'etale groupoid.  
To do this, if $\overline{\mathcal O}$ is the
closure of the orbit of 
$g(\cdot)$ under the action of the parabolic rescaling semigroup $\R^{\ge 1}$
then we define a stratification
of $\overline{\mathcal O}$ in terms of the number of local symmetries.
We let $k_0$ be the maximal number of local symmetries that can occur
in a rescaling limit of $g(\cdot)$. This corresponds to a maximally
collapsed limit. The first step is to show that
$k_0$ determines the Thurston type of $M$, and that there is a sequence
of rescalings of $g(\cdot)$ which approaches the corresponding
locally homogeneous expanding soliton.

In order to show that any rescaling limit
$\overline{g}(\cdot)$ is a locally homogeneous
expanding soliton (except possibly in the $\widetilde{\SL_2(\R)}$ case),
we use further arguments.  We show that any rescaling limit has 
$k_0$ local symmetries.
We then use a compactness argument, along with the local stability
of the space of expanders, to show that $\overline{g}(\cdot)$ is a
locally homogeneous expanding soliton.

Let $g(\cdot)$ be a Ricci flow solution on a connected closed $3$-manifold $M$, 
defined for $t \in (1, \infty)$, with
$\sup_{t \in (1, \infty)} t \parallel \Riem(g(t)) \parallel_\infty \: \le \: 
K \: < \: \infty$ and $\sup_{t \in (1, \infty)} t^{- \: \frac12} \:
\diam(g(t)) \: \le \: 
D \: < \: \infty$. From Proposition \ref{3.5}, $M$ has a single geometric
piece.
Given $s \in [1, \infty)$, put $g_s(t) \: = \: \frac{1}{s} \: g(st)$.
Then for all $s$, we have
$\sup_{t \in (1, \infty)} t \parallel \Riem(g_s(t)) \parallel_\infty \: 
\le \: K$ and $\sup_{t \in (1, \infty)} t^{- \: \frac12} \:
\diam(g_s(t)) \: \le \: D$. By Proposition \ref{3.2}, the family of Ricci flow
solutions $\{g_s(\cdot)\}_{s \in [1, \infty)}$ is sequentially 
precompact among
Ricci flow solutions on \'etale groupoids.

Let $\overline{\mathcal O}$ be the sequential closure of the forward orbit 
$\{g_s(\cdot)\}_{s \in [1, \infty)}$. Let
$\overline{\mathcal O}_{(k)}$ be the elements of $\overline{\mathcal O}$
with a $k$-dimensional local symmetry sheaf $\underline{\frak g}$.

\begin{lemma} \label{6.1}
If $\widehat{g}(\cdot) \in \overline{\mathcal O}$ then
the underlying \'etale groupoid of $\widehat{g}(\cdot)$ is locally free.
\end{lemma}
\begin{proof}
If $\widehat{g}(\cdot) \in \overline{\mathcal O}_{(0)}$ then there is
nothing to show.  If $\widehat{g}(\cdot) \in \overline{\mathcal O}_{(1)}$
then the lemma follows from the fact there is no point
$x \in {\frak G}^{(0)}$ where the local Killing vector fields vanish
simultaneously.

Suppose that $\widehat{g}(\cdot) \in \overline{\mathcal O}_{(2)}$.
Write $\widehat{g}(\cdot) = \lim_{i \rightarrow \infty} 
(M, g_{s_j^\prime}(\cdot))$
for some sequence $\{s_j^\prime\}_{j=1}^\infty$ tending to infinity.
By \cite{Cheeger-Fukaya-Gromov (1992)},
for any $\epsilon > 0$, there is an integer $J_\epsilon < \infty$ so that
if $j \ge J_\epsilon$ then there is a locally
$T^2$-invariant Riemannian metric
$g^\prime_j$ on $M$ which is $\epsilon$-close in the $C^1$-topology to
$\frac{1}{s_j^\prime} g(s_j^\prime)$. Furthermore, one can take the
sectional curvature of $g^\prime_j$ to be uniformly bounded in $\epsilon$
\cite[Theorem 2.1]{Rong (1996)}. The collapsing is along
the $T^2$ fibers. Taking a sequence of values of $\epsilon$ going to zero
and choosing $j \ge J_\epsilon$,
after passing to a subsequence we can say that
$\widehat{g}(1) = \lim_{j \rightarrow \infty} (M, g^\prime_j)$.

Let $S$ be the orbit space of the \'etale groupoid.  It is a circle or
an interval. If $S$ is a circle then $M$ is the total space of a
$T^2$-bundle over $S^1$. (The fibers cannot be Klein bottles since $M$ is
orientable.) Hence the 
local $T^2$-action on $(M, g^\prime_j)$ is free. 
Let $H \in \SL(2, \Z)$ be the holonomy of the $T^2$-bundle, defined up
to conjugation in $\SL(2, \Z)$. Given $M$, there is a finite number
of possibilities for $H$, as follows from 
\cite[pp. 439,469-470,481-482]{Scott (1983)}. 
After passing to a subsequence, we can assume that
there is a single such $H$.
For each $j$, the $T^2$-bundle
with invariant metric $g^\prime_j$ is the total space of a
twisted principal $T^2$-bundle over $S^1$, where the
twisting bundle $E$ is a flat $T^2$-bundle on $S^1$ with holonomy $H$.
From Proposition \ref{5.2}, for all $j$ these give rise to
equivalent \'etale groupoids. Looking at how one constructs the
limiting Riemannian groupoid as $j \rightarrow \infty$ 
\cite[Proposition 5.9]{Lott (2007)}, it follows that 
$\widehat{g}(\cdot)$ is defined on this same \'etale groupoid.
In particular, it is locally free.

If $S$ is an interval
then as in the proof of Proposition \ref{3.5}, the asphericity of 
$M$ implies that the local $T^2$-action on $M$ is locally free.
Then $M$ is the total space of an orbifold $T^2$-bundle over
the orbifold $S$. As $S$ is double covered by a circle, we can
take a double cover $\widehat{M}$ of $M$ which is the total space
of a $T^2$-bundle over $S^1$. Applying the preceding argument
$\Z_2$-equivariantly to $\widehat{M}$, we conclude that 
the underlying \'etale groupoid of $\widehat{g}(\cdot)$ is 
again locally free.

Finally, suppose that $\widehat{g}(\cdot) \in \overline{\mathcal O}_{(3)}$.
Write $\widehat{g}(\cdot) = \lim_{i \rightarrow \infty} 
(M, g_{s_j^\prime}(\cdot))$
for some sequence $\{s_j^\prime\}_{j=1}^\infty$ tending to infinity.
Then the orbit space $S$ of the \'etale groupoid is a point and
$\left\{ \left( M, \frac{1}{s_j^\prime} g(s_j^\prime) \right) 
\right\}_{j=1}^\infty$ Gromov-Hausdorff 
converges, with bounded sectional curvature, to a point. 
That is, $M$ is almost flat
and so is an infranilmanifold \cite{Gromov (1978)}. 
There is a finite normal cover ${M}_0$ of $M$
which is diffeomorphic
to a flat manifold or a nilmanifold. Let 
$g_0(\cdot)$ be the lift of $g(\cdot)$ to $M_0$ and let
$\widehat{g}_0(\cdot)$ be the corresponding limiting Ricci flow on an
\'etale groupoid, with $\widehat{g}(\cdot)$ as a finite quotient.
By \cite{Cheeger-Fukaya-Gromov (1992)},
for any $\epsilon > 0$, there is an integer $J_\epsilon < \infty$ so that
if $j \ge J_\epsilon$ then there is a left-invariant
Riemannian metric
$g^\prime_j$ on ${M}_0$, of $\R^3$ or $\Nil$-type, 
which is $\epsilon$-close in the $C^1$-topology to
$\frac{1}{s_j^\prime} g_0(s_j^\prime)$. Furthermore, one can take the
sectional curvature of $g^\prime_j$ to be uniformly bounded in $\epsilon$
\cite[Theorem 2.1]{Rong (1996)}. The collapsing is along
all of $M_0$. Taking a sequence of values of $\epsilon$ going to zero
and choosing $j \ge J_\epsilon$,
after passing to a subsequence we can say that
$\widehat{g}_0(1) = \lim_{j \rightarrow \infty} (M_0, g^\prime_j)$.
Looking at how one constructs the
limiting Riemannian groupoid as $j \rightarrow \infty$ 
\cite[Proposition 5.9]{Lott (2007)}, it follows that 
the underlying \'etale groupoid of $\widehat{g}_0(1)$ is a
cross-product groupoid
$\R^3 \rtimes \R^3_\delta$ or $\Nil \rtimes \Nil_\delta$, where 
$\delta$ denotes the discrete topology. Hence the
underlying \'etale groupoid of $\widehat{g}(\cdot)$ is locally free.
\end{proof}

The relevance of Proposition \ref{6.1} is that it allows us
to use Proposition \ref{4.76} to analyze
blowdown limits of $\widehat{g}(\cdot)$.

Let $k_0$ be the largest $k$ so that 
$\overline{\mathcal O}_{(k)}$ is nonempty. For simplicity of terminology,
we will say that a Ricci flow on an \'etale groupoid is a locally
homogeneous expanding soliton if there is some homogeneous expanding soliton
to which the Ricci flow on the unit space of the \'etale groupoid is
locally isometric.

\begin{proposition} \label{6.2}
If $k_0 \: = \: 0$ then $M$ admits an $H^3$-structure.
If $k_0 \: = \: 1$
then $M$ admits an $H^2 \times \R$ or $\widetilde{\SL_2(\R)}$-structure. If
$k_0 \: = \: 2$ then $M$ admits a $\Sol$-structure.
If $k_0 \: = \: 3$ then $M$ admits an $\R^3$ or $\Nil$-structure.

In any case, there is a sequence $\{s_j\}_{j=1}^\infty$ tending to
infinity so that $\lim_{j \rightarrow \infty} (M, g_{s_j}(\cdot))$ 
exists as a Ricci flow solution on an \'etale groupoid, and is a 
locally homogeneous expanding soliton of type
\begin{itemize}
\item $H^3$ if $k_0 = 0$,
\item $H^2 \times \R$ if $k_0 = 1$,
\item $\Sol$ if $k_0 = 2$,
\item $\R^3$ or $\Nil$ if $k_0 = 3$.
\end{itemize}
\end{proposition}
\begin{proof}
Given $\widehat{g}(\cdot) \in \overline{\mathcal O}_{(k_0)}$, put
$\widehat{g}_s(t) \: = \: \frac{1}{s} \: \widehat{g}(st)$.
We claim that
the forward orbit $\{ \widehat{g}_s(\cdot)\}_{s \in [1,\infty)}$ is relatively
sequentially compact in $\overline{\mathcal O}_{(k_0)}$. To see this,
suppose that
there is a sequence $\{ \widehat{g}_{s_i}(\cdot)\}_{i=1}^\infty$
having a limit $\widehat{g}^\prime(\cdot)$. We can find a subsequence of
$\{(M, g_s(\cdot))\}_{s \in [1, \infty)}$ that converges to
$\widehat{g}^\prime(\cdot)$. Thus $\widehat{g}^\prime(\cdot) \in
\overline{\mathcal O}$. However, the number of local symmetries cannot
decrease in the limit. Hence $\widehat{g}^\prime(\cdot) \in
\overline{\mathcal O}_{(k)}$ for some $k \ge k_0$. We must have
$k = k_0$, by the definition of $k_0$, which proves the claim.

Let $\{s_i\}_{i=1}^\infty$
be a sequence tending to infinity such that $\lim_{i \rightarrow \infty} 
\widehat{g}_{s_i}(\cdot) \: = \: \widehat{g}_\infty(\cdot)$ for some
$\widehat{g}_\infty(\cdot) \in \overline{\mathcal O}_{(k_0)}$.
Let $S$ denote the underlying orbit space of
$\widehat{g}_\infty(1)$.
There is a sequence $\{s_j^\prime\}_{j=1}^\infty$ tending to infinity
so that $\lim_{j \rightarrow \infty} 
\left( M, g_{s_j^\prime}(\cdot) \right) = \widehat{g}_\infty(\cdot)$.
In particular, 
$\lim_{j \rightarrow \infty} 
\left( M, \frac{1}{s_j^\prime} g(s_j^\prime) \right) \stackrel{GH}{=} S$.

If $k_0 = 0$ then by Proposition \ref{4.77}, $(M, \widehat{g}_\infty(\cdot))$
is the Ricci flow on a manifold of constant negative sectional
curvature.

If $k_0 = 1$ then $S$ is a closed
two-dimensional orbifold. Taking a double cover if
necessary, we can assume that $S$ is orientable.
From Proposition \ref{5.2},
we can assume that the underlying \'etale groupoid comes from an
orbifold principal $S^1$-bundle on $S$.
(The triviality of $e$ 
comes from its identification with $\HH^1$ of the circle
fiber of the orbifold bundle $M \rightarrow S$.)
By Proposition \ref{4.77}, $\widehat{g}_\infty(\cdot)$ has
$(H^2 \times \R)$-type and $S$ has a metric of constant 
curvature $- \: \frac{1}{2t}$.
As $M$ is the total space
of an orbifold circle bundle over $S$, it follows that $M$ admits an
$H^2 \times \R$ or $\widetilde{\SL_2(\R)}$-structure
(using \cite{Meeks-Scott (1986)} if we took a double cover).

If $k_0 = 2$ then $S$ is $S^1$ or an interval $[0,L]$. 
Suppose first that $S = S^1$. Then
$M$ is the total space of a $T^2$-fiber bundle over $S$.
Let $H \in \SL(2, \Z)$ be the holonomy of the fiber bundle,
defined up to conjugacy. 
As in the proof of Lemma \ref{6.1}, 
the \'etale groupoid of $\widehat{g}_\infty(\cdot)$ arises from a
(twisted) principal $T^2$-bundle
on $S^1$. The flat bundle $e$ over $S^1$ has holonomy $H \in \SL(2,\Z)$.
By Proposition \ref{4.77}, $\widehat{g}_\infty(\cdot)$ has $\Sol$-type and
$H$ is a hyperbolic element of
$\SL(2, \Z)$. Thus $M$ admits a $\Sol$-structure.

Suppose now that $S = [0,L]$.
As in the proof of Lemma \ref{6.1}, $M$ is the total
space of an orbifold $T^2$-bundle over the orbifold $[0,L]$.
A double cover $\widehat{M}$ of $M$ fibers over $S^1$. 
Running the previous argument on $\widehat{M}$ with the
pullback metric, we conclude that
$\widehat{M}$ admits a $\Sol$-structure.
Hence $M$ admits a $\Sol$-structure \cite{Meeks-Scott (1986)}.

If $k_0 = 3$ then $S$ is a point. Hence 
$\left\{ \left( M, \frac{1}{s_j^\prime} g(s_j^\prime) \right) 
\right\}_{j=1}^\infty$ Gromov-Hausdorff 
converges, with bounded sectional curvature, to a point. 
As in the proof of Lemma \ref{6.1},
$\widehat{g}_\infty(\cdot)$ is locally homogeneous and has
$\R^3$ or $\Nil$ as its local symmetry group. 
Such a Ricci flow solution is automatically
a locally homogeneous expanding soliton.
\end{proof}

We have shown that there is some sequence $\{s_j\}_{j=1}^\infty$
tending to infinity so that $\lim_{j \rightarrow \infty}
(M, g_{s_j}(\cdot))$ exists and is a locally homogeneous expanding
soliton. We now
wish to show that this is true for any sequence $\{s_j\}_{j=1}^\infty$
tending to infinity, at least if the Thurston type of $M$ is not
$\widetilde{\SL_2(\R)}$. The first step is to show that under a compactness
assumption,
there is a parameter $T$ so that if we take any rescaling limit
$\overline{g}(\cdot)$ then
upon further rescaling of
$\overline{g}(\cdot)$, the result 
is near a locally homogeneous expanding soliton
for some rescaling parameter $s \in [1,T]$.

\begin{proposition} \label{6.3}
Given $k$, let $C$ be a sequentially compact subset of
$\overline{\mathcal O}_{(k)}$. Let $U$ be a neighborhood of
\begin{itemize}
\item The $H^3$-type locally homogeneous expanding solitons in $\overline{\mathcal O}_{(0)}$ 
if $k = 0$,
\item The $(H^2 \times \R)$-type locally homogeneous 
expanding solitons in $\overline{\mathcal O}_{(1)}$ 
if $k = 1$,
\item
The $\Sol$-type locally homogeneous expanding solitons in $\overline{\mathcal O}_{(2)}$ 
if $k = 2$,
\item The $\R^3$-type and $\Nil$-type locally homogeneous expanding 
solitons in
$\overline{\mathcal O}_{(3)}$ if $k = 3$.
\end{itemize}

Then there is a $T = T(k,C,U) \in [1, \infty)$ so that for any
$\overline{g}(\cdot) \in C$, if 
$\overline{g}_s(\cdot) \in C$ for all $s \in [1, T]$ then there is some
$s \in [1, T]$ such that $\overline{g}_s(\cdot) \in U$.
\end{proposition}
\begin{proof}
Given $k$, $C$ and $U$, 
suppose that the proposition is not true.  
Then for each $j \in \Z^+$, there is
some 
$\overline{g}^{(j)}(\cdot) \in C$ 
so that for each $s \in [1, j]$, $\overline{g}^{(j)}_s(\cdot) \in C$ and
$\overline{g}^{(j)}_s(\cdot) \notin U$.
Take
a convergent subsequence of the $\{\overline{g}^{(j)}(\cdot)\}_{j=1}^\infty$ 
with limit 
$\overline{g}^{(\infty)}(\cdot)$. Then for all $s \in [1, \infty)$, we have
$\overline{g}_s^{(\infty)}(\cdot) \in C$ and
$\overline{g}_s^{(\infty)}(\cdot) \notin U$.
By sequential compactness, 
there is a sequence $\{t_k\}_{k=1}^\infty$ in $\Z^+$ tending
to infinity so that $\lim_{k \rightarrow \infty} 
\overline{g}^{(\infty)}_{t_k}(\cdot)$
exists and equals some 
$\overline{g}^{(\infty)}_\infty(\cdot) \in C$.
By Proposition \ref{4.77}, $\overline{g}^{(\infty)}_\infty(\cdot)$ is a
locally homogeneous expanding
soliton as in the statement of the present proposition. Then for large
$k$, we have $\overline{g}^{(\infty)}_{t_k}(\cdot) \in U$,
which is a contradiction.
\end{proof}

The next step is to use local stability to say that after
rescaling $\overline{g}(\cdot)$
by the parameter $T$, the result is definitely near a
locally homogeneous expanding soliton solution.

\begin{proposition} \label{6.4}
Suppose that $M$ does not have
Thurston type $\widetilde{\SL_2(\R)}$. Then there are decreasing open sets
$\{U_l\}_{l=1}^\infty$ 
of the type described in Proposition \ref{6.3}, 
whose intersection is the corresponding
set of locally homogeneous
expanding soliton solutions, so that under the hypotheses of
Proposition \ref{6.3}, if $T_l = T(k,C,U_l)$ then
we are ensured that $\overline{g}_{T_l}(\cdot) \in U_l$.
(In the case $k = 1$ we restrict to 
Ricci flow solutions on an \'etale groupoid with vanishing Euler class,
so $U_l$ is a neighborhood in the relative topology.) 
\end{proposition}
\begin{proof}
This follows from the local stability of the expanding solitons in
$\overline{\mathcal O}_{(k)}$. That is, there is a sequence
$\{U_l\}_{l=1}^\infty$ of such neighborhoods so that 
$\overline{g}_{s}(\cdot) \in U_l$ implies that 
$\overline{g}_{s^\prime}(\cdot) \in U_l$ whenever $s^\prime \ge s$.
(In fact, one has exponential convergence to the set of expanding
solitons.)
The case $k=0$ appears in \cite{Ye (1993)}. The case $k=2$ appears in 
\cite{Knopf}. In the case $k = 1$, recall from Example \ref{5.10}
that there are two relevant types of \'etale groupoids,
one with vanishing Euler class and one with nonvanishing Euler class.
The locally homogeneous expanding solitons live on \'etale groupoids
with vanishing Euler class.
Their local stability (modulo the center manifold), among Ricci flows
on \'etale groupoids with vanishing Euler class, is shown in
\cite{Knopf}. We remark that if $M$ has Thurston type $H^2 \times \R$ then
a limit $\lim_{j \rightarrow \infty} (M, g_{s_j}(\cdot))$ can only
be a Ricci flow on an \'etale groupoid with vanishing Euler class.

Note that if $k = 1$ then there may be a moduli space of locally
homogeneous
expanding
solitons of type $H^2 \times \R$ in $\overline{\mathcal O}_{(1)}$, corresponding to
various
metrics of constant curvature $- \: \frac12$ on the orbit space.  However, because
of our diameter bound, the moduli space is compact. 
Comparing with \cite{Knopf}, it may appear that there
is also a factor in the moduli space consisting of harmonic $1$-forms
on the orbit space. However, by Lemma \ref{5.11}, the various harmonic
$1$-forms all give equivalent geometries.
\end{proof}

\begin{remark} \label{6.5}
There is no locally homogeneous expanding soliton solution on
a three-dimensional \'etale groupoid of the type considered in 
Proposition \ref{6.2}
if it has an orbifold surface base with negative Euler characteristic, 
and a nonvanishing
Euler class.
A Ricci flow on such an \'etale groupoid will have a rescaling sequence
that converges to an $(H^2 \times \R)$-type expander on an \'etale
groupoid with vanishing Euler class.

In order to show convergence of the Ricci flow on a $3$-manifold with
Thurston type $\widetilde{\SL_2(\R)}$, at least by our methods, one
would have to show that the expanding solitons of type $H^2 \times \R$ are
also locally stable if one considers neighborhoods that include
\'etale groupoids with nonvanishing Euler class. The difficulty is that
the nearby Ricci flows live on an inequivalent groupoid and so one cannot
just linearize around the $(H^2 \times \R)$-type expanding solitons.
One approach would be to instead consider Ricci flows with 
$\widetilde{\SL_2(\R)}$-symmetry on \'etale groupoids with 
nonzero Euler class and show that 
this finite-dimensional family is an attractor.
\end{remark}

We now show if $k_0 < 3$ and $M$ does not have Thurston type
$\widetilde{\SL_2(\R)}$ then any rescaling limit $\overline{g}(\cdot)$
is a locally
homogenous expanding soliton. The method of proof is to
show that we can rescale $\overline{g}(\cdot)$ backward by a factor $T$, and
then apply the previous proposition.

\begin{proposition} \label{6.6}
If $k_0 < 3$ and $M$ does not have Thurston type
$\widetilde{\SL_2(\R)}$
then for any sequence $\{s_j\}_{j=1}^\infty$ tending to infinity,
as $j \rightarrow \infty$, 
$(M, g_{s_j}(\cdot))$ approaches the set of locally homogeneous expanding
solitons of the type listed in Proposition \ref{6.3}, with $k = k_0$.
\end{proposition}
\begin{proof}
If the proposition is not true then there is a sequence 
$\{s_j\}_{j=1}^\infty$ tending to infinity
and a neighborhood $U_l$ as in Proposition \ref{6.4} so that for all $j$, 
$g_{s_j}(\cdot) \notin U_l$.
After passing to a further subsequence, we can assume that $\lim_{j \rightarrow
\infty} g_{s_j}(\cdot) \: = \: \overline{g}(\cdot)$ for some
$\overline{g}(\cdot) \in \overline{\mathcal O}$. 

If $k_0 = 2$ then from Proposition \ref{6.2}, $M$ admits a $\Sol$-structure.
As $M$ cannot collapse 
with bounded curvature and bounded diameter to something
of dimension other than one, $\overline{\mathcal O}_{(0)} = 
\overline{\mathcal O}_{(1)} = 
\overline{\mathcal O}_{(3)} = \emptyset$. Then $C = \overline{\mathcal O}_{(2)}$
is sequentially compact and 
$\overline{g}(\cdot) \in \overline{\mathcal O}_{(2)}$. A similar argument applies in the
other cases when $k_0 < 3$ to show that $C = \overline{\mathcal O}_{(k_0)}$ is
sequentially compact and
$\overline{g}(\cdot) \in \overline{\mathcal O}_{(k_0)}$.

For $s \ge 1$, let $\overline{g}^{(s^{-1})}(\cdot)$ be the limit in $\overline{\mathcal O}$
of a convergent subsequence of $\{g_{s^{-1} s_j}(\cdot)\}_{j=1}^\infty$. Then
$\overline{g}(\cdot) \: = \:  \overline{g}^{(s^{-1})}_s(\cdot)$.
Note that $\overline{g}^{(s^{-1})}(\cdot) \in \overline{\mathcal O}_{k_0}$.
By Proposition \ref{6.4}, there is a number
$T_l \ge 1$ so
that for each $s \ge 1$, $\overline{g}^{(s^{-1})}_{T_l}(\cdot) \in U_l$.
Taking $s = T_l$, we conclude that
$\overline{g}(\cdot) \in U_l$.
This is a contradiction.
\end{proof}

\begin{corollary} \label{6.7}
If $k_0 =0$ or $k_0 = 2$ 
then 
$\lim_{s \rightarrow \infty} (M, g_{s}(\cdot))$ exists and is one of
the locally homogeneous expanding 
solitons of the type listed in Proposition \ref{6.3}, with $k = k_0$.
\end{corollary}
\begin{proof}
In these cases, given $M$, there is a unique locally homogeneous expanding
soliton of the
type listed in  Proposition \ref{6.3}, with $k = k_0$.
The relationship between the topology of $M$ and the
equivalence class of the \'etale groupoid comes from the
proof of Lemma \ref{6.1}. 
If $k_0 = 0$ then
$M$ admits a hyperbolic metric and the expander is the
solution $\overline{g}(t) = 4t g_{hyp}$, where $g_{hyp}$ is the metric
of constant sectional curvature $-1$ on $M$. If $k_0 = 2$ then
$M$ is a $\Sol$-manifold. Suppose first that $M$ is the total space
of a $T^2$-bundle over $S^1$, with hyperbolic holonomy $H \in \SL(2, \Z)$.
Then by Remark \ref{4.78}, the expander can be written
$\overline{g} = \frac{t}{2} \Tr(X^2) \: db^2 \: + \: (dy)^T e^{bX} dy$,
where $b \in [0,1]$ and $e^X = H^T H$. If $M$ fibers over the
orbifold $[0,1]$ then the expander is a $\Z_2$-quotient thereof.
\end{proof}

In the case $k_0 = 3$, we must show that any rescaling limit 
$\overline{g}(\cdot)$ has
three local symmetries. This does not follow just from topological
arguments. The method of proof is to rescale backwards
and then apply the monotonicity arguments of Section \ref{section4} to
a backward limit.

\begin{proposition} \label{6.8}
If $k_0 = 3$ and $\overline{g}(\cdot) \in \overline{\mathcal O}$ is a
limit $\lim_{j \rightarrow \infty} (M, g_{s_j}(\cdot))$, for some
sequence $\{s_j\}_{j=1}^\infty$ tending to infinity, then 
$\overline{g}(\cdot) \in \overline{\mathcal O}_{(3)}$.
\end{proposition}
\begin{proof}
Suppose that $\overline{g}(\cdot) \in \overline{\mathcal O}_{(k)}$
with $k < 3$. As in the proof of Proposition \ref{6.6},
for $s \ge 1$, let $\overline{g}^{(s^{-1})}(\cdot)$ be the limit in $\overline{\mathcal O}$
of a convergent subsequence of $\{g_{s^{-1} s_j}(\cdot)\}_{j=1}^\infty$. Then
$\overline{g}(\cdot) \: = \:  \overline{g}^{(s^{-1})}_s(\cdot)$.
More precisely, for each $s \in [1, \infty)$ there is an equivalence
$\phi_s$ of groupoids so that 
\begin{equation} \label{6.9}
\overline{g}(t) \: = \:  \frac{1}{s} \: \phi_s^* \:
\overline{g}^{(s^{-1})}(st).
\end{equation}
In particular, 
$\overline{g}^{(s^{-1})}(\cdot) \in \overline{\mathcal O}_{(k)}$.
Using (\ref{6.9}), we can extend the domain of definition
of $\overline{g}(\cdot)$ to $[s^{-1}, \infty)$ for all $s \ge 1$, and hence
to all $t \in (0, \infty)$. We still have the bounds
$\sup_{t \in (0, \infty)} t \parallel \Riem(\overline{g}(t)) 
\parallel_\infty \: \le \: K$ and 
$\sup_{t \in (0, \infty)} t^{- \: \frac12} \:
\diam(\overline{g}(t)) \: \le \: D$.

As in the proof of Proposition \ref{4.76}, we construct a solution
$f(t)$ of the conjugate heat equation on the orbit space $S$: 
\begin{equation} \label{6.10}
\frac{\partial(e^{-f})}{\partial t} \: = \:
- \: \nabla^2 \: e^{-f} \: + \:
\left( R \: - \: \frac14 \:
g^{\alpha \beta} \: G^{ij} \: G_{jk, \alpha} \: G^{kl} \:
G_{li, \beta} \: - \: \frac12 \: g^{\alpha \gamma} \: g^{\beta \delta} \:
G_{ij} \: F^i_{\alpha \beta} \: F^j_{\gamma \delta}
\: + \: \frac{n}{2t} \right) e^{-f},
\end{equation}
where $n = \dim(S) = 3-k$,
that satisfies $(4 \pi t)^{- \: \frac{n}{2}} \: \int_S e^{-f} \: \dvol_S
\: = \: 1$ for all $t \in (0, \infty)$. Then 
${\mathcal W}_+(G_{ij}(t), A^i_\alpha(t), g_{\alpha \beta}(t),
f(t), t)$ is
nondecreasing in $t$. From Lemma \ref{4.58}, for $t < 1$ 
there is a uniform positive lower bound on $t^{- \: \frac{n}{2}} \: 
\vol(S, g_{\alpha \beta}(t))$. By O'Neill's theorem, the lower sectional
curvature bound on $\overline{g}(t)$ implies the same lower sectional
curvature bound on $g_{\alpha \beta}(t)$. Hence the (orbifolds)
$(S, t^{-1} g_{\alpha \beta}(t))$ are noncollapsing in the
Gromov-Hausdorff sense as $t \rightarrow 0$. It follows that
$\{ \overline{g}^{s^{-1}}(\cdot) \}_{s \ge 1}$ lies in a 
sequentially compact subset of $\overline{\mathcal O}_{(k)}$,
since if a sequence $\{ \overline{g}^{{s_r}^{-1}}(\cdot) \}_{r=1}^\infty$
with $\lim_{r \rightarrow \infty} s_r = \infty$
converged to an element of $\overline{\mathcal O}_{(k^\prime)}$ with
$k^\prime > k$ then the orbit spaces 
$\{(S, s_r \: g_{\alpha \beta}(s_r^{-1}))\}_{j=1}^\infty$ would 
Gromov-Hausdorff converge to something of dimension $3 - k^\prime < 3-k$,
which contradicts the noncollapsing.
In particular, $t^{- \: \frac{n}{2}} \: \vol(S, g_{\alpha \beta}(t))$
is uniformly bounded above as $t \rightarrow 0$ (as also follows
from the diameter and lower curvature bounds).  Then from
Lemma \ref{4.58}, 
${\mathcal W}_+(G_{ij}(t), A^i_\alpha(t), g_{\alpha \beta}(t), f(t), t)$ is
uniformly bounded from below as $t \rightarrow 0$.
There is a sequence of times $t_j \rightarrow 0$ so that
$\lim_{j \rightarrow \infty} t_j \: \frac{d}{dt} \Big|_{t = t_j}
{\mathcal W}_+(G_{ij}(t), A^i_\alpha(t), g_{\alpha \beta}(t), f(t), 
t) \: = \: 0$.
After passing to a subsequence, we can assume that
$\lim_{j \rightarrow \infty} \overline{g}^{t_j}(\cdot) \: = \:
\overline{g}_0(\cdot)$ for some 
$\overline{g}_0(\cdot) \in \overline{\mathcal O}_{(k)}$, defined
for $t \in (0, \infty)$. As in the
proof of Proposition \ref{4.76}, for any $t \in (0, \infty)$ the measures
$(4 \pi t_j t)^{- \:
\frac{n}{2}} \: e^{- f(t_j t)} \: \dvol(S, g_{\alpha \beta}(t_j t))$
will subconverge to a smooth positive probability measure on $S$. 
Using (\ref{4.63}), we get that $\overline{g}_0(\cdot)$ satisfies
the conclusion of Proposition \ref{4.64} at time $t=1$. 
It follows that $\overline{g}_0(\cdot)$ satisfies the conclusion of
Proposition \ref{4.64} for all $t \ge 1$. 
In particular, $M$ admits a geometric structure
other than an $\R^3$ or a $\Nil$-structure (see the proof of
Proposition \ref{6.2}), which is a contradiction.
\end{proof}

\begin{proposition} \label{6.11}
If $k_0 = 3$ then for any sequence $\{s_j\}_{j=1}^\infty$ tending
to infinity,
$\lim_{j \rightarrow \infty} g_{s_j}(\cdot)$ exists and is a locally
homogeneous expanding
soliton of the $\R^3$ or $\Nil$-type.
\end{proposition}
\begin{proof}
If the proposition is not true then there is a sequence
$\{s_j\}_{j=1}^\infty$ tending to infinity 
such that $\lim_{j \rightarrow
\infty} g_{s_j}(\cdot) \: = \: \overline{g}(\cdot)$ for some
$\overline{g}(\cdot) \in \overline{\mathcal O}$, but 
$\overline{g}(\cdot)$ is not an expander of type
$\R^3$ or $\Nil$.
From Proposition \ref{6.8}, 
$\overline{g}(\cdot) \in \overline{\mathcal O}_{(3)}$. In
particular, $\overline{g}(\cdot)$ is locally homogeneous.
If $\underline{\frak g}$ is a local system of $\R^3$ Lie algebras
then $\overline{g}(\cdot)$ must be flat.  
If $\underline{\frak g}$ is a local system of $\nil$ Lie algebras
then $\overline{g}(\cdot)$ 
is automatically a locally homogeneous 
expanding soliton, with respect to some origin of time. 
{\it A priori}, the equation for
$\overline{g}(\cdot)$ could differ
from the expanding $\Nil$ soliton in Theorem \ref{1.2} by an additive
change of the time parameter. We can rule this out by using stability
arguments as before, which are simpler in this case because we
are now talking about dynamics on the finite-dimensional space of
locally homogenous $\Nil$-solutions. First, we argue that there
is some sequence $s_j \rightarrow \infty$ so that 
$\lim_{j \rightarrow \infty} g_{s_j}(\cdot)$ is a locally homogeneous
expanding soliton modeled on the $\Nil$ expanding soliton of 
Theorem \ref{1.2}. Then we use the fact that this expanding soliton
is an attractor for the $\R^{\ge 1}$-semigroup action on
the locally homogeneous $\Nil$-solutions
\cite[Section 3.3.3]{Lott (2007)}. Finally, we use a backward
rescaling, as in the proof of Proposition \ref{6.8}, to show that
for any sequence $\{s_j\}_{j=1}^\infty$ tending
to infinity, $\lim_{j \rightarrow \infty} g_{s_j}(\cdot)$ is 
a locally homogeneous
expanding soliton modeled on the $\Nil$ expanding soliton of 
Theorem \ref{1.2}
\end{proof}

\begin{remark} \label{6.12}
Some of the results of this subsection extend to higher dimension.
Suppose that $(M, g(\cdot))$ is a Ricci flow on a closed $n$-dimensional
manifold that exists for $t \in (1, \infty)$, with
sectional curvatures that are uniformly $O(t^{-1})$ and diameter that
grows at most like $O(t^{\frac12})$. (If $n > 3$ then not all compact
Ricci flows satisfy these assumptions, as seen by
the static solution on a Ricci-flat $K3$ surface.) If 
$\{ s_j \}_{j=1}^\infty$ is any sequence tending to infinity then
after passing to a subsequence, there is a limit Ricci flow
$\overline{g}(\cdot)$ on an $n$-dimensional
\'etale groupoid ${\frak G}$. If $M$ is aspherical
then ${\frak G}$ is locally free.

If $n$ is greater than three then the first point
is that the local symmetry sheaf 
$\underline{\frak g}$ may be a sheaf of nonabelian nilpotent Lie
algebras.  (This could also happen in dimension $3$, but then
$\overline{g}(\cdot)$ is locally homogeneous with respect to the
three-dimensional Heisenberg group.)
Thus the analysis of Subsection \ref{subsection4.2} 
would have to be extended to the
case of twisted ${\mathcal G}$-bundles where ${\mathcal G}$ is a
nilpotent Lie group. 

If we do assume that $\underline{\frak g}$ is
abelian then Proposition \ref{4.76} says that any blowdown
limit of $\overline{g}(\cdot)$ satisfies the harmonic-Einstein
equations (\ref{4.65}). 
Proposition \ref{4.77} describes the blowdown limit of
a Ricci flow solution $(M, g(\cdot))$
on an aspherical $4$-manifold, defined for
$t \in (1, \infty)$, with sectional curvatures that are uniformly
$O(t^{-1})$ and diameter which is $O(t^{\frac12})$, provided that
$\underline{\frak g}$ is abelian.
\end{remark}

\subsection{Proof of Theorem \ref{1.2}} \label{subsection6.2}

In this subsection we use the fact that $M$ is aspherical
in order to extend the convergence result of Subsection \ref{subsection6.1} 
from a statement about a limiting Ricci flow on an \'etale
groupoid to 
a statement about a limiting Ricci flow on $\widetilde{M}$.

By Proposition \ref{3.5}, $M$ is irreducible, aspherical and
has a single geometric piece in its geometric decomposition.
We assume first that $M$ does not have Thurston type
$\widetilde{\SL_2(\R)}$.

Suppose that for a sequence
$\{s_j\}_{j=1}^\infty$ tending to infinity, the limit
$\lim_{j \rightarrow \infty} \left( M, g_{s_j}(\cdot) \right)$
exists and equals a Ricci flow $\overline{g}(\cdot)$ on an
\'etale groupoid ${\frak G}$. If $S$ is the orbit space of
$({\frak G}, \overline{g}(1))$ 
then $\lim_{j \rightarrow \infty} \left(M, \frac{g(s_j)}{s_j}
\right) \stackrel{GH}{=} S$. 

From Propositions \ref{6.6} and \ref{6.11}, $\overline{g}(\cdot)$ is
a locally homogeneous expanding soliton of the type listed in Proposition 
\ref{6.3}.
There is an orbifold fiber bundle $M \rightarrow S$. 
Now $S$ is a very good orbifold, i.e. $S$ is the quotient of
a manifold $\widehat{S}$ by a finite group action. Taking the
corresponding finite cover $\widehat{M}$ of $M$,
if we are interested in what happens on the universal cover
$\widetilde{M}$ then we can assume that $S$ is a closed manifold.

Suppose that $M$ is not of $\Nil$-type.
For large $j$, we know that 
$\left(M, \frac{g(s_j)}{s_j}
\right)$ is the total space of
a $T^{k_0}$-bundle over $S$ which defines an $F$-structure,
where ${k_0} = \dim(M) - \dim(S)$. As $M$ is aspherical, the map
$\pi_1(T^{k_0}) \rightarrow \pi_1(M)$ is injective
\cite[Remark 0.9]{Cheeger-Rong (1995)}.

Choose $\delta \in 
\left( 0, \min \left( \frac{\inj(S)}{10}, 
\frac{1}{10\sqrt{K}} \right) \right)$
and take a finite collection  $\{x_i\}$ of points in $S$ with
the property that $\{B(x_i, \delta)\}$ covers $S$. 
For large $j$, let $\{ p_{i,j} \}$ be points in 
$\left(M, \frac{g(s_j)}{s_j}
\right)$ that are the image of $\{x_i\}$ under a Gromov-Hausdorff
approximation.  Then for such $j$, $\{B(p_{i,j}, 5\delta)\}$ covers 
$\left(M, \frac{g(s_j)}{s_j}
\right)$. Each $B(p_{i,j}, \delta)$ is homeomorphic to
$B^{3-{k_0}} \times T^{k_0}$ and its lift $\widetilde{B(p_{i,j}, \delta)}$
to $\widetilde{M}$ is homeomorphic to $B^{3-{k_0}} \times \R^{k_0}$.

Suppose that $\widetilde{p}_{i,j} \in \widetilde{M}$ is a preimage of
$p_{i,j}$. Then the $5\delta$-ball $B(0, 5\delta) \subset T_{p_{i,j}} M$, 
with respect to the metric $\exp_{p_{i,j}}^* \frac{g(s_j)}{s_j}$, is 
isometric to $B(\widetilde{p}_{i,j}, 5\delta) \subset 
\left( \widetilde{M}, \frac{\widetilde{g}(s_j)}{s_j}
\right)$.
From the construction of the Riemannian groupoid 
$({\frak G}, \overline{g}(1))$
\cite[Proposition 5.9]{Lott (2007)}, $\lim_{j \rightarrow \infty} 
B(\widetilde{p}_{i,j}, 5\delta)$ is isometric to
a $5\delta$-ball in the time-$1$ slice of the
(homogeneous) expanding soliton solution on the manifold $\R^3$. 

Let $\widetilde{m} \in \widetilde{M}$ be a basepoint.  Given $R > 0$,
consider $B(\widetilde{m}, R) \subset  
\left(\widetilde{M}, \frac{\widetilde{g}(s_j)}{s_j}
\right)$. For large $j$, we can find a finite collection of points
$\{ \widetilde{p}_{r} \}$ (depending on $j$) in 
$\left(\widetilde{M}, \frac{\widetilde{g}(s_j)}{s_j}
\right)$, where each $\widetilde{p}_{r}$ projects to some element of
$\{ p_{i,j} \} \subset M$, 
so that the cardinality of $\{ \widetilde{p}_{r} \}$
is uniformly bounded in $j$ and 
$B(\widetilde{m}, R) \subset  
\left(\widetilde{M}, \frac{\widetilde{g}(s_j)}{s_j}
\right)$ is covered by
$\{ B(\widetilde{p}_{r}, 5\delta) \}$. Namely,
for each $i$ and $j$, take points in the strip
$\widetilde{B(p_{i,j}, \delta)} \subset
\left( \widetilde{M}, \frac{\widetilde{g}(s_j)}{s_j} \right)$ that 
lie in $B(\widetilde{m}, R)$, cover $p_{i,j}$ and form a separated net of
size approximately $\delta$.

After relabeling the
indices if necessary, suppose that
$\widetilde{m} \in B(\widetilde{p}_{1}, 5\delta)$ with
$\widetilde{p}_1 \in \widetilde{M}$ 
projectioning to $p_{1,j} \in M$. For large
$j$, fix an almost-isometry from $B(\widetilde{p}_{1}, 5\delta)$
to a $5\delta$-ball in the time-$1$ slice of the
(homogeneous) expanding soliton solution on $\R^3$. 
Taking the union of the balls $B(\widetilde{p}_r, 5 \delta)$ whose
centers project to $p_{1,j} \in M$,
it follows that for large $j$, the metric 
$\frac{\widetilde{g}(s_j)}{s_j}$ on 
$\widetilde{B(p_{1,j}, \delta)} \cap B(\widetilde{m}, R)$
approaches the homogeneous expanding soliton metric on the strip.
We do the same procedure
for the other values of $i$, on the strips 
$\widetilde{B(p_{i,j}, \delta)} \cap B(\widetilde{m}, R) \subset
\left(\widetilde{M}, \frac{\widetilde{g}(s_j)}{s_j}
\right)$.
Then taking the union of these strips for the various $i$,
it follows that for large $j$, the metric 
$\frac{\widetilde{g}(s_j)}{s_j}$ on 
$B(\widetilde{m}, R)$
approaches an $R$-ball in the time-one slice of the
homogeneous expanding soliton solution
on $\R^3$. Finally, we can perform the argument with the time
parameter added, to conclude $\left\{
\left( \widetilde{M}, \widetilde{m},
\widetilde{g}_{s_j}(\cdot) \right) \right\}_{j=1}^\infty$
converges to
the expanding soliton solution on $\R^3$, in the topology of
pointed smooth convergence. We can perform a similar argument
in the $\Nil$-case, where $S$ is a point.

To prove Theorem \ref{1.2}, suppose first that $M$ has Thurston type
$\R^3$, $\Nil$, $\Sol$ or $H^3$. Suppose that the theorem is not true.
Then there is a sequence $\{s_j\}_{j=1}^\infty$ tending to infinity
so that for any subsequence
$\{s_{j_r}\}_{r=1}^\infty$, either
$\left\{ \left( M, g_{s_{j_r}}(1) \right) \right\}_{r=1}^\infty$
does not converge in the Gromov-Hausdorff topology to the limit stated in the
theorem or 
$\left\{ \left( \widetilde{M}, \widetilde{m}, \widetilde{g}_{s_{j_r}}(\cdot) 
\right) \right\}_{r=1}^\infty$ does not converge 
in the pointed smooth topology to the
homogeneous expanding soliton solution stated in the theorem.
After passing to a subsequence,
there is a limit $\lim_{j \rightarrow \infty} 
\left( M, g_{s_j}(\cdot) \right)$ as a Ricci flow on an
\'etale groupoid, whose time-one orbit space will be the
Gromov-Hausdorff limit
$\lim_{j \rightarrow \infty} 
\left( M, \frac{g(s_j)}{s_j} \right)$. 
The limit is characterized by
Corollary \ref{6.7} and Proposition \ref{6.11}.
From the preceding discussion, 
there is a pointed smooth limit $\lim_{j \rightarrow \infty} 
\left( \widetilde{M}, \widetilde{m}, \widetilde{g}_{s_j}(\cdot) \right)$ as a
Ricci flow on $\widetilde{M}$, which is a homogeneous expanding soliton
solution on $\R^3$ of the corresponding type.  In any case, we get
a contradiction.

Suppose now that $M$ has Thurston type $H^2 \times \R$. We can
apply the same argument.  The only difference is that we can no
longer say that $\lim_{s \rightarrow \infty} (M, g_s(1))$ exists
in the Gromov-Hausdorff topology.  All that we get from Proposition
\ref{6.6} is that for any sequence $\{ t_j \}_{j=1}^\infty$ tending
to infinity, there is a subsequence $\{ t_{j_r} \}_{r=1}^\infty$
for which $\lim_{r \rightarrow \infty} 
\left( M, \frac{g(t_{j_r})}{t_{j_r}} \right)$ exists and equals
a closed $2$-dimensional orbifold with constant sectional curvature
$- \: \frac12$. {\em A priori}, different subsequences could give rise to
to different constant-curvature orbifolds. (From our
diameter bound, for a given $M$ 
there is a compact set of such orbifolds that can
arise). However, we claim that on the universal cover
$\widetilde{M}$ we do get pointed smooth convergence 
$\lim_{s \rightarrow \infty} \left( \widetilde{M}, \widetilde{m},
\widetilde{g}_s(\cdot) \right)$ to
the Ricci flow $(H^2 \times \R, 2tg_{hyp} \: + \: g_{\R})$. To see this,
suppose that $\{s_j\}_{j=1}^\infty$ is a sequence tending to infinity
so that for any subsequence
$\{s_{j_r}\}_{r=1}^\infty$, 
$\left\{ \left( \widetilde{M}, \widetilde{g}_{s_{j_r}}(\cdot) 
\right) \right\}_{r=1}^\infty$ does not converge to 
$(H^2 \times \R, 2tg_{hyp} \: + \: g_{\R})$ in the pointed smooth topology.
We know that there is a subsequence $\{s_{j_r}\}_{r=1}^\infty$
for which 
$\lim_{r \rightarrow \infty} 
\left( {M}, {g}_{s_{j_r}}(\cdot) 
\right)$ exists and is a Ricci flow $\overline{g}(\cdot)$ 
on an \'etale groupoid.
From Proposition \ref{6.6}, $\overline{g}(\cdot)$ is a Ricci flow
solution of type $H^2 \times \R$. From the preceding discussion,
$\lim_{r \rightarrow \infty} 
\left( \widetilde{M}, \widetilde{m}, \widetilde{g}_{s_{j_r}}(\cdot) 
\right)$ exists in the pointed smooth topology and 
equals the expanding soliton solution
$(H^2 \times \R, 2tg_{hyp} \: + \: g_{\R})$. This is a contradiction.

Finally, if $M$ has Thurston type
$\widetilde{\SL_2(\R)}$ then the theorem follows from Proposition \ref{6.2}.

\begin{remark} \label{6.13}
The assumption $\diam(M, g(t)) = O(t^{\frac12})$ of Theorem \ref{1.2}, 
together with the curvature assumption, ensures that $M$ has a
single geometric piece.  One can ask what happens if one
removes the diameter assumption but keeps
the curvature assumption.
In such a case one would clearly have to consider pointed limits
$\lim_{j \rightarrow \infty} \left( M, m, g_{s_j}(\cdot) \right)$
of the Ricci flow solution.  After passing to a subsequence,
there will be convergence to a Ricci flow solution
$\overline{g}(\cdot)$ on a pointed
\'etale groupoid. However, the analysis of Subsection \ref{4.2}
does not immediately extend to the pointed noncompact setting.
For example, $\overline{g}(t)$ need not have finite
volume in any reasonable sense; see Example \ref{2.3}.
\end{remark}

\end{document}